
\catcode`\@=11


\message{Loading jyTeX fonts...}



\font\vptrm=cmr5 \font\vptmit=cmmi5 \font\vptsy=cmsy5 \font\vptbf=cmbx5

\skewchar\vptmit='177 \skewchar\vptsy='60 \fontdimen16
\vptsy=\the\fontdimen17 \vptsy

\def\vpt{\ifmmode\err@badsizechange\else
     \@mathfontinit
     \textfont0=\vptrm  \scriptfont0=\vptrm  \scriptscriptfont0=\vptrm
     \textfont1=\vptmit \scriptfont1=\vptmit \scriptscriptfont1=\vptmit
     \textfont2=\vptsy  \scriptfont2=\vptsy  \scriptscriptfont2=\vptsy
     \textfont3=\xptex  \scriptfont3=\xptex  \scriptscriptfont3=\xptex
     \textfont\bffam=\vptbf
     \scriptfont\bffam=\vptbf
     \scriptscriptfont\bffam=\vptbf
     \@fontstyleinit
     \def\rm{\vptrm\fam=\z@}%
     \def\bf{\vptbf\fam=\bffam}%
     \def\oldstyle{\vptmit\fam=\@ne}%
     \rm\fi}


\font\viptrm=cmr6 \font\viptmit=cmmi6 \font\viptsy=cmsy6
\font\viptbf=cmbx6

\skewchar\viptmit='177 \skewchar\viptsy='60 \fontdimen16
\viptsy=\the\fontdimen17 \viptsy

\def\vipt{\ifmmode\err@badsizechange\else
     \@mathfontinit
     \textfont0=\viptrm  \scriptfont0=\vptrm  \scriptscriptfont0=\vptrm
     \textfont1=\viptmit \scriptfont1=\vptmit \scriptscriptfont1=\vptmit
     \textfont2=\viptsy  \scriptfont2=\vptsy  \scriptscriptfont2=\vptsy
     \textfont3=\xptex   \scriptfont3=\xptex  \scriptscriptfont3=\xptex
     \textfont\bffam=\viptbf
     \scriptfont\bffam=\vptbf
     \scriptscriptfont\bffam=\vptbf
     \@fontstyleinit
     \def\rm{\viptrm\fam=\z@}%
     \def\bf{\viptbf\fam=\bffam}%
     \def\oldstyle{\viptmit\fam=\@ne}%
     \rm\fi}

\font\viiptrm=cmr7 \font\viiptmit=cmmi7 \font\viiptsy=cmsy7
\font\viiptit=cmti7 \font\viiptbf=cmbx7

\skewchar\viiptmit='177 \skewchar\viiptsy='60 \fontdimen16
\viiptsy=\the\fontdimen17 \viiptsy

\def\viipt{\ifmmode\err@badsizechange\else
     \@mathfontinit
     \textfont0=\viiptrm  \scriptfont0=\vptrm  \scriptscriptfont0=\vptrm
     \textfont1=\viiptmit \scriptfont1=\vptmit \scriptscriptfont1=\vptmit
     \textfont2=\viiptsy  \scriptfont2=\vptsy  \scriptscriptfont2=\vptsy
     \textfont3=\xptex    \scriptfont3=\xptex  \scriptscriptfont3=\xptex
     \textfont\itfam=\viiptit
     \scriptfont\itfam=\viiptit
     \scriptscriptfont\itfam=\viiptit
     \textfont\bffam=\viiptbf
     \scriptfont\bffam=\vptbf
     \scriptscriptfont\bffam=\vptbf
     \@fontstyleinit
     \def\rm{\viiptrm\fam=\z@}%
     \def\it{\viiptit\fam=\itfam}%
     \def\bf{\viiptbf\fam=\bffam}%
     \def\oldstyle{\viiptmit\fam=\@ne}%
     \rm\fi}


\font\viiiptrm=cmr8 \font\viiiptmit=cmmi8 \font\viiiptsy=cmsy8
\font\viiiptit=cmti8
\font\viiiptbf=cmbx8

\skewchar\viiiptmit='177 \skewchar\viiiptsy='60 \fontdimen16
\viiiptsy=\the\fontdimen17 \viiiptsy

\def\viiipt{\ifmmode\err@badsizechange\else
     \@mathfontinit
     \textfont0=\viiiptrm  \scriptfont0=\viptrm  \scriptscriptfont0=\vptrm
     \textfont1=\viiiptmit \scriptfont1=\viptmit \scriptscriptfont1=\vptmit
     \textfont2=\viiiptsy  \scriptfont2=\viptsy  \scriptscriptfont2=\vptsy
     \textfont3=\xptex     \scriptfont3=\xptex   \scriptscriptfont3=\xptex
     \textfont\itfam=\viiiptit
     \scriptfont\itfam=\viiptit
     \scriptscriptfont\itfam=\viiptit
     \textfont\bffam=\viiiptbf
     \scriptfont\bffam=\viptbf
     \scriptscriptfont\bffam=\vptbf
     \@fontstyleinit
     \def\rm{\viiiptrm\fam=\z@}%
     \def\it{\viiiptit\fam=\itfam}%
     \def\bf{\viiiptbf\fam=\bffam}%
     \def\oldstyle{\viiiptmit\fam=\@ne}%
     \rm\fi}


\def\getixpt{%
     \font\ixptrm=cmr9
     \font\ixptmit=cmmi9
     \font\ixptsy=cmsy9
     \font\ixptit=cmti9
     \font\ixptbf=cmbx9
     \skewchar\ixptmit='177 \skewchar\ixptsy='60
     \fontdimen16 \ixptsy=\the\fontdimen17 \ixptsy}

\def\ixpt{\ifmmode\err@badsizechange\else
     \@mathfontinit
     \textfont0=\ixptrm  \scriptfont0=\viiptrm  \scriptscriptfont0=\vptrm
     \textfont1=\ixptmit \scriptfont1=\viiptmit \scriptscriptfont1=\vptmit
     \textfont2=\ixptsy  \scriptfont2=\viiptsy  \scriptscriptfont2=\vptsy
     \textfont3=\xptex   \scriptfont3=\xptex    \scriptscriptfont3=\xptex
     \textfont\itfam=\ixptit
     \scriptfont\itfam=\viiptit
     \scriptscriptfont\itfam=\viiptit
     \textfont\bffam=\ixptbf
     \scriptfont\bffam=\viiptbf
     \scriptscriptfont\bffam=\vptbf
     \@fontstyleinit
     \def\rm{\ixptrm\fam=\z@}%
     \def\it{\ixptit\fam=\itfam}%
     \def\bf{\ixptbf\fam=\bffam}%
     \def\oldstyle{\ixptmit\fam=\@ne}%
     \rm\fi}


\font\xptrm=cmr10 \font\xptmit=cmmi10 \font\xptsy=cmsy10
\font\xptex=cmex10 \font\xptit=cmti10 \font\xptsl=cmsl10
\font\xptbf=cmbx10 \font\xpttt=cmtt10 \font\xptss=cmss10
\font\xptsc=cmcsc10 \font\xptbfs=cmb10 \font\xptbmit=cmmib10

\skewchar\xptmit='177 \skewchar\xptbmit='177 \skewchar\xptsy='60
\fontdimen16 \xptsy=\the\fontdimen17 \xptsy

\def\xpt{\ifmmode\err@badsizechange\else
     \@mathfontinit
     \textfont0=\xptrm  \scriptfont0=\viiptrm  \scriptscriptfont0=\vptrm
     \textfont1=\xptmit \scriptfont1=\viiptmit \scriptscriptfont1=\vptmit
     \textfont2=\xptsy  \scriptfont2=\viiptsy  \scriptscriptfont2=\vptsy
     \textfont3=\xptex  \scriptfont3=\xptex    \scriptscriptfont3=\xptex
     \textfont\itfam=\xptit
     \scriptfont\itfam=\viiptit
     \scriptscriptfont\itfam=\viiptit
     \textfont\bffam=\xptbf
     \scriptfont\bffam=\viiptbf
     \scriptscriptfont\bffam=\vptbf
     \textfont\bfsfam=\xptbfs
     \scriptfont\bfsfam=\viiptbf
     \scriptscriptfont\bfsfam=\vptbf
     \textfont\bmitfam=\xptbmit
     \scriptfont\bmitfam=\viiptmit
     \scriptscriptfont\bmitfam=\vptmit
     \@fontstyleinit
     \def\rm{\xptrm\fam=\z@}%
     \def\it{\xptit\fam=\itfam}%
     \def\sl{\xptsl}%
     \def\bf{\xptbf\fam=\bffam}%
     \def\tt{\xpttt}%
     \def\ss{\xptss}%
     \def\sc{\xptsc}%
     \def\bfs{\xptbfs\fam=\bfsfam}%
     \def\bmit{\fam=\bmitfam}%
     \def\oldstyle{\xptmit\fam=\@ne}%
     \rm\fi}


\def\getxipt{%
     \font\xiptrm=cmr10  scaled\magstephalf
     \font\xiptmit=cmmi10 scaled\magstephalf
     \font\xiptsy=cmsy10 scaled\magstephalf
     \font\xiptex=cmex10 scaled\magstephalf
     \font\xiptit=cmti10 scaled\magstephalf
     \font\xiptsl=cmsl10 scaled\magstephalf
     \font\xiptbf=cmbx10 scaled\magstephalf
     \font\xipttt=cmtt10 scaled\magstephalf
     \font\xiptss=cmss10 scaled\magstephalf
     \skewchar\xiptmit='177 \skewchar\xiptsy='60
     \fontdimen16 \xiptsy=\the\fontdimen17 \xiptsy}

\def\xipt{\ifmmode\err@badsizechange\else
     \@mathfontinit
     \textfont0=\xiptrm  \scriptfont0=\viiiptrm  \scriptscriptfont0=\viptrm
     \textfont1=\xiptmit \scriptfont1=\viiiptmit \scriptscriptfont1=\viptmit
     \textfont2=\xiptsy  \scriptfont2=\viiiptsy  \scriptscriptfont2=\viptsy
     \textfont3=\xiptex  \scriptfont3=\xptex     \scriptscriptfont3=\xptex
     \textfont\itfam=\xiptit
     \scriptfont\itfam=\viiiptit
     \scriptscriptfont\itfam=\viiptit
     \textfont\bffam=\xiptbf
     \scriptfont\bffam=\viiiptbf
     \scriptscriptfont\bffam=\viptbf
     \@fontstyleinit
     \def\rm{\xiptrm\fam=\z@}%
     \def\it{\xiptit\fam=\itfam}%
     \def\sl{\xiptsl}%
     \def\bf{\xiptbf\fam=\bffam}%
     \def\tt{\xipttt}%
     \def\ss{\xiptss}%
     \def\oldstyle{\xiptmit\fam=\@ne}%
     \rm\fi}


\font\xiiptrm=cmr12 \font\xiiptmit=cmmi12 \font\xiiptsy=cmsy10
scaled\magstep1 \font\xiiptex=cmex10  scaled\magstep1 \font\xiiptit=cmti12
\font\xiiptsl=cmsl12 \font\xiiptbf=cmbx12
\font\xiiptss=cmss12 \font\xiiptsc=cmcsc10 scaled\magstep1
\font\xiiptbfs=cmb10  scaled\magstep1 \font\xiiptbmit=cmmib10
scaled\magstep1

\skewchar\xiiptmit='177 \skewchar\xiiptbmit='177 \skewchar\xiiptsy='60
\fontdimen16 \xiiptsy=\the\fontdimen17 \xiiptsy

\def\xiipt{\ifmmode\err@badsizechange\else
     \@mathfontinit
     \textfont0=\xiiptrm  \scriptfont0=\viiiptrm  \scriptscriptfont0=\viptrm
     \textfont1=\xiiptmit \scriptfont1=\viiiptmit \scriptscriptfont1=\viptmit
     \textfont2=\xiiptsy  \scriptfont2=\viiiptsy  \scriptscriptfont2=\viptsy
     \textfont3=\xiiptex  \scriptfont3=\xptex     \scriptscriptfont3=\xptex
     \textfont\itfam=\xiiptit
     \scriptfont\itfam=\viiiptit
     \scriptscriptfont\itfam=\viiptit
     \textfont\bffam=\xiiptbf
     \scriptfont\bffam=\viiiptbf
     \scriptscriptfont\bffam=\viptbf
     \textfont\bfsfam=\xiiptbfs
     \scriptfont\bfsfam=\viiiptbf
     \scriptscriptfont\bfsfam=\viptbf
     \textfont\bmitfam=\xiiptbmit
     \scriptfont\bmitfam=\viiiptmit
     \scriptscriptfont\bmitfam=\viptmit
     \@fontstyleinit
     \def\rm{\xiiptrm\fam=\z@}%
     \def\it{\xiiptit\fam=\itfam}%
     \def\sl{\xiiptsl}%
     \def\bf{\xiiptbf\fam=\bffam}%
     \def\tt{\xiipttt}%
     \def\ss{\xiiptss}%
     \def\sc{\xiiptsc}%
     \def\bfs{\xiiptbfs\fam=\bfsfam}%
     \def\bmit{\fam=\bmitfam}%
     \def\oldstyle{\xiiptmit\fam=\@ne}%
     \rm\fi}


\def\getxiiipt{%
     \font\xiiiptrm=cmr12  scaled\magstephalf
     \font\xiiiptmit=cmmi12 scaled\magstephalf
     \font\xiiiptsy=cmsy9  scaled\magstep2
     \font\xiiiptit=cmti12 scaled\magstephalf
     \font\xiiiptsl=cmsl12 scaled\magstephalf
     \font\xiiiptbf=cmbx12 scaled\magstephalf
     \font\xiiipttt=cmtt12 scaled\magstephalf
     \font\xiiiptss=cmss12 scaled\magstephalf
     \skewchar\xiiiptmit='177 \skewchar\xiiiptsy='60
     \fontdimen16 \xiiiptsy=\the\fontdimen17 \xiiiptsy}

\def\xiiipt{\ifmmode\err@badsizechange\else
     \@mathfontinit
     \textfont0=\xiiiptrm  \scriptfont0=\xptrm  \scriptscriptfont0=\viiptrm
     \textfont1=\xiiiptmit \scriptfont1=\xptmit \scriptscriptfont1=\viiptmit
     \textfont2=\xiiiptsy  \scriptfont2=\xptsy  \scriptscriptfont2=\viiptsy
     \textfont3=\xivptex   \scriptfont3=\xptex  \scriptscriptfont3=\xptex
     \textfont\itfam=\xiiiptit
     \scriptfont\itfam=\xptit
     \scriptscriptfont\itfam=\viiptit
     \textfont\bffam=\xiiiptbf
     \scriptfont\bffam=\xptbf
     \scriptscriptfont\bffam=\viiptbf
     \@fontstyleinit
     \def\rm{\xiiiptrm\fam=\z@}%
     \def\it{\xiiiptit\fam=\itfam}%
     \def\sl{\xiiiptsl}%
     \def\bf{\xiiiptbf\fam=\bffam}%
     \def\tt{\xiiipttt}%
     \def\ss{\xiiiptss}%
     \def\oldstyle{\xiiiptmit\fam=\@ne}%
     \rm\fi}


\font\xivptrm=cmr12   scaled\magstep1 \font\xivptmit=cmmi12
scaled\magstep1 \font\xivptsy=cmsy10  scaled\magstep2
\font\xivptex=cmex10  scaled\magstep2 \font\xivptit=cmti12 scaled\magstep1
\font\xivptsl=cmsl12  scaled\magstep1 \font\xivptbf=cmbx12
scaled\magstep1
\font\xivptss=cmss12  scaled\magstep1 \font\xivptsc=cmcsc10
scaled\magstep2 \font\xivptbfs=cmb10  scaled\magstep2
\font\xivptbmit=cmmib10 scaled\magstep2

\skewchar\xivptmit='177 \skewchar\xivptbmit='177 \skewchar\xivptsy='60
\fontdimen16 \xivptsy=\the\fontdimen17 \xivptsy

\def\xivpt{\ifmmode\err@badsizechange\else
     \@mathfontinit
     \textfont0=\xivptrm  \scriptfont0=\xptrm  \scriptscriptfont0=\viiptrm
     \textfont1=\xivptmit \scriptfont1=\xptmit \scriptscriptfont1=\viiptmit
     \textfont2=\xivptsy  \scriptfont2=\xptsy  \scriptscriptfont2=\viiptsy
     \textfont3=\xivptex  \scriptfont3=\xptex  \scriptscriptfont3=\xptex
     \textfont\itfam=\xivptit
     \scriptfont\itfam=\xptit
     \scriptscriptfont\itfam=\viiptit
     \textfont\bffam=\xivptbf
     \scriptfont\bffam=\xptbf
     \scriptscriptfont\bffam=\viiptbf
     \textfont\bfsfam=\xivptbfs
     \scriptfont\bfsfam=\xptbfs
     \scriptscriptfont\bfsfam=\viiptbf
     \textfont\bmitfam=\xivptbmit
     \scriptfont\bmitfam=\xptbmit
     \scriptscriptfont\bmitfam=\viiptmit
     \@fontstyleinit
     \def\rm{\xivptrm\fam=\z@}%
     \def\it{\xivptit\fam=\itfam}%
     \def\sl{\xivptsl}%
     \def\bf{\xivptbf\fam=\bffam}%
     \def\tt{\xivpttt}%
     \def\ss{\xivptss}%
     \def\sc{\xivptsc}%
     \def\bfs{\xivptbfs\fam=\bfsfam}%
     \def\bmit{\fam=\bmitfam}%
     \def\oldstyle{\xivptmit\fam=\@ne}%
     \rm\fi}


\font\xviiptrm=cmr17 \font\xviiptmit=cmmi12 scaled\magstep2
\font\xviiptsy=cmsy10 scaled\magstep3 \font\xviiptex=cmex10
scaled\magstep3 \font\xviiptit=cmti12 scaled\magstep2
\font\xviiptbf=cmbx12 scaled\magstep2 \font\xviiptbfs=cmb10
scaled\magstep3

\skewchar\xviiptmit='177 \skewchar\xviiptsy='60 \fontdimen16
\xviiptsy=\the\fontdimen17 \xviiptsy

\def\xviipt{\ifmmode\err@badsizechange\else
     \@mathfontinit
     \textfont0=\xviiptrm  \scriptfont0=\xiiptrm  \scriptscriptfont0=\viiiptrm
     \textfont1=\xviiptmit \scriptfont1=\xiiptmit \scriptscriptfont1=\viiiptmit
     \textfont2=\xviiptsy  \scriptfont2=\xiiptsy  \scriptscriptfont2=\viiiptsy
     \textfont3=\xviiptex  \scriptfont3=\xiiptex  \scriptscriptfont3=\xptex
     \textfont\itfam=\xviiptit
     \scriptfont\itfam=\xiiptit
     \scriptscriptfont\itfam=\viiiptit
     \textfont\bffam=\xviiptbf
     \scriptfont\bffam=\xiiptbf
     \scriptscriptfont\bffam=\viiiptbf
     \textfont\bfsfam=\xviiptbfs
     \scriptfont\bfsfam=\xiiptbfs
     \scriptscriptfont\bfsfam=\viiiptbf
     \@fontstyleinit
     \def\rm{\xviiptrm\fam=\z@}%
     \def\it{\xviiptit\fam=\itfam}%
     \def\bf{\xviiptbf\fam=\bffam}%
     \def\bfs{\xviiptbfs\fam=\bfsfam}%
     \def\oldstyle{\xviiptmit\fam=\@ne}%
     \rm\fi}


\font\xxiptrm=cmr17  scaled\magstep1


\def\xxipt{\ifmmode\err@badsizechange\else
     \@mathfontinit
     \@fontstyleinit
     \def\rm{\xxiptrm\fam=\z@}%
     \rm\fi}


\font\xxvptrm=cmr17  scaled\magstep2


\def\xxvpt{\ifmmode\err@badsizechange\else
     \@mathfontinit
     \@fontstyleinit
     \def\rm{\xxvptrm\fam=\z@}%
     \rm\fi}




\message{Loading jyTeX macros...}

\message{modifications to plain.tex,}


\def\newcount{\alloc@0\count\countdef\insc@unt}
\def\newdimen{\alloc@1\dimen\dimendef\insc@unt}
\def\newskip{\alloc@2\skip\skipdef\insc@unt}
\def\newmuskip{\alloc@3\muskip\muskipdef\@cclvi}
\def\newbox{\alloc@4\box\chardef\insc@unt}
\def\newtoks{\alloc@5\toks\toksdef\@cclvi}
\def\newhelp#1#2{\newtoks#1\global#1\expandafter{\csname#2\endcsname}}
\def\newread{\alloc@6\read\chardef\sixt@@n}
\def\newwrite{\alloc@7\write\chardef\sixt@@n}
\def\newfam{\alloc@8\fam\chardef\sixt@@n}
\def\newinsert#1{\global\advance\insc@unt by\m@ne
     \ch@ck0\insc@unt\count
     \ch@ck1\insc@unt\dimen
     \ch@ck2\insc@unt\skip
     \ch@ck4\insc@unt\box
     \allocationnumber=\insc@unt
     \global\chardef#1=\allocationnumber
     \wlog{\string#1=\string\insert\the\allocationnumber}}
\def\newif#1{\count@\escapechar \escapechar\m@ne
     \expandafter\expandafter\expandafter
          \xdef\@if#1{true}{\let\noexpand#1=\noexpand\iftrue}%
     \expandafter\expandafter\expandafter
          \xdef\@if#1{false}{\let\noexpand#1=\noexpand\iffalse}%
     \global\@if#1{false}\escapechar=\count@}


\newlinechar=`\^^J
\overfullrule=0pt




\let\itfam=\undefined

\let\bffam=\undefined

\count18=3


\chardef\sharps="19


\mathchardef\alpha="710B \mathchardef\beta="710C
\mathchardef\gamma="710D \mathchardef\delta="710E
\mathchardef\epsilon="710F \mathchardef\zeta="7110
\mathchardef\eta="7111 \mathchardef\theta="7112 \mathchardef\iota="7113
\mathchardef\kappa="7114 \mathchardef\lambda="7115
\mathchardef\mu="7116 \mathchardef\nu="7117 \mathchardef\xi="7118
\mathchardef\pi="7119 \mathchardef\rho="711A \mathchardef\sigma="711B
\mathchardef\tau="711C \mathchardef\upsilon="711D
\mathchardef\phi="711E \mathchardef\chi="711F \mathchardef\psi="7120
\mathchardef\omega="7121 \mathchardef\varepsilon="7122
\mathchardef\vartheta="7123 \mathchardef\varpi="7124
\mathchardef\varrho="7125 \mathchardef\varsigma="7126
\mathchardef\varphi="7127 \mathchardef\imath="717B
\mathchardef\jmath="717C \mathchardef\ell="7160 \mathchardef\wp="717D
\mathchardef\partial="7140 \mathchardef\flat="715B
\mathchardef\natural="715C \mathchardef\sharp="715D



\def\angle{{\vbox{\ialign{$\m@th\scriptstyle##$\crcr
     \not\mathrel{\mkern14mu}\crcr
     \noalign{\nointerlineskip}
     \mkern2.5mu\leaders\hrule height.34\rp@\hfill\mkern2.5mu\crcr}}}}
\def\vdots{\vbox{\baselineskip4\rp@ \lineskiplimit\z@
     \kern6\rp@\hbox{.}\hbox{.}\hbox{.}}}
\def\ddots{\mathinner{\mkern1mu\raise7\rp@\vbox{\kern7\rp@\hbox{.}}\mkern2mu
     \raise4\rp@\hbox{.}\mkern2mu\raise\rp@\hbox{.}\mkern1mu}}
\def\overrightarrow#1{\vbox{\ialign{##\crcr
     \rightarrowfill\crcr
     \noalign{\kern-\rp@\nointerlineskip}
     $\hfil\displaystyle{#1}\hfil$\crcr}}}
\def\overleftarrow#1{\vbox{\ialign{##\crcr
     \leftarrowfill\crcr
     \noalign{\kern-\rp@\nointerlineskip}
     $\hfil\displaystyle{#1}\hfil$\crcr}}}
\def\overbrace#1{\mathop{\vbox{\ialign{##\crcr
     \noalign{\kern3\rp@}
     \downbracefill\crcr
     \noalign{\kern3\rp@\nointerlineskip}
     $\hfil\displaystyle{#1}\hfil$\crcr}}}\limits}
\def\underbrace#1{\mathop{\vtop{\ialign{##\crcr
     $\hfil\displaystyle{#1}\hfil$\crcr
     \noalign{\kern3\rp@\nointerlineskip}
     \upbracefill\crcr
     \noalign{\kern3\rp@}}}}\limits}
\def\big#1{{\hbox{$\left#1\vbox to8.5\rp@ {}\right.\n@space$}}}
\def\Big#1{{\hbox{$\left#1\vbox to11.5\rp@ {}\right.\n@space$}}}
\def\bigg#1{{\hbox{$\left#1\vbox to14.5\rp@ {}\right.\n@space$}}}
\def\Bigg#1{{\hbox{$\left#1\vbox to17.5\rp@ {}\right.\n@space$}}}
\def\@vereq#1#2{\lower.5\rp@\vbox{\baselineskip\z@skip\lineskip-.5\rp@
     \ialign{$\m@th#1\hfil##\hfil$\crcr#2\crcr=\crcr}}}
\def\rlh@#1{\vcenter{\hbox{\ooalign{\raise2\rp@
     \hbox{$#1\rightharpoonup$}\crcr
     $#1\leftharpoondown$}}}}
\def\bordermatrix#1{\begingroup\m@th
     \setbox\z@\vbox{%
          \def\cr{\crcr\noalign{\kern2\rp@\global\let\cr\endline}}%
          \ialign{$##$\hfil\kern2\rp@\kern\p@renwd
               &\thinspace\hfil$##$\hfil&&\quad\hfil$##$\hfil\crcr
               \omit\strut\hfil\crcr
               \noalign{\kern-\baselineskip}%
               #1\crcr\omit\strut\cr}}%
     \setbox\tw@\vbox{\unvcopy\z@\global\setbox\@ne\lastbox}%
     \setbox\tw@\hbox{\unhbox\@ne\unskip\global\setbox\@ne\lastbox}%
     \setbox\tw@\hbox{$\kern\wd\@ne\kern-\p@renwd\left(\kern-\wd\@ne
          \global\setbox\@ne\vbox{\box\@ne\kern2\rp@}%
          \vcenter{\kern-\ht\@ne\unvbox\z@\kern-\baselineskip}%
          \,\right)$}%
     \null\;\vbox{\kern\ht\@ne\box\tw@}\endgroup}
\def\endinsert{\egroup
     \if@mid\dimen@\ht\z@
          \advance\dimen@\dp\z@
          \advance\dimen@12\rp@
          \advance\dimen@\pagetotal
          \ifdim\dimen@>\pagegoal\@midfalse\p@gefalse\fi
     \fi
     \if@mid\bigskip\box\z@
          \bigbreak
     \else\insert\topins{\penalty100 \splittopskip\z@skip
               \splitmaxdepth\maxdimen\floatingpenalty\z@
               \ifp@ge\dimen@\dp\z@
                    \vbox to\vsize{\unvbox\z@\kern-\dimen@}%
               \else\box\z@\nobreak\bigskip
               \fi}%
     \fi
     \endgroup}


\def\cases#1{\left\{\,\vcenter{\m@th
     \ialign{$##\hfil$&\quad##\hfil\crcr#1\crcr}}\right.}
\def\matrix#1{\null\,\vcenter{\m@th
     \ialign{\hfil$##$\hfil&&\quad\hfil$##$\hfil\crcr
          \mathstrut\crcr
          \noalign{\kern-\baselineskip}
          #1\crcr
          \mathstrut\crcr
          \noalign{\kern-\baselineskip}}}\,}


\newif\ifraggedbottom

\def\raggedbottom{\ifraggedbottom\else
     \advance\topskip by\z@ plus60pt \raggedbottomtrue\fi}%
\def\normalbottom{\ifraggedbottom
     \advance\topskip by\z@ plus-60pt \raggedbottomfalse\fi}

\message{hacks,}


\toksdef\toks@i=1 \toksdef\toks@ii=2


\def\TeX{T\kern-.1667em \lower.5ex \hbox{E}\kern-.125em X\null}
\def\jyTeX{{\leavevmode
     \raise.587ex \hbox{\it\j}\kern-.1em \lower.048ex \hbox{\it y}\kern-.12em
     \TeX}}

\let\then=\iftrue
\def\ifnoarg#1\then{\def\hack@{#1}\ifx\hack@\empty}
\def\ifundefined#1\then{%
     \expandafter\ifx\csname\expandafter\blank\string#1\endcsname\relax}
\def\useif#1\then{\csname#1\endcsname}
\def\usename#1{\csname#1\endcsname}
\def\useafter#1#2{\expandafter#1\csname#2\endcsname}

\long\def\loop#1\repeat{\def\@iterate{#1\expandafter\@iterate\fi}\@iterate
     \let\@iterate=\relax}

\let\TeXend=\end
\def\begin#1{\begingroup\def\@@blockname{#1}\usename{begin#1}}
\def\end#1{\usename{end#1}\def\hack@{#1}%
     \ifx\@@blockname\hack@
          \endgroup
     \else\err@badgroup\hack@\@@blockname
     \fi}
\def\@@blockname{}

\def\defaultoption[#1]#2{%
     \def\hack@{\ifx\hack@ii[\toks@={#2}\else\toks@={#2[#1]}\fi\the\toks@}%
     \futurelet\hack@ii\hack@}

\def\markup#1{\let\@@marksf=\empty
     \ifhmode\edef\@@marksf{\spacefactor=\the\spacefactor\relax}\/\fi
     ${}^{\hbox{\subscriptfonts#1}}$\@@marksf}


\newtoks\shortyear
\newtoks\militaryhour
\newtoks\standardhour
\newtoks\minute
\newtoks\amorpm

\def\settime{\count@=\time\divide\count@ by60
     \militaryhour=\expandafter{\number\count@}%
     {\multiply\count@ by-60 \advance\count@ by\time
          \xdef\hack@{\ifnum\count@<10 0\fi\number\count@}}%
     \minute=\expandafter{\hack@}%
     \ifnum\count@<12
          \amorpm={am}
     \else\amorpm={pm}
          \ifnum\count@>12 \advance\count@ by-12 \fi
     \fi
     \standardhour=\expandafter{\number\count@}%
     \def\hack@19##1##2{\shortyear={##1##2}}%
          \expandafter\hack@\the\year}

\def\monthword#1{%
     \ifcase#1
          $\bullet$\err@badcountervalue{monthword}%
          \or January\or February\or March\or April\or May\or June%
          \or July\or August\or September\or October\or November\or December%
     \else$\bullet$\err@badcountervalue{monthword}%
     \fi}

\def\monthabbr#1{%
     \ifcase#1
          $\bullet$\err@badcountervalue{monthabbr}%
          \or Jan\or Feb\or Mar\or Apr\or May\or Jun%
          \or Jul\or Aug\or Sep\or Oct\or Nov\or Dec%
     \else$\bullet$\err@badcountervalue{monthabbr}%
     \fi}

\def\militarytime{\the\militaryhour:\the\minute}
\def\standardtime{\the\standardhour:\the\minute}


\def\@setnumstyle#1#2{\expandafter\global\expandafter\expandafter
     \expandafter\let\expandafter\expandafter
     \csname @\expandafter\blank\string#1style\endcsname
     \csname#2\endcsname}
\def\numstyle#1{\usename{@\expandafter\blank\string#1style}#1}
\def\ifblank#1\then{\useafter\ifx{@\expandafter\blank\string#1}\blank}

\def\blank#1{}

\def\Roman#1{\expandafter\uppercase\expandafter{\romannumeral#1}}
\def\alphabetic#1{%
     \ifcase#1
          $\bullet$\err@badcountervalue{alphabetic}%
          \or a\or b\or c\or d\or e\or f\or g\or h\or i\or j\or k\or l\or m%
          \or n\or o\or p\or q\or r\or s\or t\or u\or v\or w\or x\or y\or z%
     \else$\bullet$\err@badcountervalue{alphabetic}%
     \fi}
\def\Alphabetic#1{\expandafter\uppercase\expandafter{\alphabetic{#1}}}
\def\symbols#1{%
     \ifcase#1
          $\bullet$\err@badcountervalue{symbols}%
          \or*\or\dag\or\ddag\or\S\or$\|$%
          \or**\or\dag\dag\or\ddag\ddag\or\S\S\or$\|\|$%
     \else$\bullet$\err@badcountervalue{symbols}%
     \fi}


\catcode`\^^?=13 \def^^?{\relax}

\def\trimleading#1\to#2{\edef#2{#1}%
     \expandafter\@trimleading\expandafter#2#2^^?^^?}
\def\@trimleading#1#2#3^^?{\ifx#2^^?\def#1{}\else\def#1{#2#3}\fi}

\def\trimtrailing#1\to#2{\edef#2{#1}%
     \expandafter\@trimtrailing\expandafter#2#2^^? ^^?\relax}
\def\@trimtrailing#1#2 ^^?#3{\ifx#3\relax\toks@={}%
     \else\def#1{#2}\toks@={\trimtrailing#1\to#1}\fi
     \the\toks@}

\def\trim#1\to#2{\trimleading#1\to#2\trimtrailing#2\to#2}

\catcode`\^^?=15


\long\def\additemL#1\to#2{\toks@={\^^\{#1}}\toks@ii=\expandafter{#2}%
     \xdef#2{\the\toks@\the\toks@ii}}

\long\def\additemR#1\to#2{\toks@={\^^\{#1}}\toks@ii=\expandafter{#2}%
     \xdef#2{\the\toks@ii\the\toks@}}

\def\getitemL#1\to#2{\expandafter\@getitemL#1\hack@#1#2}
\def\@getitemL\^^\#1#2\hack@#3#4{\def#4{#1}\def#3{#2}}

\message{font macros,}


\newdimen\rp@
\newcount\@@sizeindex \@@sizeindex=0
\newcount\@@factori
\newcount\@@factorii
\newcount\@@factoriii
\newcount\@@factoriv

\countdef\maxfam=18
\newfam\itfam
\newfam\bffam
\newfam\bfsfam
\newfam\bmitfam

\def\@mathfontinit{\count@=4
     \loop\textfont\count@=\nullfont
          \scriptfont\count@=\nullfont
          \scriptscriptfont\count@=\nullfont
          \ifnum\count@<\maxfam\advance\count@ by\@ne
     \repeat}

\def\@fontstyleinit{%
     \def\it{\err@fontnotavailable\it}%
     \def\bf{\err@fontnotavailable\bf}%
     \def\bfs{\err@bfstobf}%
     \def\bmit{\err@fontnotavailable\bmit}%
     \def\sc{\err@fontnotavailable\sc}%
     \def\sl{\err@sltoit}%
     \def\ss{\err@fontnotavailable\ss}%
     \def\tt{\err@fontnotavailable\tt}}

\def\@parameterinit#1{\rm\rp@=.1em \@getscaling{#1}%
     \let\^^\=\@doscaling\scalingskipslist
     \setbox\strutbox=\hbox{\vrule
          height.708\baselineskip depth.292\baselineskip width\z@}}

\def\@getfactor#1#2#3#4{\@@factori=#1 \@@factorii=#2
     \@@factoriii=#3 \@@factoriv=#4}

\def\@getscaling#1{\count@=#1 \advance\count@ by-\@@sizeindex\@@sizeindex=#1
     \ifnum\count@<0
          \let\@mulordiv=\divide
          \let\@divormul=\multiply
          \multiply\count@ by\m@ne
     \else\let\@mulordiv=\multiply
          \let\@divormul=\divide
     \fi
     \edef\@@scratcha{\ifcase\count@                {1}{1}{1}{1}\or
          {1}{7}{23}{3}\or     {2}{5}{3}{1}\or      {9}{89}{13}{1}\or
          {6}{25}{6}{1}\or     {8}{71}{14}{1}\or    {6}{25}{36}{5}\or
          {1}{7}{53}{4}\or     {12}{125}{108}{5}\or {3}{14}{53}{5}\or
          {6}{41}{17}{1}\or    {13}{31}{13}{2}\or   {9}{107}{71}{2}\or
          {11}{139}{124}{3}\or {1}{6}{43}{2}\or     {10}{107}{42}{1}\or
          {1}{5}{43}{2}\or     {5}{69}{65}{1}\or    {11}{97}{91}{2}\fi}%
     \expandafter\@getfactor\@@scratcha}

\def\@doscaling#1{\@mulordiv#1by\@@factori\@divormul#1by\@@factorii
     \@mulordiv#1by\@@factoriii\@divormul#1by\@@factoriv}


\newskip\headskip
\newskip\footskip

\def\typesize=#1pt{\count@=#1 \advance\count@ by-10
     \ifcase\count@
          \@setsizex\or\err@badtypesize\or
          \@setsizexii\or\err@badtypesize\or
          \@setsizexiv
     \else\err@badtypesize
     \fi}

\def\@setsizex{\getixpt
     \def\subsubscriptfonts{\vpt}%
          \def\subsubscriptsize{\vpt\@parameterinit{-8}}%
     \def\subscriptfonts{\viipt}\def\subscriptsize{\viipt\@parameterinit{-4}}%
     \def\footnotefonts{\viiipt}\def\footnotesize{\viiipt\@parameterinit{-2}}%
     \def\smallfonts{\ixpt}\def\smallsize{\ixpt\@parameterinit{-1}}%
     \def\normalfonts{\xpt}\def\normalsize{\xpt\@parameterinit{0}}%
     \def\bigfonts{\xiipt}\def\bigsize{\xiipt\@parameterinit{2}}%
     \def\Bigfonts{\xivpt}\def\Bigsize{\xivpt\@parameterinit{4}}%
     \def\biggfonts{\xviipt}\def\biggsize{\xviipt\@parameterinit{6}}%
     \def\Biggfonts{\xxipt}\def\Biggsize{\xxipt\@parameterinit{8}}%
     \def\tinyfonts{\vpt}\def\tinysize{\vpt\@parameterinit{-8}}%
     \def\HUGEFONTS{\xxvpt}\def\HUGESIZE{\xxvpt\@parameterinit{10}}%
     \normalsize\fixedskipslist}

\def\@setsizexii{\getxipt
     \def\subsubscriptfonts{\vipt}%
          \def\subsubscriptsize{\vipt\@parameterinit{-6}}%
     \def\subscriptfonts{\viiipt}%
          \def\subscriptsize{\viiipt\@parameterinit{-2}}%
     \def\footnotefonts{\xpt}\def\footnotesize{\xpt\@parameterinit{0}}%
     \def\smallfonts{\xipt}\def\smallsize{\xipt\@parameterinit{1}}%
     \def\normalfonts{\xiipt}\def\normalsize{\xiipt\@parameterinit{2}}%
     \def\bigfonts{\xivpt}\def\bigsize{\xivpt\@parameterinit{4}}%
     \def\Bigfonts{\xviipt}\def\Bigsize{\xviipt\@parameterinit{6}}%
     \def\biggfonts{\xxipt}\def\biggsize{\xxipt\@parameterinit{8}}%
     \def\Biggfonts{\xxvpt}\def\Biggsize{\xxvpt\@parameterinit{10}}%
     \def\tinyfonts{\vpt}\def\tinysize{\vpt\@parameterinit{-8}}%
     \def\HUGEFONTS{\xxvpt}\def\HUGESIZE{\xxvpt\@parameterinit{10}}%
     \normalsize\fixedskipslist}

\def\@setsizexiv{\getxiiipt
     \def\subsubscriptfonts{\viipt}%
          \def\subsubscriptsize{\viipt\@parameterinit{-4}}%
     \def\subscriptfonts{\xpt}\def\subscriptsize{\xpt\@parameterinit{0}}%
     \def\footnotefonts{\xiipt}\def\footnotesize{\xiipt\@parameterinit{2}}%
     \def\smallfonts{\xiiipt}\def\smallsize{\xiiipt\@parameterinit{3}}%
     \def\normalfonts{\xivpt}\def\normalsize{\xivpt\@parameterinit{4}}%
     \def\bigfonts{\xviipt}\def\bigsize{\xviipt\@parameterinit{6}}%
     \def\Bigfonts{\xxipt}\def\Bigsize{\xxipt\@parameterinit{8}}%
     \def\biggfonts{\xxvpt}\def\biggsize{\xxvpt\@parameterinit{10}}%
     \def\Biggfonts{\err@sizetoolarge\Biggfonts\HUGEFONTS}%
          \def\Biggsize{\err@sizetoolarge\Biggsize\HUGESIZE}%
     \def\tinyfonts{\vpt}\def\tinysize{\vpt\@parameterinit{-8}}%
     \def\HUGEFONTS{\xxvpt}\def\HUGESIZE{\xxvpt\@parameterinit{10}}%
     \normalsize\fixedskipslist}

\def\subsubscriptfonts{\vpt} \def\subsubscriptsize{\vpt\@parameterinit{-8}}
\def\subscriptfonts{\viipt}  \def\subscriptsize{\viipt\@parameterinit{-4}}
\def\footnotefonts{\viiipt}  \def\footnotesize{\viiipt\@parameterinit{-2}}
\def\smallfonts{\err@sizenotavailable\smallfonts}
                             \def\smallsize{\ixpt\@parameterinit{-1}}
\def\normalfonts{\xpt}       \def\normalsize{\xpt\@parameterinit{0}}
\def\bigfonts{\xiipt}        \def\bigsize{\xiipt\@parameterinit{2}}
\def\Bigfonts{\xivpt}        \def\Bigsize{\xivpt\@parameterinit{4}}
\def\biggfonts{\xviipt}      \def\biggsize{\xviipt\@parameterinit{6}}
\def\Biggfonts{\xxipt}       \def\Biggsize{\xxipt\@parameterinit{8}}
\def\tinyfonts{\vpt}         \def\tinysize{\vpt\@parameterinit{-8}}
\def\HUGEFONTS{\xxvpt}       \def\HUGESIZE{\xxvpt\@parameterinit{10}}

\message{document layout,}


\newtoks\everyoutput \everyoutput={}
\newdimen\depthofpage
\newcount\pagenum \pagenum=0

\newdimen\oddtopmargin  \newdimen\eventopmargin
\newdimen\oddleftmargin \newdimen\evenleftmargin
\newtoks\oddhead        \newtoks\evenhead
\newtoks\oddfoot        \newtoks\evenfoot

\def\topmargin{\afterassignment\@seteventop\oddtopmargin}
\def\leftmargin{\afterassignment\@setevenleft\oddleftmargin}
\def\head{\afterassignment\@setevenhead\oddhead}
\def\foot{\afterassignment\@setevenfoot\oddfoot}

\def\@seteventop{\eventopmargin=\oddtopmargin}
\def\@setevenleft{\evenleftmargin=\oddleftmargin}
\def\@setevenhead{\evenhead=\oddhead}
\def\@setevenfoot{\evenfoot=\oddfoot}

\def\pagenumstyle#1{\@setnumstyle\pagenum{#1}}

\newif\ifdraft
\def\draft{\drafttrue\leftmargin=.5in \overfullrule=5pt }

\def\outputstyle#1{\global\expandafter\let\expandafter
          \@outputstyle\csname#1output\endcsname
     \usename{#1setup}}

\output={\@outputstyle}

\def\normaloutput{\the\everyoutput
     \global\advance\pagenum by\@ne
     \ifodd\pagenum
          \voffset=\oddtopmargin \hoffset=\oddleftmargin
     \else\voffset=\eventopmargin \hoffset=\evenleftmargin
     \fi
     \advance\voffset by-1in  \advance\hoffset by-1in
     \count0=\pagenum
     \expandafter\shipout\pagebox
     \ifnum\outputpenalty>-\@MM\else\dosupereject\fi}

\newdimen\fullhsize
\newbox\leftpage
\newcount\leftpagenum
\newcount\outputpagenum \outputpagenum=0
\let\leftorright=L

\def\twoupoutput{\the\everyoutput
     \global\advance\pagenum by\@ne
     \if L\leftorright
          \global\setbox\leftpage=\leftline{\pagebox}%
          \global\leftpagenum=\pagenum
          \global\let\leftorright=R%
     \else\global\advance\outputpagenum by\@ne
          \ifodd\outputpagenum
               \voffset=\oddtopmargin \hoffset=\oddleftmargin
          \else\voffset=\eventopmargin \hoffset=\evenleftmargin
          \fi
          \advance\voffset by-1in  \advance\hoffset by-1in
          \count0=\leftpagenum \count1=\pagenum
          \shipout\vbox{\hbox to\fullhsize
               {\box\leftpage\hfil\leftline{\pagebox}}}%
          \global\let\leftorright=L%
     \fi
     \ifnum\outputpenalty>-\@MM
     \else\dosupereject
          \if R\leftorright
               \globaldefs=\@ne\head={\hfil}\foot={\hfil}\globaldefs=\z@
               \null\newpage
          \fi
     \fi}

\def\pagebox{\vbox{\makeheadline\pagebody\makefootline}}

\def\makeheadline{%
     \vbox to\z@{\baselinestretch=\@m
          \vskip\topskip\vskip-.708\baselineskip\vskip-\headskip
          \line{\vbox to\ht\strutbox{}%
               \ifodd\pagenum\the\oddhead\else\the\evenhead\fi}%
          \vss}%
     \nointerlineskip}

\def\pagebody{\vbox to\vsize{%
     \boxmaxdepth\maxdepth
     \ifvoid\topins\else\unvbox\topins\fi
     \depthofpage=\dp255
     \unvbox255
     \ifraggedbottom\kern-\depthofpage\vfil\fi
     \ifvoid\footins
     \else\vskip\skip\footins
          \footnoterule
          \unvbox\footins
          \vskip-\footnoteskip
     \fi}}

\def\makefootline{\baselineskip=\footskip
     \line{\ifodd\pagenum\the\oddfoot\else\the\evenfoot\fi}}


\newskip\abovechapterskip
\newskip\belowchapterskip
\newskip\abovesectionskip
\newskip\belowsectionskip
\newskip\abovesubsectionskip
\newskip\belowsubsectionskip

\def\chapterstyle#1{\global\expandafter\let\expandafter\@chapterstyle
     \csname#1text\endcsname}
\def\sectionstyle#1{\global\expandafter\let\expandafter\@sectionstyle
     \csname#1text\endcsname}
\def\subsectionstyle#1{\global\expandafter\let\expandafter\@subsectionstyle
     \csname#1text\endcsname}

\def\chapter#1{%
     \ifdim\lastskip=17sp \else\chapterbreak\vskip\abovechapterskip\fi
     \@chapterstyle{\ifblank\chapternumstyle\then
          \else\newchapternum=\next\chapternumformat\ \fi#1}%
     \nobreak\vskip\belowchapterskip\vskip17sp }

\def\section#1{%
     \ifdim\lastskip=17sp \else\sectionbreak\vskip\abovesectionskip\fi
     \@sectionstyle{\ifblank\sectionnumstyle\then
          \else\newsectionnum=\next\sectionnumformat\ \fi#1}%
     \nobreak\vskip\belowsectionskip\vskip17sp }

\def\subsection#1{%
     \ifdim\lastskip=17sp \else\subsectionbreak\vskip\abovesubsectionskip\fi
     \@subsectionstyle{\ifblank\subsectionnumstyle\then
          \else\newsubsectionnum=\next\subsectionnumformat\ \fi#1}%
     \nobreak\vskip\belowsubsectionskip\vskip17sp }


\let\TeXunderline=\underline
\let\TeXoverline=\overline
\def\underline#1{\relax\ifmmode\TeXunderline{#1}\else
     $\TeXunderline{\hbox{#1}}$\fi}
\def\overline#1{\relax\ifmmode\TeXoverline{#1}\else
     $\TeXoverline{\hbox{#1}}$\fi}

\def\baselinestretch{\afterassignment\@baselinestretch\count@}
\def\@baselinestretch{\baselineskip=\normalbaselineskip
     \divide\baselineskip by\@m\baselineskip=\count@\baselineskip
     \setbox\strutbox=\hbox{\vrule
          height.708\baselineskip depth.292\baselineskip width\z@}%
     \bigskipamount=\the\baselineskip
          plus.25\baselineskip minus.25\baselineskip
     \medskipamount=.5\baselineskip
          plus.125\baselineskip minus.125\baselineskip
     \smallskipamount=.25\baselineskip
          plus.0625\baselineskip minus.0625\baselineskip}

\def\\{\ifhmode\ifnum\lastpenalty=-\@M\else\hfil\penalty-\@M\fi\fi
     \ignorespaces}
\def\newpage{\vfil\break}

\def\lefttext#1{\par{\@text\leftskip=\z@\rightskip=\centering
     \noindent#1\par}}
\def\righttext#1{\par{\@text\leftskip=\centering\rightskip=\z@
     \noindent#1\par}}
\def\centertext#1{\par{\@text\leftskip=\centering\rightskip=\centering
     \noindent#1\par}}
\def\@text{\parindent=\z@ \parfillskip=\z@ \everypar={}%
     \spaceskip=.3333em \xspaceskip=.5em
     \def\\{\ifhmode\ifnum\lastpenalty=-\@M\else\penalty-\@M\fi\fi
          \ignorespaces}}

\def\beginleft{\par\@text\leftskip=\z@ \rightskip=\centering}
     
\def\beginright{\par\@text\leftskip=\centering\rightskip=\z@ }
     
\def\begincenter{\par\@text\leftskip=\centering\rightskip=\centering}

\def\beginnarrow{\defaultoption[\parindent]\@beginnarrow}
\def\@beginnarrow[#1]{\par\advance\leftskip by#1\advance\rightskip by#1}

\begingroup
\catcode`\[=1 \catcode`\{=11 \gdef\beginignore[\endgroup\bgroup
     \catcode`\e=0 \catcode`\\=12 \catcode`\{=11 \catcode`\f=12 \let\or=\relax
     \let\nd{ignor=\fi \let\}=\egroup
     \iffalse}
\endgroup

\long\def\marginnote#1{\leavevmode
     \edef\@marginsf{\spacefactor=\the\spacefactor\relax}%
     \ifdraft\strut\vadjust{%
          \hbox to\z@{\hskip\hsize\hskip.1in
               \vbox to\z@{\vskip-\dp\strutbox
                    \marginnoteformat
                    \vskip-\ht\strutbox
                    \noindent\strut#1\par
                    \vss}%
               \hss}}%
     \fi
     \@marginsf}


\newtoks\everybye \everybye={\par\vfil}
\outer\def\bye{\the\everybye
     \footnotecheck
     \prelabelcheck
     \streamcheck
     \supereject
     \TeXend}

\message{footnotes,}

\newcount\footnotenum \footnotenum=0
\newskip\footnoteskip
\let\@footnotelist=\empty

\def\footnotenumstyle#1{\@setnumstyle\footnotenum{#1}%
     \useafter\ifx{@footnotenumstyle}\symbols
          \global\let\@footup=\empty
     \else\global\let\@footup=\markup
     \fi}

\def\footnote{\footnotecheck\defaultoption[]\@footnote}
\def\@footnote[#1]{\@footnotemark[#1]\@footnotetext}

\def\footnotemark{\defaultoption[]\@footnotemark}
\def\@footnotemark[#1]{\let\@footsf=\empty
     \ifhmode\edef\@footsf{\spacefactor=\the\spacefactor\relax}\/\fi
     \ifnoarg#1\then
          \global\advance\footnotenum by\@ne
          \@footup{\footnotenumformat}%
          \edef\@@foota{\footnotenum=\the\footnotenum\relax}%
          \expandafter\additemR\expandafter\@footup\expandafter
               {\@@foota\footnotenumformat}\to\@footnotelist
          \global\let\@footnotelist=\@footnotelist
     \else\markup{#1}%
          \additemR\markup{#1}\to\@footnotelist
          \global\let\@footnotelist=\@footnotelist
     \fi
     \@footsf}

\def\footnotetext{%
     \ifx\@footnotelist\empty\err@extrafootnotetext\else\@footnotetext\fi}
\def\@footnotetext{%
     \getitemL\@footnotelist\to\@@foota
     \global\let\@footnotelist=\@footnotelist
     \insert\footins\bgroup
     \footnoteformat
     \splittopskip=\ht\strutbox\splitmaxdepth=\dp\strutbox
     \interlinepenalty=\interfootnotelinepenalty\floatingpenalty=\@MM
     \noindent\llap{\@@foota}\strut
     \bgroup\aftergroup\@footnoteend
     \let\@@scratcha=}
\def\@footnoteend{\strut\par\vskip\footnoteskip\egroup}

\def\footnoterule{\normalfonts
     \kern-.3em \hrule width2in height.04em \kern .26em }

\def\footnotecheck{%
     \ifx\@footnotelist\empty
     \else\err@extrafootnotemark
          \global\let\@footnotelist=\empty
     \fi}

\message{labels,}

\let\@@labeldef=\xdef
\newif\if@labelfile
\newwrite\@labelfile
\let\@prelabellist=\empty

\def\label#1#2{\trim#1\to\@@labarg\edef\@@labtext{#2}%
     \edef\@@labname{lab@\@@labarg}%
     \useafter\ifundefined\@@labname\then\else\@yeslab\fi
     \useafter\@@labeldef\@@labname{#2}%
     \ifstreaming
          \expandafter\toks@\expandafter\expandafter\expandafter
               {\csname\@@labname\endcsname}%
          \immediate\write\streamout{\noexpand\label{\@@labarg}{\the\toks@}}%
     \fi}
\def\@yeslab{%
     \useafter\ifundefined{if\@@labname}\then
          \err@labelredef\@@labarg
     \else\useif{if\@@labname}\then
               \err@labelredef\@@labarg
          \else\global\usename{\@@labname true}%
               \useafter\ifundefined{pre\@@labname}\then
               \else\useafter\ifx{pre\@@labname}\@@labtext
                    \else\err@badlabelmatch\@@labarg
                    \fi
               \fi
               \if@labelfile
               \else\global\@labelfiletrue
                    \immediate\write\sixt@@n{--> Creating file \jobname.lab}%
                    \immediate\openout\@labelfile=\jobname.lab
               \fi
               \immediate\write\@labelfile
                    {\noexpand\prelabel{\@@labarg}{\@@labtext}}%
          \fi
     \fi}

\def\putlab#1{\trim#1\to\@@labarg\edef\@@labname{lab@\@@labarg}%
     \useafter\ifundefined\@@labname\then\@nolab\else\usename\@@labname\fi}
\def\@nolab{%
     \useafter\ifundefined{pre\@@labname}\then
          \undefinedlabelformat
          \err@needlabel\@@labarg
          \useafter\xdef\@@labname{\undefinedlabelformat}%
     \else\usename{pre\@@labname}%
          \useafter\xdef\@@labname{\usename{pre\@@labname}}%
     \fi
     \useafter\newif{if\@@labname}%
     \expandafter\additemR\@@labarg\to\@prelabellist}

\def\prelabel#1{\useafter\gdef{prelab@#1}}

\def\ifundefinedlabel#1\then{%
     \expandafter\ifx\csname lab@#1\endcsname\relax}
\def\useiflab#1\then{\csname iflab@#1\endcsname}

\def\prelabelcheck{{%
     \def\^^\##1{\useiflab{##1}\then\else\err@undefinedlabel{##1}\fi}%
     \@prelabellist}}

\message{equation numbering,}

\newcount\chapternum
\newcount\sectionnum
\newcount\subsectionnum
\newcount\equationnum
\newcount\subequationnum
\newcount\figurenum
\newcount\subfigurenum
\newcount\tablenum
\newcount\subtablenum

\newif\if@subeqncount
\newif\if@subfigcount
\newif\if@subtblcount

\def\newchapternum{\newsectionnum=\z@\@resetnum\chapternum}
\def\newsectionnum{\newsubsectionnum=\z@\@resetnum\sectionnum}
\def\newsubsectionnum{\newequationnum=\z@\newfigurenum=\z@\newtablenum=\z@
     \@resetnum\subsectionnum}
\def\newequationnum{\newsubequationnum=\z@\@resetnum\equationnum}
\def\newsubequationnum{\@resetnum\subequationnum}
\def\newfigurenum{\newsubfigurenum=\z@\@resetnum\figurenum}
\def\newsubfigurenum{\@resetnum\subfigurenum}
\def\newtablenum{\newsubtablenum=\z@\@resetnum\tablenum}
\def\newsubtablenum{\@resetnum\subtablenum}

\def\@resetnum#1{\global\advance#1by1 \edef\next{\the#1\relax}\global#1}

\newchapternum=0

\def\chapternumstyle#1{\@setnumstyle\chapternum{#1}}
\def\sectionnumstyle#1{\@setnumstyle\sectionnum{#1}}
\def\subsectionnumstyle#1{\@setnumstyle\subsectionnum{#1}}
\def\equationnumstyle#1{\@setnumstyle\equationnum{#1}}
\def\subequationnumstyle#1{\@setnumstyle\subequationnum{#1}%
     \ifblank\subequationnumstyle\then\global\@subeqncountfalse\fi
     \ignorespaces}
\def\figurenumstyle#1{\@setnumstyle\figurenum{#1}}
\def\subfigurenumstyle#1{\@setnumstyle\subfigurenum{#1}%
     \ifblank\subfigurenumstyle\then\global\@subfigcountfalse\fi
     \ignorespaces}
\def\tablenumstyle#1{\@setnumstyle\tablenum{#1}}
\def\subtablenumstyle#1{\@setnumstyle\subtablenum{#1}%
     \ifblank\subtablenumstyle\then\global\@subtblcountfalse\fi
     \ignorespaces}

\def\eqnlabel#1{%
     \if@subeqncount
          \newsubequationnum=\next
     \else\newequationnum=\next
          \ifblank\subequationnumstyle\then
          \else\global\@subeqncounttrue
               \newsubequationnum=\@ne
          \fi
     \fi
     \label{#1}{\puteqnformat}(\puteqn{#1})%
     \ifdraft\rlap{\hskip.1in{\tt#1}}\fi}

\let\puteqn=\putlab

\def\equation#1#2{\useafter\gdef{eqn@#1}{#2\eqno\eqnlabel{#1}}}
\def\Equation#1{\useafter\gdef{eqn@#1}}

\def\putequation#1{\useafter\ifundefined{eqn@#1}\then
     \err@undefinedeqn{#1}\else\usename{eqn@#1}\fi}

\def\eqnseriesstyle#1{\gdef\@eqnseriesstyle{#1}}
\def\begineqnseries{\subequationnumstyle{\@eqnseriesstyle}%
     \defaultoption[]\@begineqnseries}
\def\@begineqnseries[#1]{\edef\@@eqnname{#1}}
\def\endeqnseries{\subequationnumstyle{blank}%
     \expandafter\ifnoarg\@@eqnname\then
     \else\label\@@eqnname{\puteqnformat}%
     \fi
     \aftergroup\ignorespaces}

\def\figlabel#1{%
     \if@subfigcount
          \newsubfigurenum=\next
     \else\newfigurenum=\next
          \ifblank\subfigurenumstyle\then
          \else\global\@subfigcounttrue
               \newsubfigurenum=\@ne
          \fi
     \fi
     \label{#1}{\putfigformat}\putfig{#1}%
     {\def\marginnoteformat{\tt}\marginnote{#1}}}

\let\putfig=\putlab

\def\figseriesstyle#1{\gdef\@figseriesstyle{#1}}
\def\beginfigseries{\subfigurenumstyle{\@figseriesstyle}%
     \defaultoption[]\@beginfigseries}
\def\@beginfigseries[#1]{\edef\@@figname{#1}}
\def\endfigseries{\subfigurenumstyle{blank}%
     \expandafter\ifnoarg\@@figname\then
     \else\label\@@figname{\putfigformat}%
     \fi
     \aftergroup\ignorespaces}

\def\tbllabel#1{%
     \if@subtblcount
          \newsubtablenum=\next
     \else\newtablenum=\next
          \ifblank\subtablenumstyle\then
          \else\global\@subtblcounttrue
               \newsubtablenum=\@ne
          \fi
     \fi
     \label{#1}{\puttblformat}\puttbl{#1}%
     {\def\marginnoteformat{\tt}\marginnote{#1}}}

\let\puttbl=\putlab

\def\tblseriesstyle#1{\gdef\@tblseriesstyle{#1}}
\def\begintblseries{\subtablenumstyle{\@tblseriesstyle}%
     \defaultoption[]\@begintblseries}
\def\@begintblseries[#1]{\edef\@@tblname{#1}}
\def\endtblseries{\subtablenumstyle{blank}%
     \expandafter\ifnoarg\@@tblname\then
     \else\label\@@tblname{\puttblformat}%
     \fi
     \aftergroup\ignorespaces}

\message{reference numbering,}

\newcount\referencenum \referencenum=0
\newcount\@@prerefcount \@@prerefcount=0
\newcount\@@thisref
\newcount\@@lastref
\newcount\@@loopref
\newcount\@@refseq
\newdimen\refnumindent
\let\@undefreflist=\empty

\def\referencenumstyle#1{\@setnumstyle\referencenum{#1}}

\def\referencestyle#1{\usename{@ref#1}}

\def\@refsequential{%
     \gdef\@refpredef##1{\global\advance\referencenum by\@ne
          \let\^^\=0\label{##1}{\^^\{\the\referencenum}}%
          \useafter\gdef{ref@\the\referencenum}{{##1}{\undefinedlabelformat}}}%
     \gdef\@reference##1##2{%
          \ifundefinedlabel##1\then
          \else\def\^^\####1{\global\@@thisref=####1\relax}\putlab{##1}%
               \useafter\gdef{ref@\the\@@thisref}{{##1}{##2}}%
          \fi}%
     \gdef\endputreferences{%
          \loop\ifnum\@@loopref<\referencenum
                    \advance\@@loopref by\@ne
                    \expandafter\expandafter\expandafter\@printreference
                         \csname ref@\the\@@loopref\endcsname
          \repeat
          \par}}

\def\@refpreordered{%
     \gdef\@refpredef##1{\global\advance\referencenum by\@ne
          \additemR##1\to\@undefreflist}%
     \gdef\@reference##1##2{%
          \ifundefinedlabel##1\then
          \else\global\advance\@@loopref by\@ne
               {\let\^^\=0\label{##1}{\^^\{\the\@@loopref}}}%
               \@printreference{##1}{##2}%
          \fi}
     \gdef\endputreferences{%
          \def\^^\####1{\useiflab{####1}\then
               \else\reference{####1}{\undefinedlabelformat}\fi}%
          \@undefreflist
          \par}}

\def\beginprereferences{\par
     \def\reference##1##2{\global\advance\referencenum by1\@ne
          \let\^^\=0\label{##1}{\^^\{\the\referencenum}}%
          \useafter\gdef{ref@\the\referencenum}{{##1}{##2}}}}
\def\endprereferences{\global\@@prerefcount=\the\referencenum\par}

\def\beginputreferences{\par
     \refnumindent=\z@\@@loopref=\z@
     \loop\ifnum\@@loopref<\referencenum
               \advance\@@loopref by\@ne
               \setbox\z@=\hbox{\referencenum=\@@loopref
                    \referencenumformat\enskip}%
               \ifdim\wd\z@>\refnumindent\refnumindent=\wd\z@\fi
     \repeat
     \putreferenceformat
     \@@loopref=\z@
     \loop\ifnum\@@loopref<\@@prerefcount
               \advance\@@loopref by\@ne
               \expandafter\expandafter\expandafter\@printreference
                    \csname ref@\the\@@loopref\endcsname
     \repeat
     \let\reference=\@reference}

\def\@printreference#1#2{\ifx#2\undefinedlabelformat\err@undefinedref{#1}\fi
     \noindent\ifdraft\rlap{\hskip\hsize\hskip.1in \tt#1}\fi
     \llap{\referencenum=\@@loopref\referencenumformat\enskip}#2\par}

\def\reference#1#2{{\par\refnumindent=\z@\putreferenceformat\noindent#2\par}}

\def\putref#1{\trim#1\to\@@refarg
     \expandafter\ifnoarg\@@refarg\then
          \toks@={\relax}%
     \else\@@lastref=-\@m\def\@@refsep{}\def\@more{\@nextref}%
          \toks@={\@nextref#1,,}%
     \fi\the\toks@}
\def\@nextref#1,{\trim#1\to\@@refarg
     \expandafter\ifnoarg\@@refarg\then
          \let\@more=\relax
     \else\ifundefinedlabel\@@refarg\then
               \expandafter\@refpredef\expandafter{\@@refarg}%
          \fi
          \def\^^\##1{\global\@@thisref=##1\relax}%
          \global\@@thisref=\m@ne
          \setbox\z@=\hbox{\putlab\@@refarg}%
     \fi
     \advance\@@lastref by\@ne
     \ifnum\@@lastref=\@@thisref\advance\@@refseq by\@ne\else\@@refseq=\@ne\fi
     \ifnum\@@lastref<\z@
     \else\ifnum\@@refseq<\thr@@
               \@@refsep\def\@@refsep{,}%
               \ifnum\@@lastref>\z@
                    \advance\@@lastref by\m@ne
                    {\referencenum=\@@lastref\putrefformat}%
               \else\undefinedlabelformat
               \fi
          \else\def\@@refsep{--}%
          \fi
     \fi
     \@@lastref=\@@thisref
     \@more}

\message{streaming,}

\newif\ifstreaming

\def\streamto{\defaultoption[\jobname]\@streamto}
\def\@streamto[#1]{\global\streamingtrue
     \immediate\write\sixt@@n{--> Streaming to #1.str}%
     \newwrite\streamout\immediate\openout\streamout=#1.str }

\def\streamfrom{\defaultoption[\jobname]\@streamfrom}
\def\@streamfrom[#1]{\newread\streamin\openin\streamin=#1.str
     \ifeof\streamin
          \expandafter\err@nostream\expandafter{#1.str}%
     \else\immediate\write\sixt@@n{--> Streaming from #1.str}%
          \let\@@labeldef=\gdef
          \ifstreaming
               \edef\@elc{\endlinechar=\the\endlinechar}%
               \endlinechar=\m@ne
               \loop\read\streamin to\@@scratcha
                    \ifeof\streamin
                         \streamingfalse
                    \else\toks@=\expandafter{\@@scratcha}%
                         \immediate\write\streamout{\the\toks@}%
                    \fi
                    \ifstreaming
               \repeat
               \@elc
               \input #1.str
               \streamingtrue
          \else\input #1.str
          \fi
          \let\@@labeldef=\xdef
     \fi}

\def\streamcheck{\ifstreaming
     \immediate\write\streamout{\pagenum=\the\pagenum}%
     \immediate\write\streamout{\footnotenum=\the\footnotenum}%
     \immediate\write\streamout{\referencenum=\the\referencenum}%
     \immediate\write\streamout{\chapternum=\the\chapternum}%
     \immediate\write\streamout{\sectionnum=\the\sectionnum}%
     \immediate\write\streamout{\subsectionnum=\the\subsectionnum}%
     \immediate\write\streamout{\equationnum=\the\equationnum}%
     \immediate\write\streamout{\subequationnum=\the\subequationnum}%
     \immediate\write\streamout{\figurenum=\the\figurenum}%
     \immediate\write\streamout{\subfigurenum=\the\subfigurenum}%
     \immediate\write\streamout{\tablenum=\the\tablenum}%
     \immediate\write\streamout{\subtablenum=\the\subtablenum}%
     \immediate\closeout\streamout
     \fi}


\def\err@badtypesize{%
     \errhelp={The limited availability of certain fonts requires^^J%
          that the base type size be 10pt, 12pt, or 14pt.^^J}%
     \errmessage{--> Illegal base type size}}

\def\err@badsizechange{\immediate\write\sixt@@n
     {--> Size change not allowed in math mode, ignored}}

\def\err@sizetoolarge#1{\immediate\write\sixt@@n
     {--> \noexpand#1 too big, substituting HUGE}}

\def\err@sizenotavailable#1{\immediate\write\sixt@@n
     {--> Size not available, \noexpand#1 ignored}}

\def\err@fontnotavailable#1{\immediate\write\sixt@@n
     {--> Font not available, \noexpand#1 ignored}}

\def\err@sltoit{\immediate\write\sixt@@n
     {--> Style \noexpand\sl not available, substituting \noexpand\it}%
     \it}

\def\err@bfstobf{\immediate\write\sixt@@n
     {--> Style \noexpand\bfs not available, substituting \noexpand\bf}%
     \bf}

\def\err@badgroup#1#2{%
     \errhelp={The block you have just tried to close was not the one^^J%
          most recently opened.^^J}%
     \errmessage{--> \noexpand\end{#1} doesn't match \noexpand\begin{#2}}}

\def\err@badcountervalue#1{\immediate\write\sixt@@n
     {--> Counter (#1) out of bounds}}

\def\err@extrafootnotemark{\immediate\write\sixt@@n
     {--> \noexpand\footnotemark command
          has no corresponding \noexpand\footnotetext}}

\def\err@extrafootnotetext{%
     \errhelp{You have given a \noexpand\footnotetext command without first
          specifying^^Ja \noexpand\footnotemark.^^J}%
     \errmessage{--> \noexpand\footnotetext command has no corresponding
          \noexpand\footnotemark}}

\def\err@labelredef#1{\immediate\write\sixt@@n
     {--> Label "#1" redefined}}

\def\err@badlabelmatch#1{\immediate\write\sixt@@n
     {--> Definition of label "#1" doesn't match value in \jobname.lab}}

\def\err@needlabel#1{\immediate\write\sixt@@n
     {--> Label "#1" cited before its definition}}

\def\err@undefinedlabel#1{\immediate\write\sixt@@n
     {--> Label "#1" cited but never defined}}

\def\err@undefinedeqn#1{\immediate\write\sixt@@n
     {--> Equation "#1" not defined}}

\def\err@undefinedref#1{\immediate\write\sixt@@n
     {--> Reference "#1" not defined}}

\def\err@nostream#1{%
     \errhelp={You have tried to input a stream file that doesn't exist.^^J}%
     \errmessage{--> Stream file #1 not found}}

\message{jyTeX initialization}

\everyjob{\immediate\write16{--> jyTeX version \fmtversion}%
     \edef\@@jobname{\jobname}%
     \edef\jobname{\@@jobname}%
     \settime
     \openin0=\jobname.lab
     \ifeof0
     \else\closein0
          \immediate\write16{--> Getting labels from file \jobname.lab}%
          \input\jobname.lab
     \fi}


\def\fixedskipslist{%
     \^^\{\topskip}%
     \^^\{\splittopskip}%
     \^^\{\maxdepth}%
     \^^\{\skip\topins}%
     \^^\{\skip\footins}%
     \^^\{\headskip}%
     \^^\{\footskip}}

\def\scalingskipslist{%
     \^^\{\p@renwd}%
     \^^\{\delimitershortfall}%
     \^^\{\nulldelimiterspace}%
     \^^\{\scriptspace}%
     \^^\{\jot}%
     \^^\{\normalbaselineskip}%
     \^^\{\normallineskip}%
     \^^\{\normallineskiplimit}%
     \^^\{\baselineskip}%
     \^^\{\lineskip}%
     \^^\{\lineskiplimit}%
     \^^\{\bigskipamount}%
     \^^\{\medskipamount}%
     \^^\{\smallskipamount}%
     \^^\{\parskip}%
     \^^\{\parindent}%
     \^^\{\abovedisplayskip}%
     \^^\{\belowdisplayskip}%
     \^^\{\abovedisplayshortskip}%
     \^^\{\belowdisplayshortskip}%
     \^^\{\abovechapterskip}%
     \^^\{\belowchapterskip}%
     \^^\{\abovesectionskip}%
     \^^\{\belowsectionskip}%
     \^^\{\abovesubsectionskip}%
     \^^\{\belowsubsectionskip}}


\def\twoupsetup{
     \topmargin=.75in
     \leftmargin=.5in
     \vsize=6.9in
     \hsize=4.75in
     \fullhsize=10in
     \let\draft=\relax}

\outputstyle{normal}                             

\def\marginnoteformat{\subscriptsize             
     \hsize=1in \baselinestretch=1000 \everypar={}%
     \tolerance=5000 \hbadness=5000 \parskip=0pt \parindent=0pt
     \leftskip=0pt \rightskip=0pt \raggedright}

\head={\ifdraft\normalfonts\it\hfil DRAFT\hfil   
     \llap{\number\day\ \monthword\month\ \militarytime}\else\hfil\fi}
\foot={\hfil\normalfonts\numstyle\pagenum\hfil}  

\normalbaselineskip=12pt                         
\normallineskip=0pt                              
\normallineskiplimit=0pt                         
\normalbaselines                                 

\topskip=.85\baselineskip \splittopskip=\topskip \headskip=2\baselineskip
\footskip=\headskip

\pagenumstyle{arabic}                            

\parskip=0pt                                     
\parindent=20pt                                  

\baselinestretch=1000                            


\chapterstyle{left}                              
\chapternumstyle{blank}                          
\def\chapterbreak{\newpage}                      
\abovechapterskip=0pt                            
\belowchapterskip=1.5\baselineskip               
     plus.38\baselineskip minus.38\baselineskip
\def\chapternumformat{\numstyle\chapternum.}     

\sectionstyle{left}                              
\sectionnumstyle{blank}                          
\def\sectionbreak{\vskip0pt plus4\baselineskip\penalty-100
     \vskip0pt plus-4\baselineskip}              
\abovesectionskip=1.5\baselineskip               
     plus.38\baselineskip minus.38\baselineskip
\belowsectionskip=\the\baselineskip              
     plus.25\baselineskip minus.25\baselineskip
\def\sectionnumformat{
     \ifblank\chapternumstyle\then\else\numstyle\chapternum.\fi
     \numstyle\sectionnum.}

\subsectionstyle{left}                           
\subsectionnumstyle{blank}                       
\def\subsectionbreak{\vskip0pt plus4\baselineskip\penalty-100
     \vskip0pt plus-4\baselineskip}              
\abovesubsectionskip=\the\baselineskip           
     plus.25\baselineskip minus.25\baselineskip
\belowsubsectionskip=.75\baselineskip            
     plus.19\baselineskip minus.19\baselineskip
\def\subsectionnumformat{
     \ifblank\chapternumstyle\then\else\numstyle\chapternum.\fi
     \ifblank\sectionnumstyle\then\else\numstyle\sectionnum.\fi
     \numstyle\subsectionnum.}


\footnotenumstyle{symbols}                       
\footnoteskip=0pt                                
\def\footnotenumformat{\numstyle\footnotenum}    
\def\footnoteformat{\footnotesize                
     \everypar={}\parskip=0pt \parfillskip=0pt plus1fil
     \leftskip=1em \rightskip=0pt
     \spaceskip=0pt \xspaceskip=0pt
     \def\\{\ifhmode\ifnum\lastpenalty=-10000
          \else\hfil\penalty-10000 \fi\fi\ignorespaces}}


\def\undefinedlabelformat{$\bullet$}             


\equationnumstyle{arabic}                        
\subequationnumstyle{blank}                      
\figurenumstyle{arabic}                          
\subfigurenumstyle{blank}                        
\tablenumstyle{arabic}                           
\subtablenumstyle{blank}                         

\eqnseriesstyle{alphabetic}                      
\figseriesstyle{alphabetic}                      
\tblseriesstyle{alphabetic}                      

\def\puteqnformat{\hbox{
     \ifblank\chapternumstyle\then\else\numstyle\chapternum.\fi
     \ifblank\sectionnumstyle\then\else\numstyle\sectionnum.\fi
     \ifblank\subsectionnumstyle\then\else\numstyle\subsectionnum.\fi
     \numstyle\equationnum
     \numstyle\subequationnum}}
\def\putfigformat{\hbox{
     \ifblank\chapternumstyle\then\else\numstyle\chapternum.\fi
     \ifblank\sectionnumstyle\then\else\numstyle\sectionnum.\fi
     \ifblank\subsectionnumstyle\then\else\numstyle\subsectionnum.\fi
     \numstyle\figurenum
     \numstyle\subfigurenum}}
\def\puttblformat{\hbox{
     \ifblank\chapternumstyle\then\else\numstyle\chapternum.\fi
     \ifblank\sectionnumstyle\then\else\numstyle\sectionnum.\fi
     \ifblank\subsectionnumstyle\then\else\numstyle\subsectionnum.\fi
     \numstyle\tablenum
     \numstyle\subtablenum}}


\referencestyle{sequential}                      
\referencenumstyle{arabic}                       
\def\putrefformat{\numstyle\referencenum}        
\def\referencenumformat{\numstyle\referencenum.} 
\def\putreferenceformat{
     \everypar={\hangindent=1em \hangafter=1 }%
     \def\\{\hfil\break\null\hskip-1em \ignorespaces}%
     \leftskip=\refnumindent\parindent=0pt \interlinepenalty=1000 }


\normalsize


\def\fmtversion{2.6M (June 1992)}

\catcode`\@=12

\typesize=10pt \magnification=1200 \baselineskip17truept
\footnotenumstyle{arabic} \hsize=6truein\vsize=8.5truein
\input epsf
\sectionnumstyle{blank}
\chapternumstyle{blank}
\chapternum=1
\sectionnum=1
\pagenum=0

\def\begintitle{\pagenumstyle{blank}\parindent=0pt
\begin{narrow}[0.4in]}
\def\endtitle{\end{narrow}\newpage\pagenumstyle{arabic}}


\def\beginexercise{\vskip 20truept\parindent=0pt\begin{narrow}[10
truept]}
\def\endexercise{\vskip 10truept\end{narrow}}


\def\eql#1{\eqno\eqnlabel{#1}}
\def\ref{\reference}
\def\peq{\puteqn}
\def\pref{\putref}

\def\mgn{\marginnote}
\def\bex{\begin{exercise}}
\def\eex{\end{exercise}}


\font\open=msbm10 

 
\def\StretchRtArr#1{{\count255=0\loop\relbar\joinrel\advance\count255 by1
\ifnum\count255<#1\repeat\rightarrow}}
\def\StretchLtArr#1{\,{\leftarrow\!\!\count255=0\loop\relbar
\joinrel\advance\count255 by1\ifnum\count255<#1\repeat}}

\def\StretchLRtArr#1{\,{\leftarrow\!\!\count255=0\loop\relbar\joinrel\advance
\count255 by1\ifnum\count255<#1\repeat\rightarrow\,\,}}

\def\mbox#1{{\leavevmode\hbox{#1}}}

\def\hspace#1{{\phantom{\mbox#1}}}
\def\oZ{\mbox{\open\char90}}

\def\al{\alpha}
\def\Up{\Upsilon}
\def\bal{{\bmit\alpha}} 
\def\bbe{{\bmit\beta}} 
\def\be{\beta}

\def\de{\delta}
\def\Ga{\Gamma}

\def\la{\lambda}

\def\Om{\Omega}

\def\th{\theta}

\def\De{\Delta}

\def\det{{\rm det\,}}

\def\sc{{\rm sc }}


\def\frac#1/#2{\leavevmode\kern.1em
\raise.5ex\hbox{\the\scriptfont0 #1}\kern-.1em/\kern-.15em
\lower.25ex\hbox{\the\scriptfont0 #2}}
\def\sfrac#1/#2{\leavevmode\kern.1em
\raise.5ex\hbox{\the\scriptscriptfont0 #1}\kern-.1em/\kern-.15em
\lower.25ex\hbox{\the\scriptscriptfont0 #2}}

\def\gtorder{\mathrel{\raise.3ex\hbox{$>$}\mkern-14mu
             \lower0.6ex\hbox{$\sim$}}}
\def\ltorder{\mathrel{\raise.3ex\hbox{$<$}\mkern-14mu
             \lower0.6ex\hbox{$\sim$}}}

\def\semidirprod{\rlap{\ss C}\raise1pt\hbox{$\mkern.75mu\times$}}
\def\for{\lower6pt\hbox{$\Big|$}}
\def\fish{\kern-.25em{\phantom{abcde}\over \phantom{abcde}}\kern-.25em}


\def\boxit#1{\vbox{\hrule\hbox{\vrule\kern3pt
        \vbox{\kern3pt#1\kern3pt}\kern3pt\vrule}\hrule}}
\def\dalemb#1#2{{\vbox{\hrule height .#2pt
        \hbox{\vrule width.#2pt height#1pt \kern#1pt \vrule
                width.#2pt} \hrule height.#2pt}}}

\def\ol{\overline}
\def\frac#1#2{{{#1}\over{#2}}}

\def\noin{\noindent}

\def\comb#1#2{{\left(#1\atop#2\right)}}

\def\cot{{\rm cot\,}}

\def\eg{{\it e.g.}}
\def\ie{{\it i.e. }}
\def\cf{{\it cf }}
\def\pa{\partial}
\def\bra#1{\langle#1\mid}
\def\ket#1{\mid#1\rangle}

\def\br#1#2{\langle{#1}\mid{#2}\rangle}   


\def\Tr{{\rm Tr\,}}
\def\Ai{{\rm Ai\,}}
\def\Bi{{\rm Bi\,}}

\def\wt{\widetilde}

\def\3j#1#2#3#4#5#6{\left\lgroup\matrix{#1&#2&#3\cr#4&#5&#6\cr}
\right\rgroup}

\def\caI{{\cal I}}

\def\m?{\mgn{?}}

\def\pa{\partial}

\def\beq{\begin{eqnarray}}
\def\eeq{\end{eqnarray}}


\def\aop#1#2#3{{\it Ann. Phys.} {\bf {#1}} ({#2}) #3}
\def\cjp#1#2#3{{\it Can. J. Phys.} {\bf {#1}} ({#2}) #3}
\def\cmp#1#2#3{{\it Comm. Math. Phys.} {\bf {#1}} ({#2}) #3}
\def\cqg#1#2#3{{\it Class. Quant. Grav.} {\bf {#1}} ({#2}) #3}

\def\ijmp#1#2#3{{\it Int. J. Mod. Phys.} {\bf {#1}} ({#2}) #3}

\def\jmp#1#2#3{{\it J. Math. Phys.} {\bf {#1}} ({#2}) #3}
\def\jpa#1#2#3{{\it J. Phys.} {\bf A{#1}} ({#2}) #3}
\def\jpc#1#2#3{{\it J. Phys.} {\bf C{#1}} ({#2}) #3}
\def\lnm#1#2#3{{\it Lect. Notes Math.} {\bf {#1}} ({#2}) #3}

\def\np#1#2#3{{\it Nucl. Phys.} {\bf B{#1}} ({#2}) #3}
\def\npa#1#2#3{{\it Nucl. Phys.} {\bf A{#1}} ({#2}) #3}
\def\pl#1#2#3{{\it Phys. Lett.} {\bf {#1}} ({#2}) #3}
\def\phm#1#2#3{{\it Phil.Mag.} {\bf {#1}} ({#2}) #3}
\def\prp#1#2#3{{\it Phys. Rep.} {\bf {#1}} ({#2}) #3}
\def\pr#1#2#3{{\it Phys. Rev.} {\bf {#1}} ({#2}) #3}
\def\prA#1#2#3{{\it Phys. Rev.} {\bf A{#1}} ({#2}) #3}

\def\prD#1#2#3{{\it Phys. Rev.} {\bf D{#1}} ({#2}) #3}
\def\prE#1#2#3{{\it Phys. Rev.} {\bf E{#1}} ({#2}) #3}
\def\prl#1#2#3{{\it Phys. Rev. Lett.} {\bf #1} ({#2}) #3}

\def\rmp#1#2#3{{\it Rev. Mod. Phys.} {\bf {#1}} ({#2}) #3}

\def\zfp#1#2#3{{\it Z. f. Phys.} {\bf {#1}} ({#2}) #3}

\def\cras#1#2#3{{\it Comptes Rend. Acad. Sci. (Paris)} {\bf{#1}} (#2) #3}
\def\prs#1#2#3{{\it Proc. Roy. Soc.} {\bf A{#1}} ({#2}) #3}
\def\pcps#1#2#3{{\it Proc. Camb. Phil. Soc.} {\bf{#1}} ({#2}) #3}
\def\mpcps#1#2#3{{\it Math. Proc. Camb. Phil. Soc.} {\bf{#1}} ({#2}) #3}

\def\amsh#1#2#3{{\it Abh. Math. Sem. Ham.} {\bf {#1}} ({#2}) #3}
\def\am#1#2#3{{\it Acta Mathematica} {\bf {#1}} ({#2}) #3}
\def\aim#1#2#3{{\it Adv. in Math.} {\bf {#1}} ({#2}) #3}
\def\ajm#1#2#3{{\it Am. J. Math.} {\bf {#1}} ({#2}) #3}
\def\amm#1#2#3{{\it Am. Math. Mon.} {\bf {#1}} ({#2}) #3}

\def\aom#1#2#3{{\it Ann. of Math.} {\bf {#1}} ({#2}) #3}
\def\cjm#1#2#3{{\it Can. J. Math.} {\bf {#1}} ({#2}) #3}
\def\bams#1#2#3{{\it Bull.Am.Math.Soc.} {\bf {#1}} ({#2}) #3}

\def\cmh#1#2#3{{\it Comm. Math. Helv.} {\bf {#1}} ({#2}) #3}

\def\dmj#1#2#3{{\it Duke Math. J.} {\bf {#1}} ({#2}) #3}
\def\invm#1#2#3{{\it Invent. Math.} {\bf {#1}} ({#2}) #3}

\def\jdg#1#2#3{{\it J. Diff. Geom.} {\bf {#1}} ({#2}) #3}

\def\joa#1#2#3{{\it J. of Algebra} {\bf {#1}} ({#2}) #3}
\def\jram#1#2#3{{\it J. f. Reine u. Angew. Math.} {\bf {#1}} ({#2}) #3}
\def\jims#1#2#3{{\it J. Indian. Math. Soc.} {\bf {#1}} ({#2}) #3}
\def\jlms#1#2#3{{\it J. Lond. Math. Soc.} {\bf {#1}} ({#2}) #3}
\def\jmpa#1#2#3{{\it J. Math. Pures. Appl.} {\bf {#1}} ({#2}) #3}
\def\ma#1#2#3{{\it Math. Ann.} {\bf {#1}} ({#2}) #3}

\def\mz#1#2#3{{\it Math. Zeit.} {\bf {#1}} ({#2}) #3}
\def\ojm#1#2#3{{\it Osaka J.Math.} {\bf {#1}} ({#2}) #3}
\def\pams#1#2#3{{\it Proc. Am. Math. Soc.} {\bf {#1}} ({#2}) #3}
\def\pems#1#2#3{{\it Proc. Edin. Math. Soc.} {\bf {#1}} ({#2}) #3}

\def\plb#1#2#3{{\it Phys. Letts.} {\bf {B#1}} ({#2}) #3}
\def\pla#1#2#3{{\it Phys. Letts.} {\bf {A#1}} ({#2}) #3}
\def\plms#1#2#3{{\it Proc. Lond. Math. Soc.} {\bf {#1}} ({#2}) #3}
\def\pgma#1#2#3{{\it Proc. Glasgow Math. Ass.} {\bf {#1}} ({#2}) #3}
\def\qjm#1#2#3{{\it Quart. J. Math.} {\bf {#1}} ({#2}) #3}
\def\qjpam#1#2#3{{\it Quart. J. Pure and Appl. Math.} {\bf {#1}} ({#2}) #3}

\def\rmjm#1#2#3{{\it Rocky Mountain J. Math.} {\bf {#1}} ({#2}) #3}

\def\tams#1#2#3{{\it Trans.Am.Math.Soc.} {\bf {#1}} ({#2}) #3}

\begin{title}
\vglue 0.5truein
\vskip15truept
\centertext {\Bigfonts \bf Discrete eigenproblems} \vskip7truept
\vskip10truept\centertext{\Bigfonts }
 \vskip 20truept
\centertext{J.S.Dowker\footnote{dowker@man.ac.uk}} \vskip 7truept
\centertext{\it Theory Group,} \centertext{\it School of Physics and
Astronomy,} \centertext{\it The University of Manchester,} \centertext{\it
Manchester, England} \vskip 7truept \centertext{}

\vskip 7truept

\vskip40truept
\begin{narrow}
Schr\"odinger eigenproblems on a discrete interval are further investigated with special
attention given to test cases such as the linear and Rosen--Morse potentials. In the
former case it is shown that the characteristic function determining the eigenvalues is a
Lommel polynomial and considerable space is devoted to these objects. For example it is
shown that the continuum limit of the determinant is obtained by a transitional limit of
the Lommel polynomials for large order and argument.

Numerical comparisons between discrete approximations and (known) continuum values
for the ratio of functional determinants with and without the potential are made and show
good agreement, even for small numbers of vertices.

The zero mode problem is also briefly dealt with.
\end{narrow}
\vskip 5truept
\vskip 60truept
\vfil
\end{title}
\pagenum=0
\newpage
\section {\bf 1. Introduction}

In the course of a previous work, [\pref{Dowgelyag}], investigating the, rather basic,
problem of the analogue of the Schr\"odinger equation on a discretised interval (or
`path'), Chebyshev polynomials appeared as unperturbed propagators. The use of these
polynomials in the context of the free equation occurs in earlier works,
[\pref{Bass,CandY,BEC}]. In this communication, I wish to enlarge on the technique
employed in [\pref{Dowgelyag}] of replacing the three--term recurrence by a two by two
matrix, two--term one, as outlined in the classic Atkinson, [\pref{Atkinson}]. There is
nothing new in content, but it has a manipulative advantage, allowing basic properties
and computations to be expressed compactly.

In addition to developing this formalism (in sections 2 and 5) I present some numerical
studies of particular potentials comparing the `exact' continuum results with the discrete
`approximations'. I treat the discrete, confined, linear potential in some analytical detail
and show that the characteristic polynomial is a Lommel polynomial. These quantities
have not had many physical applications and so i present some basic, and not so basic,
facts.
\section{\bf 2. The recurrences}

The basic eigenvalue three--term recurrence, (\eg\ [\pref{Fort}]),
  $$
   y(j+1)+\big(\la-V(j)-2\big)\, y(j)+y(j-1)=0\,,
  \eql{lin3}
  $$
subject to two--point boundary conditions, say Dirichlet, Neumann or Robin,
  $$\eqalign{
  y(0)&=y(p+1)=0\,,\quad D\cr
  y(0)&=y(1),\,\,\,y(p)=y(p+1)\,,\quad N\cr
  \De y(0)&=\al\, y(0)\,,\quad \De y(p)=-\be\, y(p+1)\,,\quad R.
  }
  \eql{bcs}
  $$
at the boundary, $\pa\caI$, of the discrete interval, $\caI$, which is comprised of the
$(p+2)$ points $j=0,1,\ldots,p,p+1$. $\pa\caI= 0\cup p+1$.

This is a discrete Sturm--Liouville problem.

$V(j)$ is the potential because  (\peq{lin3}) can be rewritten as the more familar
looking Laplacian eigenvalue equation,\footnote{ $\nabla$ is the backwards difference
operator. There are many discrete approximations to the continuum $-y''(x)+\ol
V(x)\,y(x)=\la\,y(x).$}
  $$
  \big[-\nabla\De + V(j)\big]y(j)=\la\,y(j)\,.
  \eql{lap}
  $$

As well as the homogeneous equation (\peq{lap}), the inhomogeneous one
  $$
  \big[-\nabla\De + V(j)-\la\big]y(j)=\rho(j)\,.
  \eql{inhom}
  $$
is of interest where $\rho(j)$ is a source density and $\la$ a spectral parameter.

As intimated in the Introduction, it is helpful to recast these three--term relations by
two--term matrix ones and I have to use a little space to set this up.

For example (\peq{lin3}) (\ie (\peq{lap})), and (\peq{inhom}) can be replaced by the
matrix recurrences,
  $$
  \Up(j)=C(x_j)\Up(j-1)\equiv C(j)\Up(j-1)
  \eql{hom}
  $$
and
  $$
  \Up(j)= C(j)\Up(j-1)+D\,R(j-1)
  \eql{inhom2}
  $$
respectively, with the definitions
   $$
  \Upsilon(j)=\left(\matrix{y(j)\cr y(j+1)}\right)\,,\quad
  R(j)=\left(\matrix{\rho(j)\cr \rho(j+1)}\right)
  \eql{matrix}
  $$
and
  $$\eqalign{
  C(j)&=\left(\matrix{0&1\cr-1 &V(j)+2-\la}\right)
  =\left(\matrix{0&1\cr-1 &\phi(j)}\right)\,.\cr
  }
  \eql{cee2}
  $$
 The driving matrix, $C$, can be split as $C(j)=B(j)-\la D(j)$ where,
  $$\eqalign{
  B(j)&=\left(\matrix{0&1\cr-1 &V(j)+2}\right)\,,\quad
  D(j)=D=\left(\matrix{0&0\cr0 &1}\right)\,.
  }
  \eql{eba}
  $$

If $\Up_1$ and $\Up_2$ are two general solutions of the difference equation,
(\peq{hom}), the constancy of their Casoratian,
  $$
  \wt\Up_1(j)\,J\,\Up_2(j)=\det \big(\Up_1\otimes\Up_2\big)\,
  \eql{casor}
  $$
follows immediately, where $J$ is the symplectic metric,
  $$
  J=\left(\matrix {0&1\cr-1&0}\right)\,.
  $$

The replacement of a second order difference equation by two, first order ones, organised
into matrix form, is the simplest example of the more general case where an $n$th order
equation is replaced by an $n\times n$  matrix equation of the first order. Birkhoff,
[\pref{Birkhoff}], for example, makes extensive use of this formulation. Milne--Thomson,
[\pref{Milne-Thomson}],\S5.3, \S12.8,  also employs the same technique  in the
language of continued fractions. See also [\pref{Miller2}] and the compact discussion in
[\pref{LandT}] extended in  Agarwal, [\pref{Agarwal}].

Taking the discrete interval $\caI$ as the archetypal domain, the propagation
(\peq{hom}) can be `solved' as,
  $$
  \Up(j)=K(j,j')\Up(j')=K(j,0)\Up(0)\,,
  $$
where $\Up(0)$ is the given initial vector and $K(j,j')$ is the forwards propagator,
  $$
  K(j,j')=\th(j,j')\,C(j)\,C(j-1)\ldots C(j'+1)\,,
  \eql{kprod}
  $$
which I concentrate on.

Equivalent to (\peq{hom}) is
  $$
  K(j,j')=\de_{jj'}{\bf 1}+C(j)K(j-1,j'),
  \eql{simprec}
  $$
incorporating the initial condition $K(j,j)={\bf 1}$. In particular, for $K(j)\equiv K(j,0)$
one has
  $$
  K(j)=\de_{j0}{\bf 1}+C(j)K(j-1)\,.
  \eql{simprec2}
  $$

In my previous work, I treated the potential, $V(j)$, as a perturbation. This amounts to
extracting the factor of $2-\la$  in $C(j)$ which produces a rearrangement of
(\peq{simprec}) with Chebyshev polynomials  as unperturbed propagators. This is formally
attractive, but, if all one is interested in is a numerical answer, then a direct iteration of
(\peq{simprec2}) is adequate. Even though it is absolutely equivalent to (\peq{lin3}), it is
computationally advantageous. Furthermore, in order to implement Dirichlet boundary
conditions, say, it is sufficient to find the bottom right--hand component of $K(p+1)$, \ie
$\Tr\big(K(p+1)\,D\big)$, and the relevant recurrence is
  $$
  K_D(j)=D\,\de_{j0}{\bf 1}+C(j)K_D(j-1)\,,
  \eql{simprec3}
  $$
where $K_D(j)\equiv K(j)\,D$. This is easily programmed.

\section{\bf 3. An explicit example. The linear potential.}

In the continuum, Kirsten and McKane, [\pref{KandM}], as a simple illustration of the
Gel'fand--Yaglom procedure, compute the ratio of Laplacian functional determinants with
and without a linear potential, specifically for the operator,
  $$
   L(x)\equiv -{d^2\over dx^2}+\ol V(x)\,,
   \eql{Lop}
  $$
with $\ol V(x)=b^3\,x$. They chose $b=1$, but I use a general strength.

As is well known, Airy functions arise for such a potential and Kirsten and McKane give,
on this basis, $1.085(339648)$ for the determinant ratio.

I outline some details later in this section but first I give the discrete version of this
computation. I proceed to write down the relevant recurrence in scalar form, \ie\
equation (\peq{lin3}) with the potential linear in $j$, $V(j)=Bj$,
  $$
   y(j+1)+\big(\la-Bj-2\big)\, y(j)+y(j-1)=0\,,
  \eql{lin4}
  $$
with, for simplicity, Dirichlet boundary conditions at $j=0$ and $j=p+1$ according to
(\peq{bcs}).

One might consider this system as describing the motion of a charged particle in a
uniform electric field, confined to a box, \ie\ interval.

This recurrence is an example of a more general type analysed by Boole, [\pref{Boole}],
using operator methods and, more analytically, by Barnes, [\pref{Barnesrec}], but, again
I will not use their results. Rather, in the particular case of (\peq{lin4}), I refer, at the
moment, to Bleich and Melan, [\pref{BlandM}], who give a rather detailed treatment of
recursions of the form,
   $$
   y(j+1)-\phi(j)\, y(j)+y(j-1)=0\,,
  \eql{lin5}
  $$
\ie (\peq{hom}) with (\peq{cee2}).

Starting from a solution Ansatz suggested by the iteration of (\peq{lin5}), or
(\peq{hom}), they show that, for a necessarily restricted set of functions, $\phi(j)$,
explicit solutions are possible.

In particular a Dirichlet solution for $\phi(j)=A+B\,j$ is presented in equation (59) on
p.149 in [\pref{BlandM}]. In my notation $A=2-\la$ and I set their initial point, $\eta$
to 0. For the moment, I just quote their forwards solution which transcribes (with some
inessential notational adjustments) to \footnote{ I draw attention to the definition of the
product symbol in [\pref{BlandM}] on p.141 and footnote on p.151, which differs from
the usual one used here.},
  $$\eqalign{
  y(j,\la)=\sum_{k=0}^{j-1}\cos{k\pi\over2}\comb{j-k/2-1}{k/2}\prod_{l=1}^{j-k-1}
 \big[A+B(k/2+l)\big]\,.\cr
 }
 \eql{BMsol}
  $$

I make further comments on this expression in sections 6 and 7 where I relate it to
Lommel polynomials.

To give a flavour of the structure of this solution, I compute $y(4,\la)$, which corresponds
to the terminal value if I choose $p$ to equal $3$.

Because of the cosine, $k$ can only be even and therefore only $2$ or $0$ so that,
  $$\eqalign{
  y(4,\la)
  &=(2-\la+2B)\big((2-\la+B)(2-\la+3B)-2\big)\cr
&=-\la^3+\la^2(6+S_1)-\la(10+S_2+4S_1)+4+4S_1+2S_2+S_3-4B\,,
  }
  \eql{y4}
  $$
where, for comparison purposes, I have defined $v_1=1B$, $v_2=2B$ and $v_3=3B$ as
the values of the linear potential $Bj$ at $j=1,2,3$ and also
  where
  $$
  S_1=v_1+v_2+v_3\,,\quad S_2=v_1v_2+v_1v_3+v_2v_3\,,\quad S_3=v_1v_2v_3\,.
  \eql{esses}
  $$

I have done this in order to compare with the expression derived in [\pref{Dowgelyag}] for
a general potential. The comparison is exact if one takes into account the special relation
here that $\nabla\De v_2=v_1-2v_2+v_3=0$ for a linear potential.

The eigenfunction, for a given eigenvalue, say $\la_n$, is given by $y(j,\la_n)$,
(\peq{BMsol}), where $n$ ranges over $p$ values, say $1\le n\le p$, (setting  the
terminal point to $j=p+1$).

It is interesting to note that if $p$ is odd, there is always an eigenvalue, halfway along,
$\la_{(p+1)/2}=2+B(p+1)/2$ for which, as is obvious and can be checked, first order
perturbation theory is exact.

Finally, from general theory, the ratio of determinants is given by,
  $$
    \ol\det(B)\equiv{\det(B)\over \det(0)}={y(p+1,0)\over y(p+1,0)\big|_{B=0}}\,,
    \eql{detrat}
  $$
which can be evaluated from (\peq{BMsol}).

In order to compare with the {\it continuum} linear potential, $\ol V=b^3x$, the
strength, $B$ of the discrete potential has to be rescaled. Over the unit interval, one has
$B=\big(b/(p+1)\big)^3$ by dimensions.

For $b=1$, I get a determinant ratio value of $1.08533860$, choosing
$p=300$.\footnote{ In 1.4 secs using the wxMaxima CAS on an Athlon ii x 4 machine.}
This is not particularly efficient, but comparison with the continuum value,
$1.085339648$, demonstrates the validity of the methods. More interesting, perhaps, is
that for just one interpolating step, \ie $p=1$, the value is $1.0625$, for $p=3$ it is
$1.07947$ while for ten, I find $1.08456$, a reasonable approximation.

These results suggests that, if one is prepared to sacrifice a little accuracy, the determinant
ratio for any potential can be calculated almost by hand.

Incidentally, for an {\it attractive} potential ($b<0$), bound states arise and in Fig.1 I
plot the determinant ratio that demonstrates this fact by the zeros. Not so interestingly,
for a repulsive potential, the ratio shows a rapid, monotonic increase. The continuum
curve is also shown in the figure.

\epsfxsize=5truein \epsfbox{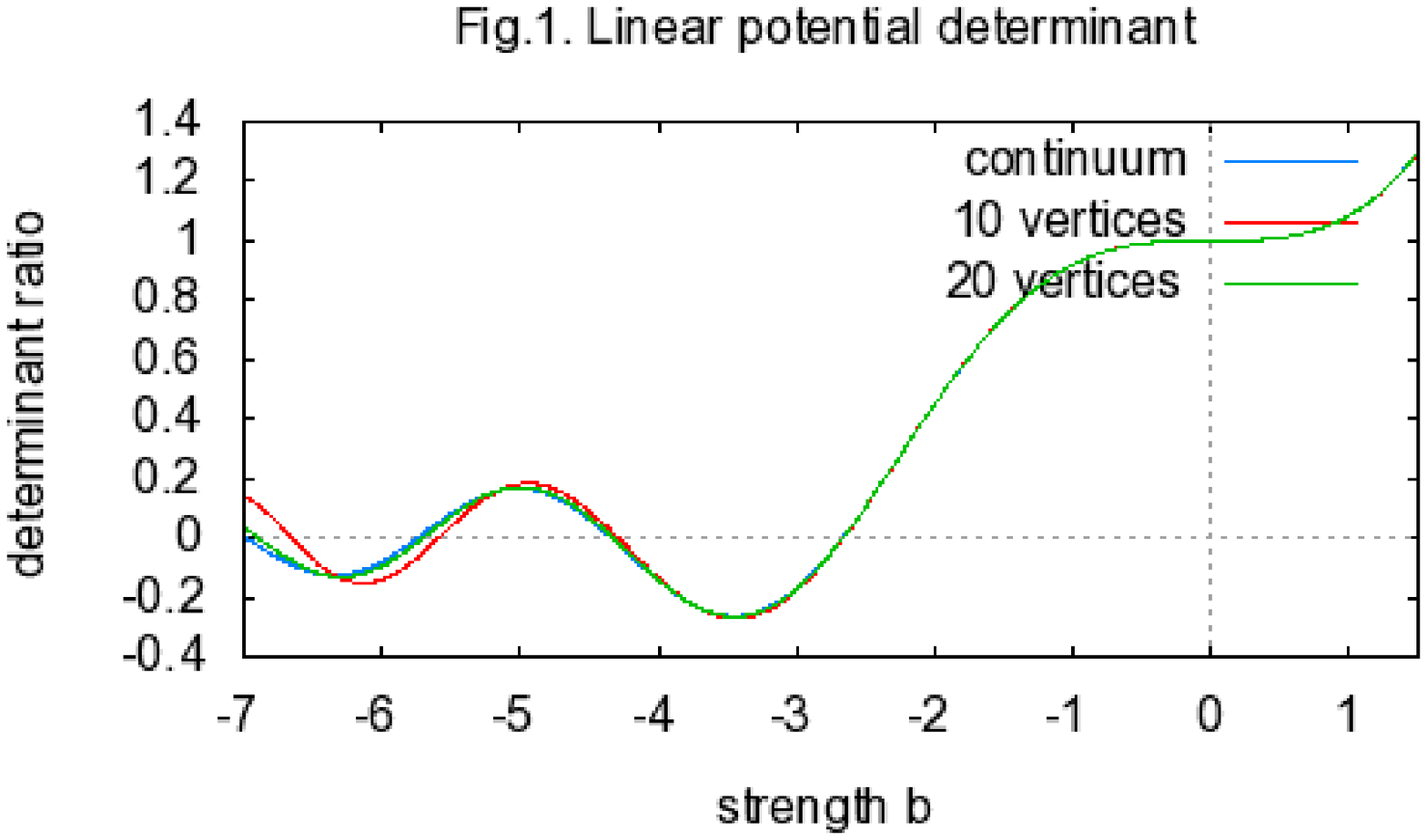}

For the continuum situation, the Gel'fand--Yaglom function for the potential $b^3\,x$ is,
\cf\ [\pref{KandM}],
  $$
   y_1(x)={\pi\over b}\big(\Bi(0)\,\Ai(b\,x)-\Ai(0)\,\Bi(b\,x)\big)\,,
   \eql{airy}
  $$
and yields the ratio of determinants,
  $$
  {\pi\over b}\big(\Bi(0)\,\Ai(b))-\Ai(0)\,\Bi(b)\big)\,,
  \eql{airy2}
  $$
which is also plotted in Fig.1. The figure also shows that the discrete approach is good for
smallish strengths, $b$, as expected since it is a perturbation--like method yielding a
finite polynomial (in $b$) approximation for the Airy expression, (\peq{airy}). This is
discussed further in sections 6 and 7.

\section{\bf 4. The Rosen--Morse potential}

Bound states also arise for the P\"oschl--Teller, hyperbolic potential,
  $$
  \ol V(x)=-{l(l+1)\over\cosh^2(x)}\,,
  \eql{PT}
  $$
confined to a box, as considered by Dunne, [\pref{Dunne}].

From its quantum mechanical origins in molecular physics, the strength, $l$, is usually
taken to be integral, but can, in fact, be any real number, in which case it might be
preferable to use the variable $\ol l=l+1/2$. I will still employ $l$ but refer to
(\peq{PT}) as the Rosen--Morse potential (\eg\ Kleinert, [\pref{Kleinert}], Fl\"ugge,
[\pref{Flugge}]).

Two particular solutions of the equation, $L(x)\,\psi(x)=0$,\footnote{For simplicity, I set
the mass equal to zero.} are the Legendre functions, $P_l(\tanh x)$ and $Q_l(\tanh x)$,
constituting a fundamental set. I choose the interval $-1/2\le x\le 1/2$ as domain so
that I can give the Gel'fand--Yaglom function as written down by Dunne,
  $$
  \psi_{GY}(x)=P_ l(-\tanh1/2)\,Q_l(\tanh x)-P_l(\tanh1/2)\,Q_ l(-\tanh x)
  $$
so that the ratio of determinants, with and without $\ol V$, is
  $$
{\ol\det_C}\,(l)=P_l(-t)\,Q_l(t)-P_ l(t)\,Q_l(-t)
  \eql{contdet2}
  $$
where $t=\tanh1/2\approx 0.46211715726$. Properties of the Legendre functions check
that $\ol\det_C(l) $ is symmetrical in $\ol l$, as must be, and one need compute it only
for positive $\ol l$, \ie\ $l\ge -1/2$.

In terms of functions of the first kind, I find,
  $$
{\ol\det_C}\,(l)={\pi\over\sin\pi l}\,\big(P_ l(-t)^2-P_ l(t)^2\big)\,
  \eql{contdet}
  $$
which is easily calculated using Gauss' hypergeometric form, say,
  $$
 P_ l(t)={}_2F_1\big(- l,1+ l;1;(1-t)/2\big)\,.
  $$

Near integral $l$, the numerics become uncertain. But for $l$ actually integral, the ratio of
determinants, (\peq{contdet2}), can be expressed purely in terms of the polynomials,
$P_l$,
  $$
{\ol\det_C}\,( l)=(-1)^l P_l\bigg(P_l-2\sum_{m=1}^l
{1\over m}P_{m-1}\,P_{n-m}\bigg)\,,\quad l\in\oZ\,,
  \eql{contdet3}
  $$
which is a polynomial in $\tanh1/2$ and provides a check of the numbers. As an exercise,
one could use these values as the basis of an interpolation.

Turning now to the discrete calculation, the formula I use is that developed in
[\pref{Dowgelyag}]. For $p=3$, an interval with three interior vertices and step size
$h=1/4$, the determinant ratio is,
  $$
\ol\det_D(l)= 1+{1\over4}(3S_1+2S_2+S_3+v_2)\,,
\eql{discdet}
  $$
where the $S_i$ are given by (\peq{esses}) with $v_j\equiv V(j)$. (Refer to
(\peq{lap})). This is to be compared, numerically,  with (\peq{contdet}).

On this side of the calculation, I use the interval $0$ to $1$ so that the discrete potential
corresponding to (\peq{PT}) is given by
  $$
  V(j)=-{h^2l(l+1)\over\cosh^2(hj-1/2)}
  $$
and the evaluation of (\peq{discdet}) is relatively elementary.  The results are presented
graphically in fig.2, where I have also displayed the 5 vertex values and the continuum,
`exact' curve, based on (\peq{contdet}).

The figure shows again that the discrete approach is good for smallish strengths, $l$. It
yields a finite polynomial (in $l$, even in $\ol l$) approximation for the Legendre
expression, (\peq{contdet2}). The exact curve oscillates for ever and, while a finite
polynomial cannot do so, it  follows fairly well, for a decent range of $l$. Note also that, for
a repulsive potential, the ratio is bigger than one.

\epsfxsize=5truein \epsfbox{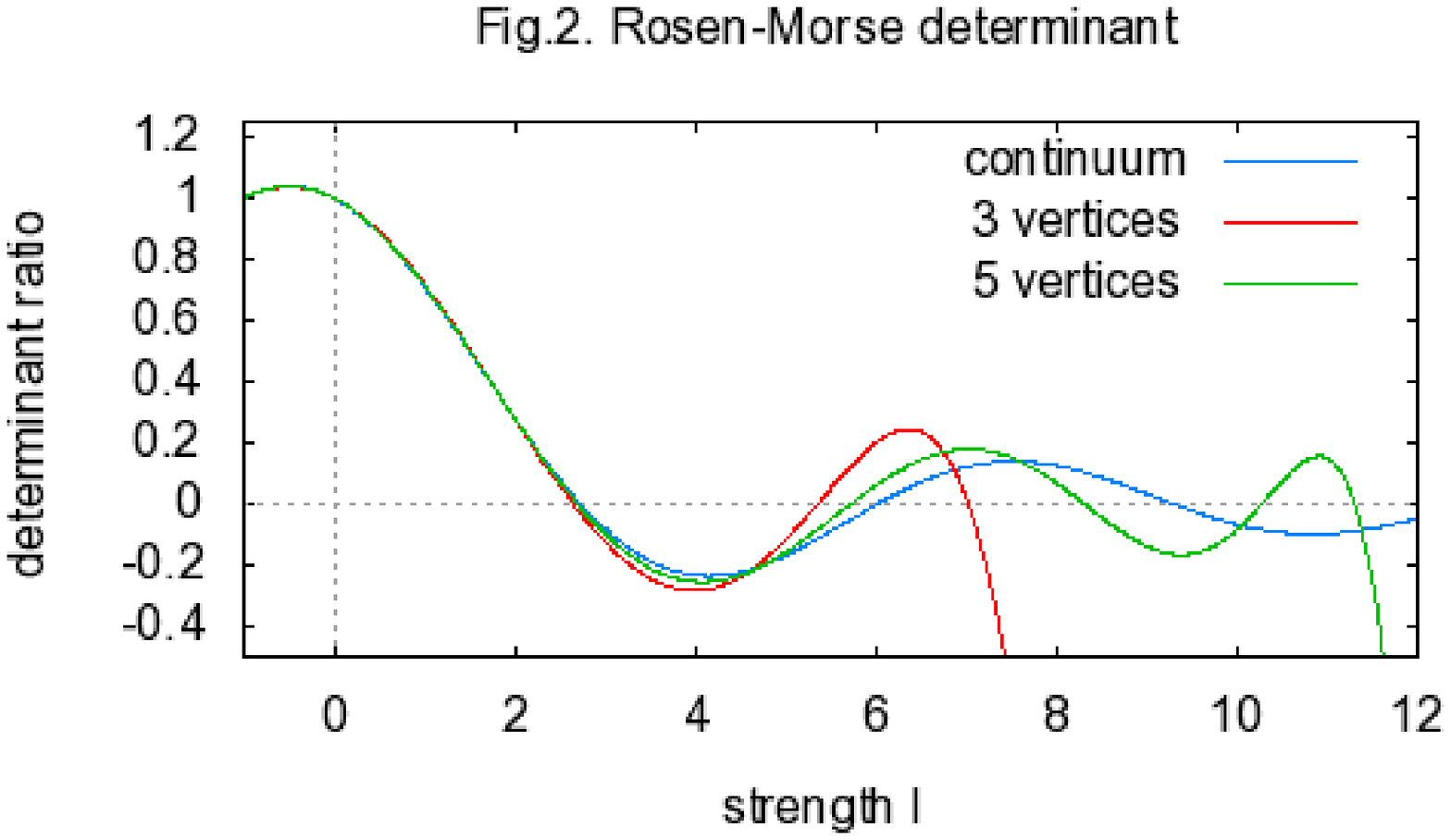}

The zeros of the curves occur when bound states appear as the well deepens.

\section{\bf 5. Analytical extensions}

So far it has been assumed that label $j$ takes just integral values. However in classical
finite difference theory, \eg\ Boole, [\pref{Boole}], this is only a particular case of a more
general situation where one deals with functions $y(x)$ of a variable, $x$ (real or
complex) satisfying a difference equation, say,
  $$
   y(x+1)-\phi(x)\, y(x)+y(x-1)=0\,,
  \eql{lin7}
  $$
as a relevant, typical example. In accordance with my previous formalism, I write this,
with, therefore, some repetition from section 2, as,
   $$
  \Up(x)=C(x)\Up(x-1)
  \eql{recurs}
  $$
with
   $$
  \Upsilon(x)=\left(\matrix{y(x)\cr y(x+1)}\right)\,,
  \eql{matrix2}
  $$
and
  $$\eqalign{
  C(x)&\equiv\left(\matrix{0&1\cr-1 &\phi(x)}\right)\,.\cr
  }
  \eql{cee3}
  $$

Assuming some starting point, $x=x_0$, the general solution of (\peq{recurs}) is
  $$\eqalign{
  \Up(x)&=C(x)C(x-1)\ldots C(x_0+1)\Up(x_0)\cr
  &=K(x,x_0)\Up(x_0)
  }
  \eql{gensol}
  $$
where $\Up(x_0)$ can be considered as an arbitrary constant vector, say $\Om$.

For a particular $x_0$, the domain of $x$ is determined, but one can consider $x_0$ to
change, and with it, the domain of $x$, (see \eg\ Boole, [\pref{Boole}], p.102, Levy and
Lessman, [\pref{LeandL}], pp.89+) (always, however, with $x-x_0$ integral). This can
be accommodated by making the constant, $\Om$, a function,
  $$
  \Om(x)=\left(\matrix{\varpi_\al(x)\cr\varpi_\be(x)}\right)\,,
  $$
where the $\varpi_i$ are arbitrary unit periodics, B\"ohmer, [\pref{Bohmer}],
[\pref{BlandM}], \S26, [\pref{LeandL}], [\pref{Milne-Thomson}], \S11.1, (\eg\
$\varpi_i(x_0)=\varpi_i(x)$) so that the general solution, (\peq{gensol}), is,
  $$\eqalign{
  \Up(x)&=C(x)C(x-1)\ldots C(x_0+1)\,\Om(x_0)\cr
  &=K(x,x_0)\,\Om(x_0)\cr
  &=\varpi_\al(x)\Up_\al(x)+\varpi_\be(x)\Up_\be(x)\,,
  }
  \eql{gensol2}
  $$
where,
  $$
\Up_\al(x)=K(x,x_0)\bal\,,\quad \Up_\be(x)=K(x,x_0)\bbe\,,
\eql{fsols}
 $$
with the standard basic vectors,
  $$
\bal=\left(\matrix{1\cr0}\right)\,,\quad \bbe=\left(\matrix{0\cr1}\right)\,.
\eql{sbv}
  $$
For example, Dirichlet initial conditions require $\varpi_\al(x_0)=0$.

The solution (\peq{gensol2}) is just the general solution to the second order difference
equation (\peq{lin7}), a standard construct, \eg\ Levy and Baggott, [\pref{LandB}],
Birkhoff, [\pref{Birkhoff}].

The basis, $\Up_\al$ and $\Up_\be$ is a {\it canonical} one (\eg\ Lakshmikantham and
Trigiante, [\pref{LandT}], Lemma 2.1.1).

$\Up_\be$ coincides with  Atkinson's {\it standard solution}, $y_n(\la)$,
[\pref{Atkinson}], equns. (4.1.4), (4.1.5) and $\Up_\al$ with the standard solution,
$z_n(\la)$, equns. (4.2.6), (4.2.8).

As a non--trivial example, the solution (\peq{BMsol}) to (\peq{lin5}), corresponds to the
Dirichlet solution $\Up_\be(j)$ because $y(0,\la)=0$ and $y(1,\la)=1$. In the following
section, I enlarge on this solution.

The constancy of the Casoratian, (\peq{casor}), now becomes the statement that it is a
unit periodic,
  $$
  \wt\Up_1(x)\,J\,\Up_2(x)=\varpi(x)\,,
  $$
where $\Up_1$ and $\Up_2$ are two {\it general} solutions of the difference equation,
(\peq{recurs}) each of the form (\peq{gensol2}), say.

\section{\bf 6. Bessel functions}

I return to a consideration of the confined linear potential of section 3, which forms an
interesting, if somewhat specific, case but with connections to some special functions.

The analysis of the second order equation by Barnes and by Bleich and Melan, proceeds
from the classic form of the recurrence,  (\peq{lin7}), so I present their results
accordingly, despite my preference for the matrix form.

The equation in question is (\cf\ (\peq{lin4})), \footnote { This is recognised as a
recurrence relation for Bessel functions but it is not necessary to do so at this time.}
  $$
   y(x+1)-(Bx+A)\, y(x)+y(x-1)=0\,,
  \eql{lin8}
  $$
which, treated using the textbook Laplace method, \footnote{ A simple discussion is
given in the elegant little book by Miller, [\pref{Miller2}].} yields two independent
solutions as contour integrals coinciding, to a simple factor, with Bessel functions,
  $$\eqalign{
  y^{(1)}(x)&=(-1)^{x+A/B}\,J_{-(x+A/B)}\big({2/B}\big)=
  J_{-(x+A/B)}\big(-{2/B}\big)\cr
  y^{(2)}(x)&=J_{(x+A/B)}\big({2/B}\big)\,,
  }
  \eql{bessls}
  $$
and the general solution is
  $$
   y(x)=\varpi_1(x)\,y^{(1)}(x)+\varpi_2(x)\,y^{(2)}(x)\,.
   \eql{gensol4}
  $$

If $\xi=x+A/B$ is an integer then these two solutions are equivalent and, in the usual
fashion, a second solution is obtained by a limiting process, $\xi\to n$,  on a Bessel
function of the second kind, due to Weber and Schl\"afli,
  $$
  Y_\xi(z)=\cot\xi\pi\,J_\xi(z)-{1\over\sin\xi\pi}\,J_{-\xi}(z)\,,
  \eql{2sol}
  $$
which satisfies the difference equation (\peq{lin7}) as it is of the form (\peq{gensol4}).

As an application, return to the conventional case when $x=j\in\oZ$ with the previously
used discrete interval, $\caI$, as domain and attempt to recover earlier results, such as
(\peq{BMsol}). Thus, apply D conditions at $j=0$ and $j=p+1$. This produces the
non--local Casorati determinant condition,
  $$
 \left|\matrix {y^{(1)}(0)&y^{(2)}(0)\cr y^{(1)}(p+1)&y^{(2)}(p+1)}\right|
  =0\,,
  \eql{nlcasdet}
  $$
which is just the vanishing of  the Dirichlet  solution,
  $$
  y(x,\la)=y^{(2)}(0)\,y^{(1)}(x)-y^{(1)}(0)\,y^{(2)}(x)
  $$
at the terminal point, $x=j=p+1$.

Inserting the explicit functions, (\peq{bessls}) gives the quantity, \footnote{ The
relation is $y(p+1,\la)\propto (-1)^\nu\,W(\nu,p)$.}

  $$
  W(\nu,p)\equiv J_{-\nu}J_{p+1+\nu}-
  {(-1)^{p+1}}J_{-(p+1+\nu)}\,J_{\nu}=0\,,\quad
 p\in\oZ\,,
 \eql{wronsk}
  $$
where $\nu\equiv A/B$ and the arguments of the Bessel functions are all $2/B\equiv z$ .
I have dropped inessential, common overall constants and leave the dependence of $W$,
which is hereby defined, on $z$ implied, until later.

The characteristic equation, (\peq{nlcasdet}), in $\la$ (remember, $A=2-\la$) provides
the eigenvalues through its roots, and is easily computable since, in this case, two
independent solutions are known.

The aim, then, is to show that the Casoratian, (\peq{nlcasdet}), is proportional to a finite
polynomial in $\nu$ of order $p$, \cf\ (\peq{BMsol}). This can be done via the series
form of the Bessel functions in (\peq{wronsk}). Of course, it is to be anticipated that the
Casoratian vanishes when $\nu$ is an integer for then the two Bessel functions are
proportional. To allow for this, it is convenient to consider the normalised
Casoratian,\footnote{ An alternative procedure is given in the next section.}
  $$
  \ol W(\nu,p)\equiv {W(\nu,p)\over W(\nu,0)}\,,
  $$
which removes the relevant factor.

A standard formula is,\footnote{$p=0$ means there ar no (internal) vertices.}
  $$
  W(\nu,0)=J_{-\nu}J_{1+\nu}+J_{-1-\nu}\,J_{\nu}=-{B\over\pi }\sin \nu\pi\,.
 \eql{wronsk3}
  $$

I give the calculation and refer to comments in sections 7 and 8. The composition relation
is used to give the series,
  $$
   J_{-\nu}J_{p+1+\nu}=\sum_{n=0}^\infty {(-1)^n B^{-(2n+p+1)}
   (p+n+2)_n\over n!\Ga(n+1-\nu)\Ga(p+\nu+n+2)}
  $$

  $$
   J_{\nu}J_{-p-1-\nu}=\sum_{n=0}^\infty {(-1)^n B^{-(2n-p-1)}
   (n-p)_n\over n!\Ga(n+1+\nu)\Ga(n-p-\nu)}\,.
  $$

In the combination, (\peq{wronsk}), most of the terms cancel, leaving the finite sum,
  $$\eqalign{
  W(\nu,p)&=\sum_{n=1}^{p+1}\comb{p+1-2n}{-n}
  {(-1)^n B^{-(p+1-2n)}
  \over\Ga(-\nu-n+1)\,\Ga(p-n+\nu+2)}
  }
  \eql{wronsksum}
  $$
Part of the summation range gives zero because
  $
  \comb{p+1-2n}{-n}=0\,,
  $
if $n\le (p+1)/2$.

I give the details that turn (\peq{wronsksum}) into (\peq{BMsol}), up to a factor. The
Gamma functions combine as,
 $$\eqalign{
      {1  \over\Ga(-\nu-n+1)\Ga(p-n+\nu+2)}&={(\nu+p-n+2)\ldots(\nu+n-1)
      \over\pi}\,\sin\nu\pi\cr
      &={\big(A+B(p-n+2)\big)\ldots\big(A+B(n-1)\big) \over \pi\,
      B^{2n-p-2}}\,\sin\nu\pi
      }
      \eql{gaprod}
 $$
and the normalised Casoratian is,
  $$\eqalign{
  \ol W(A,B,p)&=\sum_{n=1}^{p+1}\comb{p+1-2n}{-n}
  (-1)^{n+1}\big(A+B(p-n+2)\big)\ldots\big(A+B(n-1)\big)\,,
  }
  \eql{wronsksum2}
  $$
The number of brackets being $(p-2s)$.

I have changed arguments to allow for the $z$ dependence. Thus, $\ol W(A,B,p)\equiv
\ol W(\nu,p)(z)$.

The change of summation variable to,
$$
  k=2(p+1-n)\,,
$$
then yields agreement with (\peq{BMsol}) (if $j$ is chosen to be the end boundary point,
$p+1$) since,
  $$
  \comb{p+1-2n}{-n}=(-1)^{p+1-n}\comb{n-1}{2n-p-2}\,.
  $$

Bleich and Melan's solution, (\peq{BMsol}) has thus been obtained by a more direct and
particular route.

Finally the trivial transformation of summation variable $s=p+1-n$ converts
(\peq{wronsksum2}) into, the simpler--looking form,
  $$\eqalign{
  \ol W(A,B,p)&=(-1)^p\,\sum_{s=0}^{<(p+1)/2}\comb{p-s}{s}
  (-1)^{s}\big(A+B(p-s)\big)\ldots\big(A+B(s+1)\big)\,,
  }
  \eql{wronsksum3}
  $$
where the number of brackets is equal to $(p-2s)$.

The case $p=1$,
  $$
  J_{\nu+2}(z)\,J_{-\nu}(z)-J_{-(\nu+2)}(z)\,J_\nu(z)=
  -{4(\nu+1)\sin\nu\pi\over\pi z^2}\,,
  $$
is given in Gray and Matthews, [\pref{GrandM}], p.241, Ex.1. Other examples are,
\footnote{ For odd $p$, the vanishing of the overall factor $(\nu+(p+1)/2))$
corresponds to a cancellation between the two terms on the left--hand side which does
not require using $J_n=J_{-n}$.}
   $$\eqalign{
  J_{\nu+3}(z)\,J_{-\nu}(z)+J_{-(\nu+3)}(z)\,J_\nu(z)
  &={2\big(z^2-4(\nu+1)(\nu+2)\big)\sin\nu\pi\over\pi z^3}\cr
  J_{\nu+4}(z)\,J_{-\nu}(z)-J_{-(\nu+4)}(z)\,J_\nu(z)
  &=-{8(\nu+2)\big(z^2-2(\nu+1)(\nu+3)\big)\sin\nu\pi\over\pi z^4}\,.
  }
  \eql{cases}
  $$

One could regard the solution (\peq{BMsol}) as an alternative source of these results
which actually are most easily obtained by direct iteration of (\peq{lin8}) to give
$\Up_2(j)$ of (\peq{fsols}) as described in the next section.

I now place these evaluations in a more general and historic setting.
\section{\bf 7.  Lommel polynomials}

There is no obligation to identify solutions of the recurrence equation, (\peq{lin8}),
specifically with Bessel functions but making the identification saves a lot of effort due to
an extensive literature. For example, the algebra leading to (\peq{wronsksum}), and the
further extraction of the factor $\sin\nu\pi$, is well known. What I have termed the
normalised Casoratian is a Lommel polynomial (Lommel, [\pref{Lommel}], Nielsen
[\pref{Nielsen}], Graf and Gubler, [\pref{GandG}], Watson, [\pref{Watson2}]) as I now
show.

Lommel defined the polynomials, $R$, from the iteration of the particular Bessel
recurrence,
   $$
   F_{\xi+1}(z)={2\xi\over z}\,F_\xi(z)-F_{\xi-1}(z)
   \eql{bessrec}
   $$
specifically for $F_{\xi}(z)=J_{\xi}(z)$ (not the only choice).

Iteration yields the formal solution in terms of two initial values as the reduction formula,
\footnote{ I use Nielsen's notation, [\pref{Nielsen}], which is the same as that of Graf
and Gubler, [\pref{GandG}], who write $\xi\atop{P_n}$. It differs from Lommel's and
Watson's only in the ordering of the indices and a shift of unity in $\xi$.}
  $$\
  F_{\xi+n}(z)=R^{\xi-1,n}(z)\,F_\xi(z)-R^{\xi,n-1}(z)\,F_{\xi-1}(z)\,,
  $$
\ie
  $$\
  F_{\xi}(z)=R^{\,\xi-n-1,n}(z)\,F_{\xi-n}(z)-R^{\,\xi-n,n-1}(z)\,F_{\xi-n-1}(z)\,,
  \eql{frec}
  $$
where, to confirm the connection with my previous notation in (\peq{lin8}), $z=2/B$
$\xi=x+A/B$ and $y(x,A,B)=F_{x+A/B}(2/B)$. I have included the parameters $A$ and
$B$ in the $y(x)$. Either $\xi$ or $x$ can be used as the current, discrete coordinate.

Writing (\peq{frec}) in terms of the initial point $\xi_0$, \ie eliminating
$n=\xi-\xi_0-1$,
  $$\
  F_{\xi}(z)=R^{\,\xi_0,\,\xi-\xi_0-1}(z)\,F_{\xi_0+1}(z)-
  R^{\,\xi_0+1,\,\xi-\xi_0-2}(z)\,F_{\xi_0}(z)\,.
  \eql{frec2}
  $$

Comparison with (\peq{fsols}) shows that $R^{\,\xi_0,\,\xi-\xi_0-1}(z)$ corresponds to
the $\be$, or Dirichlet, solution, and $- R^{\,\xi_0+1,\,\xi-\xi_0-2}(z)$ to the $\al$
one, as functions of $\xi$. Checking this requires use of Graf's relation, [\pref{GandG}],
  $$
  R^{\,\nu,-n-1}=-R^{\,\nu-n,n-1}
  $$
which implies $R^{\,\nu,-1}=0$ and $R^{\,\nu,-2}=-1$ \eg\  [\pref{GandG}], p.103.

Lommel deduced the actual form of the polynomials, $R^{\xi,n}(z)$, by an inductive
procedure repeated by Graf and Gubler, [\pref{GandG}] vol.1 chap.1 \S2. A simpler
method, \eg\ Nielsen [\pref{Nielsen}], is that give in section 5, since easy algebra leads to
Lommel's equation,
  $$
  J_{-\nu}(z)\,J_{\nu+p+1}(z)+(-1)^p\,J_\nu(z)\,J_{-\nu-p-1}(z)=-{2\sin\pi\nu\over\pi z}
 \, R^{\nu,p}(z)\,,
 \eql{Leq}
  $$
and the left--hand side can be evaluated from the known series form of $J_\nu(z)$. The
details in section 5 are, more or less, just a rewriting of Nielsen's algebra and  I also give
the standard result from this old literature, [\pref{Lommel}], [\pref{Nielsen}],
  $$
  R^{\nu,p}(z)=\sum_{s=0}^{<\,(p+1)/2}\,{(-1)^s (p-s)!\over s!}\comb{\nu+p-s}{p-2s}
  \bigg({2\over z}\bigg)^{p-2s}\,,
  \eql{Lommp}
  $$
identical to (\peq{wronsksum3}).

A comparison of (\peq{Leq}) with (\peq{wronsk}) gives the connection,
  $$
W(\nu,p)=-{2\sin\pi\nu\over\pi z} \, R^{\nu,p}(z)\,,\quad z=2/B,
  $$
and the normalised Casoratian is then identical to the Lommel polynomial,
  $$
  \ol W(\nu,p)=R^{\nu,p}(z)\,,
  $$
for non--integral $\nu$. This extends to integral $\nu$ by continuity since $R^{n,p}$ is
{\it uniquely} defined (\cf\ Watson, [\pref{Watson2}] \S9.61).

The classic result, (\peq{wronsk3}), is a consequence of the general (\peq{Leq}) for
$p=0$ together with the explicit form of the $R^{\nu,p}$, although it can be
established independently.

In order to regain the notation of (\peq{lin8}), $\nu$ must be set equal to $A/B$, and
$z$ to $2/B$ and (\peq{Lommp}) yields the most explicit formula (\peq{wronsksum3}),
as in section 5. (See below.)

Actually, the passage to integral $\nu$ is best approached using the second kind
functions, $Y_\nu$. Simple algebra using (\peq{2sol}) shows that (\peq{Leq}) is
replaced by the more elegant,
  $$
  {2\over\pi\,z}
 \, R^{\nu,p}(z)=Y_{\nu}(z)\,J_{\nu+p+1}(z)-J_\nu(z)\,Y_{\nu+p+1}(z)=\bigg|\matrix{Y_{\nu}&Y_{\nu+p+1}
  \cr J_\nu&J_{\nu+p+1}}\bigg|\,,
 \eql{Leq2}
  $$
also due to Lommel, [\pref{Nielsen}]. This will play a role later in the continuum limit.

\section{\bf 8. Lommel polynomials as orthogonal polynomials}

First a few words regarding variables are relevant here. The variables that provide the
standard form of the Lommel polynomials, \ie\ equation (\peq{Lommp}), are $z$ and
$\nu$ related to the `physical' variables $A=2-\la$ and $B$ where $\la$ is the
eigenvalue and $B$ the strength of the (linear) potential, by $\nu=A/B$ and $z=2/B$,
as stated earlier, so that
$$
  \nu={2-\la\over B}={z\over2}(2-\la)\,\quad
$$

In terms of  $A$ and $B$ separately, (\peq{wronsksum3}) provides the expression for
Lommel polynomials. It shows that they are polynomials in both $A$ and $B$ or
extracting factors of $A$ or of $B$, polynomials in $B/A$ or $A/B=\nu$. Here I am
interested in polynomials in $A=2-\la$, or equivalently for fixed $B$, in $\nu$.

At  this point, it is convenient to list a few Lommel polynomials, \footnote{ There is a
misprint on p.102 in [\pref{GandG}]. I have temporarily changed notations so that
$R(\nu,p)\equiv R^{\nu,p}(z)$.}
  $$\eqalign{
  R(\nu,1)&=(\nu+1)\bigg({2\over z}\bigg)\cr
  R(\nu,2)&=(\nu+2)(\nu+1)\bigg({2\over z}\bigg)^2-1\cr
  R(\nu,3)&=(\nu+3)(\nu+2)(\nu+1)\bigg({2\over z}\bigg)^3-2(\nu+2)\bigg({2\over z}\bigg)\cr
  R(\nu,4)&=(\nu+4)(\nu+3)(\nu+2)(\nu+1)\bigg({2\over z}\bigg)^4-
  3(\nu+3)(\nu+2)\bigg({2\over z}\bigg)^2
+1\cr
  }
  \eql{lommelps}
  $$
together with the associated recursion (due to Lommel, [\pref{Lommel}])
[\pref{Nielsen}], p.54, [\pref{GandG}] p.102,
  $$
  R(\nu,p+2)={2(p+\nu)\over z}\,R(\nu,p+1)-R(\nu,p)\,.
  \eql{lomrec}
  $$

By classic, general theory, the $y(j,\la)$ of (\peq{BMsol}), or equivalently the $\ol W$ of
(\peq{wronsksum2}) (with $A=2-\la$ and $p+1$ reset to $j$), \ie the Lommel
polynomials, $R$, form a finite set of explicit orthogonal polynomials (in $\la$, or $\nu$)
concentrated on the eigenvalues, $\la_n$ or $\nu_n$, \eg\ Atkinson, [\pref{Atkinson}].

For $p$ vertices, orthogonality reads
  $$
   \sum_{j=1}^p R(\nu,j-1)\,R(\nu',j-1) =1+\sum_{j=2}^p R(\nu,j-1)\,R(\nu',j-1)
   =N(\nu)\,\de_{\nu\nu'}
   \eql{orth}
  $$
where $\nu$ and $\nu'$ are eigenvalues determined by the vanishing of $R(\nu,p)$.

It is interesting to confirm this relation. As a numerical example, take $p=2$. The
eigenvalue equation is then
  $$
  (\nu+1)(\nu+2)-{1\over4}z^2=0
  \eql{eig}
  $$
and (\peq{orth}) reads
  $$
  1+{4\over z^2}(\nu+1)(\nu'+1)=0\,,\quad \nu\ne\nu'
  \eql{orth2}
  $$
and
  $$
  1+{z^2\over4}(\nu+1)^2=N(\nu)\,.
  \eql{orth3}
  $$
Using $(\nu+1)$ as the variable, the product of eigenvalues is, from (\peq{eig}),
$-z^2/4$ and (\peq{orth2}) is readily verified. (\peq{orth3}) fixes the normalisation for
orthonormality.

More complicated is the  case of $p=3$.  Orthogonality is now
   $$
   1+{4\over z^2}(\nu+1)(\nu'+1)+\bigg({4\over z^2}(\nu+2)(\nu+1)-1\bigg)
  \bigg({4\over z^2}\nu'+2)(\nu'+1)-1\bigg)=0
  \eql{orth4}
   $$

From (\peq{lommelps}) the eigenvalue equation $R(\nu,3)=0$ has a solution $\nu=-2$,
as mentioned before, and, for simplicity, I set one eigenvalue, $\nu'=-2$ so that
(\peq{orth4}) reads,
   $$
   1-{4\over z^2}(\nu+1)+1-{4\over z^2}\,(\nu+2)(\nu+1)=0
   $$
which reduces to
  $$
  (\nu+3)(\nu+1){2\over z^2}-1=0\,,
  \eql{otherrts}
  $$
\ie the equation determining the other two roots of $R(\nu,3)=0$ and the confirmation is
complete. A similar mechanism works for any odd $p$.

Showing that (\peq{orth4}) holds for solutions of (\peq{otherrts}) is more complicated.

\section{\bf9. The continuum limit}

In this section I show how to derive the continuum determinant ratio, (\peq{airy2}),
from the above discrete expressions.

I denote by $h$ the path step, $h=1/(p+1)$, and seek to let $h\to0$, particularly in
(\peq{Leq2}). As noted earlier, the strength, $B$, of the discrete linear potential is
related to $b^3$, the strength of the continuum linear potential by
$B=\big(b/(p+1)\big)^3=(bh)^3$. Then the argument, $z=2/(bh)^3\to \infty$.

Since I am interested only in the determinant,  I can immediately set $\la=0$ so that
$\nu=z $.  Therefore as $h\to 0$ the {\it leading divergence} of the the orders of the
Bessel functions in (\peq{Leq2}) is just that of their arguments and I can employ
Nicholson's asymptotic relations. Starting from
  $$
  {2\over\pi\,z}
 \, R^{z,p}(z)=Y_{z}(z)\,J_{z+p+1}(z)-J_z(z)\,Y_{z+p+1}(z)
 \eql{Leq3}
  $$
I need
     $$\eqalign{
  Y_{z+\al}(z)&\sim 3^{-1/6}\bigg({\xi\over z}\bigg)^{1/3}[I_{1/3}(\xi)+I_{-1/3}(\xi)]\cr
  &=-\bigg({3\over x}\bigg)^{1/2} 3^{-1/6}\bigg({\xi\over z}\bigg)^{1/3}\,\Bi(x)
  }
  $$
and
  $$\eqalign{
  J_{z+\al}(z)&\sim \pi^{-1}3^{-1/6}\bigg({\xi\over z}\bigg)^{1/3}\,K_{1/3}(\xi)\cr
  &=\bigg({3\over x}\bigg)^{1/2}\,3^{-1/6}\bigg({\xi\over z}\bigg)^{1/3}\,\Ai(x)
  }
  $$
where $\Ai$ and $\Bi$ are Airy functions and

  $$
 \xi={2\over3}\bigg({2\over z}\bigg)^{1/2}\al^{3/2}={2\over3}\,x^{3/2}\,.
  $$

The parameter $\al$ is either 0 or $(p+1)=1/h$ and $2/z=(hb)^3$  from which it is
seen that the arguments, $x$, of the Airy functions are either 0 or $b$.

The multiplying factor reads
  $$
  \bigg({3\over x}\bigg)^{1/2}\,3^{-1/6}\bigg({\xi\over z}\bigg)^{1/3}=
  3^{1/2}\,3^{-1/6}\bigg({1\over z}\bigg)^{1/3}\bigg({2\over3}\bigg)^{1/3}=
  \bigg({2\over z}\bigg)^{1/3}\,,
  $$
and there are two of these.

Hence we have the transitional limiting behaviour of the Lommel polynomials
  $$\eqalign{
   R^{z,p}(z)&\sim \pi\bigg({z\over2}\bigg)^{1/3}\big(\Bi(0)\,\Ai(b)-\Ai(0)\,\Bi(b)\big)\cr
   &\sim {1\over h}\,{\pi\over  b}\big(\Bi(0)\,\Ai(b)-\Ai(0)\,\Bi(b)\big)\,.
   }
   \eql{asymp}
  $$

The form (\peq{asymp}) contains within it Cauchy's limit,
  $$
    J_\nu(\nu)\sim {\Ga(1/3)\over 2^{2/3}\,3^{1/6}\,\pi^{1/3}\nu^{1/3}}\,,
  $$
and also
  $$
   Y_\nu(\nu)\sim-{3^{1/3}\,\Ga(1/3)\over2^{2/3}\,\pi\nu^{1/3}}\,.
  $$

Numerically, choosing $p=9$, a not too big a number, \ie $h=1/10$, the exact
evaluation of $R(2000,9,2000)$ from the series is $10.84393086$, for $b=1$, whereas
the limiting approximation, (\peq{asymp}), gives $10.85339648$. (These are just the
same as the numbers in section (3).)

As I have not been too careful with normalisation, this quantity, (\peq{asymp}), will be
proportional to the functional determinant. This does not matter as I only want  the {\it
ratio} of the determinants, with and without $b$. As $b\to 0$, it is trivial to show that
the right--hand side of (\peq{asymp}) is just $1/(hb)$ and the desired ratio agrees with
that, (\peq{airy2}), obtained in the continuum directly by Kirsten and McKane,
[\pref{KandM}].

\section{10. \bf The discrete zero mode problem}

If one of the eigenvalues, say $\la_0$, vanishes, then the ordinary determinant is zero. It
is then customary to calculate the product of the remaining non--zero eigenvalues and
refer to this as the determinant. Kirsten and McKane, [\pref{KandM}], give a
modification of the Gel'fand Yaglom procedure that covers this case. Their formula just
involves the zero mode eigenfunction. In this section I give the discrete version.

The result drops out  of the discrete analogue of Green's theorem, or Lagrange's
identity,\footnote{ This is the method employed in [\pref{KandM}].}, or the
Christoffel--Darboux identity, \eg\ Atkinson, [\pref{Atkinson}] \S4.2. For speed, I copy
the theorem from [\pref{Atkinson}], for the recursion (\peq{lin3})
  $$
  (\la-\mu)\sum_{j=1}^k y(j,\la)\,y(j,\mu)=\bigg|\matrix{y(k+1,\la)&y(k+1,\mu)\cr
  y(k,\la)&y(k,\mu)}\bigg|\,,\quad 0\le k<p+1\,,
  \eql{ChDa}
  $$
where $\la$ and $\mu$ are any two numbers. Now choose $\mu$ to be an eigenvalue,
$\la_i$,  and set the upper limit $k=p$, the number of vertices. Then $y(p+1,\la)$ is the
characteristic polynomial and so $y(p+1,\mu)=y(p+1,\la_i)=0$. Equation (\peq{ChDa})
is rearranged into
  $$
  {y(p+1,\la)\over\la-\la_i}={1\over y(p,\la_i)}\sum_{j=1}^p y(j,\la)\,y(j,\la_i)\,.
  \eql{ChDa2}
  $$
Setting $\la$ to $\la_i$ gives
    $$
  \prod_{n\ne i}(\la_i-\la_n)= {1\over y(p,\la_i)}\sum_{j=1}^p y(j,\la_i)\,y(j,\la_i)
  =-{\br {y(\la_i)}{y(\la_i)}\over \De y(p,\la_i)}
  \eql{zmdet}
    $$
where I have again used the Dirichlet condition $y(p+1,\la_i)=0$. Choosing
$\la_i=\la_0=0$, the left--hand side is, up to a sign, the determinant, omitting the zero
mode. This result should be compared with the relevant part of equation (14) in
[\pref{KandM}].

As a non--trivial example I again turn to the linear potential case when the characteristic
polynomial is a Lommel polynomial. Specifically I choose 3 vertices ($p=3$) and select
the $\nu=-2$ root to give  the zero mode. This means that the strength, $B$, has to
equal $-1$ since $\nu=(2-\la)/B$. The remaining two roots are then easily found as
$\pm \sqrt3$ so that their determinant is $-3$. Computation of the relevant Lommel
polynomial, $R(-2,j-1)(-2)$, gives the zero mode as $(1,1,-1)$, which is consistent with
(\peq{zmdet}).

Incidentally, by dividing (\peq{ChDa2}) and (\peq{zmdet}) one obtains Legendre's
function and Legendre's polynomial interpolation of $f$ can therefore be written as,
  $$
  F(\la)=\sum_{n=0}^{p-1} f(\la_n) {\br{y(\la)}{y(\la_n)}\over\br{y(\la_n)}{y(\la_n)}}\,,
  \eql{sample}
  $$
in terms of a set of $p$ orthogonal polynomials. (To repeat, the eigenvalues, $\la_n$, are
the zeros of the polynomial, $y(p+1,\la)$.) An historic and practically important case
corresponds to the free equation, $V=0$, which yields Chebyshev polynomials. See
[\pref{MandH}] chap.6.

Equation (\peq{sample}) is related to a discrete Kramers sampling theorem, \cf\
Annaby, [\pref{Annaby}]. If $f(\la)$ is a linear combination of the $y(j,\la)$, ($1\le j\le
p$), then $F(\la)=f(\la)$ and $f(\la)$ can be reconstructed from its samples, which is,
perhaps, not surprising.

This can be more abstractly expressed as the vector space expansions,
  $$
  \ket f=\sum_{j=1}^p\ket j\br jf=\sum_{n=0}^{p-1} {\ket {\la_n}\br{\la_n}f\over \br{\la_n}{\la_n}}
  $$
(\peq{sample}) being regained on multiplying by $\bra \la$, the function label, $y$,
being suppressed. The condition on $f$ corresponds to the bandwidth limitation of the
original Whittaker-Shannon sampling theorem.

\section{ 11. \bf Comments}
It might be said that the Bessel functions are only intermediaries in the calculation of the
Lommel polynomials, which could be taken as the primary objects. Indeed, Bessel functions
can be obtained  as limits of Lommel polynomials, as shown by Hurwitz, [\pref{Hurwitz}].
(Incidentally, precisely this limit was discussed by Bleich and Melan, [\pref{BlandM}], who
make no reference to any earlier work.)

If this attitude is adopted, the form of the Lommel polynomials, \eg\ (\peq{Lommp}),
must be found without explicit use of Bessel functions, but using only the single
recursion (\peq{lin8}). Such is Lommel's original computation, [\pref{Lommel}],  but he
assumes the result and justifies it by induction which is somewhat synthetic. Bleich and
Melan's method, [\pref{BlandM}], is interesting, but a little complicated, and again
involves an assumption regarding the solution's general form. Graf and Gubler,
[\pref{GandG}] vol.2, however, use the recursion, expressed as a continued fraction, and
obtain the Lommel polynomials, which they refer to as Schl\"afli functions, directly. They
list a few examples and also give the factored form (\peq{wronsksum3}).
\footnote{Although Watson, [\pref{Watson2}], is the undoubted Bessel bible, and
Erdelyi, [\pref{EMOT2}], the best summary of results, a number of earlier, and later,
works, provide more detail and alternative approaches. I mention the rather attractive
book by Nielsen, [\pref{Nielsen}], which furnishes a more expansive treatment of the
basics and contains results not available in [\pref{Watson2}]. A little referenced book is
that by Graf and Gubler, [\pref{GandG}],  which presents the Schl\"afli, complex
integral, formulation of the theory. It is useful for the amount of algebraic detail given
and also for Lommel polynomial properties. The pioneer work of Lommel,
[\pref{Lommel}], is also valuable.}

It is possible to derive any particular Lommel polynomial by direct iteration of the
recursion (\peq{bessrec}), or (\peq{lin8}). The reduction formula, [\pref{Nielsen}],
\S26,
  $$
  R^{\,\nu,p}(z)=\cos{p\pi\over2}+{2\over z}
  \sum_{s=0}^{p-1}(\nu+s+1)\sin{p-s\over2}\pi\,R^{\,\nu,s}
  $$
can also be iterated.

\section{11. \bf Conclusion}

It has been shown that. apart from providing quite accurate approximations for smallish
numbers of lattice points to a continuum situation, discrete techniques have an intrinsic
interest and form a useful intro to the theory of orthogonal polynomials. In particular, I
have highlighted the case of Lommel polynomials which seem to have been ignored. Even
their general orthogonality properties do not seem to be available explicitly.
 \vglue 20truept

\noin{\bf References.} \vskip5truept

\begin{putreferences}
  \ref{Annaby}{Annaby,M.H. {\it Analysis} {\bf 18} (1998) 55.}
  \ref{Watson2}{Watson,G.N. {\it A Treatise on the theory of Bessel Functions} (CUP, Cambridge, 1944).}
  \ref{Dowgelyag}{Dowker,J.S. {\it Discrete Determinants and the Gel'fand--Yaglom
  formula}, ArXiv. }
  \ref{Flugge}{Fl\"ugge,S.   {\it Practical Quantum Mechanics} (Springer,Berlin,1971).}
  \ref{GandG}{Graf,J.H. and Gubler,E. {\it Einleitung in die Theorie der Bessel'schen
  Funktionen} (Wyss, Bern,1898-1900).}
  \ref{Lommel}{Lommel, E. {\it Studien \"uber die Bessel'schen Functionen} (Teubner, Leipzig,
  1868).}
  \ref{Nielsen}{Nielsen,N. {\it Handbuch der Theorie der Cylinderfunktionen} (Teubner, Leipzig,
  1904.).}
  \ref{LandT}{Lakshmikantham,V. and Trigiante,D. {\it Theory of Difference Equations}
  (Academic Press, San Diego, 1988).}
  \ref{Agarwal}{Agarwal,R.P. {\it Difference Equations and Inequalities} (Dekker, New York,
  1992).}
  \ref{Kleinert}{Kleinert, H. {\it Path Integrals in Quantum mechanics} (World Scientific,
   Singapore, 2004).}
  \ref{Miller2}{Miller,K.S. {\it Linear Difference Equations} (Benjamin, New York, 1968).}
  \ref{GrandM}{Gray,A. and Matthews, G.B. {\it A treatise on Bessel functions}, (MacMillan,
  London, 1922).}
  \ref{Birkhoff}{Birkhoff,G.D. \tams{12}{1911}{243}.}
  \ref{LeandL}{Levy,H. and Lessman,F. {\it Finite Difference Equations} (Pitman, London,
  1959).}
  \ref{EandM}{Erhardt and Mickens .}
  \ref{Bohmer}{B\"ohmer,E. {\it Differenzengleichungen und bestimmte Integrale} (Koehler,
  Leipzig, 1939).}
  \ref{Barnesrec}{Barnes,E.W. {\it Mess. Math.} {\bf34}(1904) 52.}
   \ref{Bass}{Bass,R. \jmp{26}{1985}{3068}.}
   \ref{GandK}{Gantmacher,F.R. and Krein,M.G. {\it Oszillazionsmatrizen}
   (Akad.-Verlag.Berlin, 1960).}
   \ref{Hald}{Hald,O.H. {\it Numer.Math.} {\bf27} (1977) 249.}
   \ref{CandY}{Chung,F. and Yau S.--T. {\it J.Comb.Theory} {\bf A91} (2000) 191.}
   \ref{BEC}{Bendito,E, Encinas,A.M. and Carmona,A. {\it Appl.Anal.Disc.Math.}{\bf 3}
   (2009) 282.}
   \ref{AandA}{Annaby,M.H. and Asharabi,R.M. {\it Acta.Math.Scientia} {\bf 31B} (2011)408.}
   \ref{MandH}{Mason,J.C. and Handscomb,D.C. {\it Chebyshev Polynomials} (Chapman and Hall,
   Boca Raton, 2002).}
   \ref{LandB}{Levy,H. and Baggott,E.A. \phm{18}{1934}{177}.}
   \ref{BlandM}{Bleich,Fr. and Melan,E. {\it Die Gew\"onlichen und Partiellen
   Differenzengleichungen der Baustatik}, (Springer, Berlin, 1927).}
   \ref{Spiegel}{Spiegel,M.R. {\it Schaum's Outline of Calculus of Finite Differences},
   (McGraw-Hill, New York, 1971).}
   \ref{Porter}{Porter,M.B. \aom{3}{1901}{55}.}
    \ref{Dunne}{Dunne,G. \jpa{41}{2008}{304006}.}
    \ref{MandF}{Morse,P.M. and  Feshbach,H. {\it Methods of Theoretical Physics},
    (McGraw--Hill, New York, 1953).}
    \ref{KandM}{Kirsten, K. and McKane,A. \aop{308}{2003}{502}.}
    \ref{LandS}{Levit,S. and Smilansky,U. \pams{65}{1977}{299}.}
    \ref{GeandY}{Gel'fand, I.M. and Yaglom,A.M. \jmp {1}{1960}{48}.}
    \ref{ABR}{Actor,A., Bender,C. and Reingruber,J., {\it Fortschr.Phys.} {\bf48} (2000) 4.}
    \ref{Rayleigh}{Rayleigh, Lord, {\it Theory of Sound}, (MacMillan, London, 1894).}
    \ref{Floratos}{Floratos,E.G. \plb{228}{1989}{335}.}
    \ref{deverdiere}{de Verdi\`ere, C. {\it Ann. Inst. Fourier} {\bf 49} (1999) 861.}
    \ref{Goldberg}{Goldberg,S. {\it Introduction to Difference Equations},
    (Wiley, New York, 1958).}
    \ref{Jordan}{Jordan,C. {\it Calculus of Finite Differences}, (Budapest, 1939).}
    \ref{Fort}{Fort,T. {\it Finite Differences}, (Clarendon Press, Oxford,1948). }
    \ref{Atkinson}{Atkinson,F.V. {\it Discrete and Continuous Boundary Problems},
    (Academic Press, New York,1964).}
   \ref{Elaydi}{Elaydi,S.N. {\it An Introduction to Difference Equations}, (Springer, New York,
   1999.}
   \ref{Forman}{Forman,R. \cmp{147}{1992}{485}.}
    \ref{Boole}{Boole, G. {\it Calculus of Finite Differences}, (MacMillan, Cambridge,
1860).}
     \ref{Cayley5}{Cayley,A.  {\it Phil. Trans. Roy. Soc. Lond.} {\bf 148} (1858) 47.}
     \ref{Milne-Thomson}{Milne-Thomson, L.M. {\it The Calculus of Finite Differences},
     (MacMillan, London, 1933).}
    \ref{Herschel}{Herschel, J.F.W.  {\it Phil. Trans. Roy. Soc. Lond.} {\bf 106} (1816) 25.}
    \ref{Littlewood2}{Littlewood,D.E. {\it The Theory of Group Characters}
    (Clarendon Press, Oxford, 1950).}
    \ref{Wright}{Wright,E.M. \amm{68}{1961}{144}.}
\ref{Carlitz}{Carlitz,L. \dmj{27}{1960}{401}.}
     \ref{Netto}{Netto,E. {\it Lehrbuch der Combinatorik} 2nd Edn. (Teubner, Leipzig, 1927).}
    \ref{FdeB}{Fa\`{a} de Bruno, F. {\it Th\'eorie des Formes Binaires} (Brero, Turin,1876).}
    \ref{Ehrhart}{Ehrhart,E. \jram{227}{25}{1967}.}
    \ref{Bell}{Bell,E.T.\ajm{65}{1943}{382}.}
    \ref{BandR}{Beck,M. and Robins,S. {\it Computing the Continuous Discretely,}
    (Springer, New York, 2007).}
    \ref{BandR2}{Beck, M. and Robins,S. {\it Discrete and Comp. Geom.} {\bf 27}(2002) 443.}
    \ref{Harmer}{Harmer,M. {\it J.Australian Math.Soc.} {\bf 84}(2008)217.}
    \ref{RandF}{Rubinstein, B.Y. and Fel,L.G., {\it Ramanujan J.} {\bf11}(2006)331.}
    \ref{BGK}{Beck,M., Gessel, I.M. and Komatsu,T. {\it Electronic Journal of Combinatorics}
    {\bf8}(2001) 1.}
    \ref{Sylvester}{Sylvester,J.J. \qjpam{1}{1858}{81}.}
    \ref{Sylvester2}{Sylvester,J.J. \qjpam{1}{1858}{142}.}
    \ref{Sylvester3}{Sylvester,J.J. \ajm{5}{1882}{79}.}
    \ref{Sylvester4}{Sylvester,J.J. \plms{28}{1896}{33}.}
    \ref{Dowgta}{J.S.Dowker, {\it Group theory aspects of spectral problems on spherical factors},
   ArXiv.Math.DG: 0907.1309.}
    \ref{BDR}{Beck,M., Diaz and Robins,S. {\it J.Numb.Theory} {\bf 96} (2002) 1.}
    \ref{PandS}{P\'{o}lya, G. and Szeg\H{o},G. {\it Aufgaben und Lehrs\"atze aus der Analysis}
    (Springer--Verlag, Berlin, 1925).}
    \ref{EOS}{Elizalde,E., Odintsov, S.D. and Saharian, A.A. \prD{79}{2009}{065023}.}
    \ref{Cavalcanti}{Cavalcanti,R.M. \prD{69}{2004}{065015}.}
    \ref{MWK}{Milton, K.A., Wagner,J. and Kirsten,K. \prD{80}{2009}{125028}.}
    \ref{EBM2}{Ellingsen,S.A., Brevik,I. and Milton,K.A. \prE{81}{2010}{065031}.}
    \ref{EBM}{Ellingsen,S.A., Brevik,I. and Milton,K.A. \prE{80}{2009}{021125}.}
    \ref{BEM}{Brevik,I., Ellingsen,S.A. and Milton,K.A. \prE{79}{2009}{041120}.}
    \ref{FKW}{Fulling,S.A, Kaplan L. and Wilson,J.H. \prA{76}{2007}{012118}.}
    \ref{Lukosz}{Lukosz,W, {\it Physica} {\bf 56} (1971) 109; \zfp{258}{1973}{99}
    ;\zfp{262}{1973}{327}.}
    \ref{Gromes}{Gromes, D. \mz{94}{1966}{110}.}
    \ref{FandK1}{Kirsten,K. and Fulling,S.A. \prD{79}{2009}{065019} .}
    \ref{FandK2}{Fucci,G. and Kirsten,K, JHEP (2011), 1103:016.}
    \ref{dowgjms}{Dowker,J.S. {\it Determinants and conformal anomalies
    of GJMS operators on spheres}, ArXiv: 1007.3865.}
    \ref{Dowcascone}{dowker,J.S. \prD{36}{1987}{3095}.}
    \ref{Dowcos}{dowker,J.S. \prD{36}{1987}{3742}.}
    \ref{Dowtherm}{Dowker,J.S. \prD{18}{1978}{1856}.}
    \ref{Dowgeo}{Dowker,J.S. \cqg{11}{1994}{L55}.}
    \ref{ApandD2}{Dowker,J.S. and Apps,J.S. \cqg{12}{1995}{1363}.}
   \ref{HandW}{Hertzberg,M.P. and Wilczek,F. {\it Some calculable contributions to
   Entanglement Entropy}, ArXiv:1007.0993.}
   \ref{KandB}{Kamela,M. and Burgess,C.P. \cjp{77}{1999}{85}.}
   \ref{Dowhyp}{Dowker,J.S. \jpa{43}{2010}{445402}; ArXiv:1007.3865.}
   \ref{LNST}{Lohmayer,R., Neuberger,H, Schwimmer,A. and Theisen,S.
   \plb{685}{2010}{222}.}
   \ref{Allen2}{Allen,B. PhD Thesis, University of Cambridge, 1984.}
   \ref{MyandS}{Myers,R.C. and Sinha,A. {\it Seeing a c-theorem with
   holography}, ArXiv:1006.1263}
   \ref{MyandS2}{Myers,R.C. and Sinha,A. {\it Holographic c-theorems in
   arbitrary dimensions},\break ArXiv: 1011.5819.}
   \ref{RyandT}{Ryu,S. and Takayanagi,T. JHEP {\bf 0608}(2006)045.}
   \ref{CaandH}{Casini,H. and Huerta,M. {\it Entanglement entropy
   for the n--sphere},\break arXiv:1007.1813.}
   \ref{CaandH3}{Casini,H. and Huerta,M. \jpa {42}{2009}{504007}.}
   \ref{Solodukhin}{Solodukhin,S.N. \plb{665}{2008}{305}.}
   \ref{Solodukhin2}{Solodukhin,S.N. \plb{693}{2010}{605}.}
   \ref{CaandW}{Callan,C.G. and Wilczek,F. \plb{333}{1994}{55}.}
   \ref{FandS1}{Fursaev,D.V. and Solodukhin,S.N. \plb{365}{1996}{51}.}
   \ref{FandS2}{Fursaev,D.V. and Solodukhin,S.N. \prD{52}{1995}{2133}.}
   \ref{Fursaev}{Fursaev,D.V. \plb{334}{1994}{53}.}
   \ref{Donnelly2}{Donnelly,H. \ma{224}{1976}{161}.}
   \ref{ApandD}{Apps,J.S. and Dowker,J.S. \cqg{15}{1998}{1121}.}
   \ref{FandM}{Fursaev,D.V. and Miele,G. \prD{49}{1994}{987}.}
   \ref{Dowker2}{Dowker,J.S.\cqg{11}{1994}{L137}.}
   \ref{Dowker1}{Dowker,J.S.\prD{50}{1994}{6369}.}
   \ref{FNT}{Fujita,M.,Nishioka,T. and Takayanagi,T. JHEP {\bf 0809}
   (2008) 016.}
   \ref{Hund}{Hund,F. \zfp{51}{1928}{1}.}
   \ref{Elert}{Elert,W. \zfp {51}{1928}{8}.}
   \ref{Poole2}{Poole,E.G.C. \qjm{3}{1932}{183}.}
   \ref{Bellon}{Bellon,M.P. {\it On the icosahedron: from two to three
   dimensions}, arXiv:0705.3241.}
   \ref{Bellon2}{Bellon,M.P. \cqg{23}{2006}{7029}.}
   \ref{McLellan}{McLellan,A,G. \jpc{7}{1974}{3326}.}
   \ref{Boiteaux}{Boiteaux, M. \jmp{23}{1982}{1311}.}
   \ref{HHandK}{Hage Hassan,M. and Kibler,M. {\it On Hurwitz
   transformations} in {Le probl\`eme de factorisation de Hurwitz}, Eds.,
   A.Ronveaux and D.Lambert (Fac.Univ.N.D. de la Paix, Namur, 1991),
   pp.1-29.}
   \ref{Weeks2}{Weeks,Jeffrey \cqg{23}{2006}{6971}.}
   \ref{LandW}{Lachi\`eze-Rey,M. and Weeks,Jeffrey, {\it Orbifold construction of
   the modes on the Poincar\'e dodecahedral space}, arXiv:0801.4232.}
   \ref{Cayley4}{Cayley,A. \qjpam{58}{1879}{280}.}
   \ref{JMS}{Jari\'c,M.V., Michel,L. and Sharp,R.T. {\it J.Physique}
   {\bf 45} (1984) 1. }
   \ref{AandB}{Altmann,S.L. and Bradley,C.J.  {\it Phil. Trans. Roy. Soc. Lond.}
   {\bf 255} (1963) 199.}
   \ref{CandP}{Cummins,C.J. and Patera,J. \jmp{29}{1988}{1736}.}
   \ref{Sloane}{Sloane,N.J.A. \amm{84}{1977}{82}.}
   \ref{Gordan2}{Gordan,P. \ma{12}{1877}{147}.}
   \ref{DandSh}{Desmier,P.E. and Sharp,R.T. \jmp{20}{1979}{74}.}
   \ref{Kramer}{Kramer,P., \jpa{38}{2005}{3517}.}
   \ref{Klein2}{Klein, F.\ma{9}{1875}{183}.}
   \ref{Hodgkinson}{Hodgkinson,J. \jlms{10}{1935}{221}.}
   \ref{ZandD}{Zheng,Y. and Doerschuk, P.C. {\it Acta Cryst.} {\bf A52}
   (1996) 221.}
   \ref{EPM}{Elcoro,L., Perez--Mato,J.M. and Madariaga,G.
   {\it Acta Cryst.} {\bf A50} (1994) 182.}
    \ref{PSW2}{Prandl,W., Schiebel,P. and Wulf,K.
   {\it Acta Cryst.} {\bf A52} (1999) 171.}
    \ref{FCD}{Fan,P--D., Chen,J--Q. and Draayer,J.P.
   {\it Acta Cryst.} {\bf A55} (1999) 871.}
   \ref{FCD2}{Fan,P--D., Chen,J--Q. and Draayer,J.P.
   {\it Acta Cryst.} {\bf A55} (1999) 1049.}
   \ref{Honl}{H\"onl,H. \zfp{89}{1934}{244}.}
   \ref{PSW}{Patera,J., Sharp,R.T. and Winternitz,P. \jmp{19}{1978}{2362}.}
   \ref{LandH}{Lohe,M.A. and Hurst,C.A. \jmp{12}{1971}{1882}.}
   \ref{RandSA}{Ronveaux,A. and Saint-Aubin,Y. \jmp{24}{1983}{1037}.}
   \ref{JandDeV}{Jonker,J.E. and De Vries,E. \npa{105}{1967}{621}.}
   \ref{Rowe}{Rowe, E.G.Peter. \jmp{19}{1978}{1962}.}
   \ref{KNR}{Kibler,M., N\'egadi,T. and Ronveaux,A. {\it The Kustaanheimo-Stiefel
   transformation and certain special functions} \lnm{1171}{1985}{497}.}
   \ref{GLP}{Gilkey,P.B., Leahy,J.V. and Park,J-H, \jpa{29}{1996}{5645}.}
   \ref{Kohler}{K\"ohler,K.: Equivariant Reidemeister torsion on
   symmetric spaces. Math.Ann. {\bf 307}, 57-69 (1997)}
   \ref{Kohler2}{K\"ohler,K.: Equivariant analytic torsion on ${\bf P^nC}$.
   Math.Ann.{\bf 297}, 553-565 (1993) }
   \ref{Kohler3}{K\"ohler,K.: Holomorphic analytic torsion on Hermitian
   symmetric spaces. J.Reine Angew.Math. {\bf 460}, 93-116 (1995)}
   \ref{Zagierzf}{Zagier,D. {\it Zetafunktionen und Quadratische
   K\"orper}, (Springer--Verlag, Berlin, 1981).}
   \ref{Stek}{Stekholschkik,R. {\it Notes on Coxeter transformations and the McKay
   correspondence.} (Springer, Berlin, 2008).}
   \ref{Pesce}{Pesce,H. \cmh {71}{1996}{243}.}
   \ref{Pesce2}{Pesce,H. {\it Contemp. Math} {\bf 173} (1994) 231.}
   \ref{Sutton}{Sutton,C.J. {\it Equivariant isospectrality
   and isospectral deformations on spherical orbifolds}, ArXiv:math/0608567.}
   \ref{Sunada}{Sunada,T. \aom{121}{1985}{169}.}
   \ref{GoandM}{Gornet,R, and McGowan,J. {\it J.Comp. and Math.}
   {\bf 9} (2006) 270.}
   \ref{Suter}{Suter,R. {\it Manusc.Math.} {\bf 122} (2007) 1-21.}
   \ref{Lomont}{Lomont,J.S. {\it Applications of finite groups} (Academic
   Press, New York, 1959).}
   \ref{DandCh2}{Dowker,J.S. and Chang,Peter {\it Analytic torsion on
   spherical factors and tessellations}, arXiv:math.DG/0904.0744 .}
   \ref{Mackey}{Mackey,G. {\it Induced representations}
   (Benjamin, New York, 1968).}
   \ref{Koca}{Koca, {\it Turkish J.Physics}.}
   \ref{Brylinski}{Brylinski, J-L., {\it A correspondence dual to McKay's}
    ArXiv alg-geom/9612003.}
   \ref{Rossmann}{Rossman,W. {\it McKay's correspondence
   and characters of finite subgroups of\break SU(2)} {\it Progress in Math.}
      Birkhauser  (to appear) .}
   \ref{JandL}{James, G. and Liebeck, M. {\it Representations and
   characters of groups} (CUP, Cambridge, 2001).}
   \ref{IandR}{Ito,Y. and Reid,M. {\it The Mckay correspondence for finite
   subgroups of SL(3,C)} Higher dimensional varieties, (Trento 1994),
   221-240, (Berlin, de Gruyter 1996).}
   \ref{BandF}{Bauer,W. and Furutani, K. {\it J.Geom. and Phys.} {\bf
   58} (2008) 64.}
   \ref{Luck}{L\"uck,W. \jdg{37}{1993}{263}.}
   \ref{LandR}{Lott,J. and Rothenberg,M. \jdg{34}{1991}{431}.}
   \ref{DoandKi} {Dowker.J.S. and Kirsten, K. {\it Analysis and Appl.}
   {\bf 3} (2005) 45.}
   \ref{dowtess1}{Dowker,J.S. \cqg{23}{2006}{1}.}
   \ref{dowtess2}{Dowker,J.S. {\it J.Geom. and Phys.} {\bf 57} (2007) 1505.}
   \ref{MHS}{De Melo,T., Hartmann,L. and Spreafico,M. {\it Reidemeister
   Torsion and analytic torsion of discs}, ArXiv:0811.3196.}
   \ref{Vertman}{Vertman, B. {\it Analytic Torsion of a  bounded
   generalized cone}, ArXiv:0808.0449.}
   \ref{WandY} {Weng,L. and You,Y., {\it Int.J. of Math.}{\bf 7} (1996)
   109.}
   \ref{ScandT}{Schwartz, A.S. and Tyupkin,Yu.S. \np{242}{1984}{436}.}
   \ref{AAR}{Andrews, G.E., Askey,R. and Roy,R. {\it Special functions}
   (CUP, Cambridge, 1999).}
   \ref{Tsuchiya}{Tsuchiya, N.: R-torsion and analytic torsion for spherical
   Clifford-Klein manifolds.: J. Fac.Sci., Tokyo Univ. Sect.1 A, Math.
   {\bf 23}, 289-295 (1976).}
   \ref{Tsuchiya2}{Tsuchiya, N. J. Fac.Sci., Tokyo Univ. Sect.1 A, Math.
   {\bf 23}, 289-295 (1976).}
  \ref{Lerch}{Lerch,M. \am{11}{1887}{19}.}
  \ref{Lerch2}{Lerch,M. \am{29}{1905}{333}.}
  \ref{TandS}{Threlfall, W. and Seifert, H. \ma{104}{1930}{1}.}
  \ref{RandS}{Ray, D.B., and Singer, I. \aim{7}{1971}{145}.}
  \ref{RandS2}{Ray, D.B., and Singer, I. {\it Proc.Symp.Pure Math.}
  {\bf 23} (1973) 167.}
  \ref{Jensen}{Jensen,J.L.W.V. \aom{17}{1915-1916}{124}.}
  \ref{Rosenberg}{Rosenberg, S. {\it The Laplacian on a Riemannian Manifold}
  (CUP, Cambridge, 1997).}
  \ref{Nando2}{Nash, C. and O'Connor, D-J. {\it Int.J.Mod.Phys.}
  {\bf A10} (1995) 1779.}
  \ref{Fock}{Fock,V. \zfp{98}{1935}{145}.}
  \ref{Levy}{Levy,M. \prs {204}{1950}{145}.}
  \ref{Schwinger2}{Schwinger,J. \jmp{5}{1964}{1606}.}
  \ref{Muller}{M\"uller, \lnm{}{}{}.}
  \ref{VMK}{Varshalovich.}
  \ref{DandWo}{Dowker,J.S. and Wolski, A. \prA{46}{1992}{6417}.}
  \ref{Zeitlin1}{Zeitlin,V. {\it Physica D} {\bf 49} (1991).  }
  \ref{Zeitlin0}{Zeitlin,V. {\it Nonlinear World} Ed by
   V.Baryakhtar {\it et al},  Vol.I p.717,  (World Scientific, Singapore, 1989).}
  \ref{Zeitlin2}{Zeitlin,V. \prl{93}{2004}{264501}. }
  \ref{Zeitlin3}{Zeitlin,V. \pla{339}{2005}{316}. }
  \ref{Groenewold}{Groenewold, H.J. {\it Physica} {\bf 12} (1946) 405.}
  \ref{Cohen}{Cohen, L. \jmp{7}{1966}{781}.}
  \ref{AandW}{Argawal G.S. and Wolf, E. \prD{2}{1970}{2161,2187,2206}.}
  \ref{Jantzen}{Jantzen,R.T. \jmp{19}{1978}{1163}.}
  \ref{Moses2}{Moses,H.E. \aop{42}{1967}{343}.}
  \ref{Carmeli}{Carmeli,M. \jmp{9}{1968}{1987}.}
  \ref{SHS}{Siemans,M., Hancock,J. and Siminovitch,D. {\it Solid State
  Nuclear Magnetic Resonance} {\bf 31}(2007)35.}
 \ref{Dowk}{Dowker,J.S. \prD{28}{1983}{3013}.}
 \ref{Heine}{Heine, E. {\it Handbuch der Kugelfunctionen}
  (G.Reimer, Berlin. 1878, 1881).}
  \ref{Pockels}{Pockels, F. {\it \"Uber die Differentialgleichung $\De
  u+k^2u=0$} (Teubner, Leipzig. 1891).}
  \ref{Hamermesh}{Hamermesh, M., {\it Group Theory} (Addison--Wesley,
  Reading. 1962).}
  \ref{Racah}{Racah, G. {\it Group Theory and Spectroscopy}
  (Princeton Lecture Notes, 1951). }
  \ref{Gourdin}{Gourdin, M. {\it Basics of Lie Groups} (Editions
  Fronti\'eres, Gif sur Yvette. 1982.)}
  \ref{Clifford}{Clifford, W.K. \plms{2}{1866}{116}.}
  \ref{Story2}{Story, W.E. \plms{23}{1892}{265}.}
  \ref{Story}{Story, W.E. \ma{41}{1893}{469}.}
  \ref{Poole}{Poole, E.G.C. \plms{33}{1932}{435}.}
  \ref{Dickson}{Dickson, L.E. {\it Algebraic Invariants} (Wiley, N.Y.
  1915).}
  \ref{Dickson2}{Dickson, L.E. {\it Modern Algebraic Theories}
  (Sanborn and Co., Boston. 1926).}
  \ref{Hilbert2}{Hilbert, D. {\it Theory of algebraic invariants} (C.U.P.,
  Cambridge. 1993).}
  \ref{Olver}{Olver, P.J. {\it Classical Invariant Theory} (C.U.P., Cambridge.
  1999.)}
  \ref{AST}{A\v{s}erova, R.M., Smirnov, J.F. and Tolsto\v{i}, V.N. {\it
  Teoret. Mat. Fyz.} {\bf 8} (1971) 255.}
  \ref{AandS}{A\v{s}erova, R.M., Smirnov, J.F. \np{4}{1968}{399}.}
  \ref{Shapiro}{Shapiro, J. \jmp{6}{1965}{1680}.}
  \ref{Shapiro2}{Shapiro, J.Y. \jmp{14}{1973}{1262}.}
  \ref{NandS}{Noz, M.E. and Shapiro, J.Y. \np{51}{1973}{309}.}
  \ref{Cayley2}{Cayley, A. {\it Phil. Trans. Roy. Soc. Lond.}
  {\bf 144} (1854) 244.}
  \ref{Cayley3}{Cayley, A. {\it Phil. Trans. Roy. Soc. Lond.}
  {\bf 146} (1856) 101.}
  \ref{Wigner}{Wigner, E.P. {\it Gruppentheorie} (Vieweg, Braunschweig. 1931).}
  \ref{Sharp}{Sharp, R.T. \ajop{28}{1960}{116}.}
  \ref{Laporte}{Laporte, O. {\it Z. f. Naturf.} {\bf 3a} (1948) 447.}
  \ref{Lowdin}{L\"owdin, P-O. \rmp{36}{1964}{966}.}
  \ref{Ansari}{Ansari, S.M.R. {\it Fort. d. Phys.} {\bf 15} (1967) 707.}
  \ref{SSJR}{Samal, P.K., Saha, R., Jain, P. and Ralston, J.P. {\it
  Testing Isotropy of Cosmic Microwave Background Radiation},
  astro-ph/0708.2816.}
  \ref{Lachieze}{Lachi\'eze-Rey, M. {\it Harmonic projection and
  multipole Vectors}. astro- \break ph/0409081.}
  \ref{CHS}{Copi, C.J., Huterer, D. and Starkman, G.D.
  \prD{70}{2003}{043515}.}
  \ref{Jaric}{Jari\'c, J.P. {\it Int. J. Eng. Sci.} {\bf 41} (2003) 2123.}
  \ref{RandD}{Roche, J.A. and Dowker, J.S. \jpa{1}{1968}{527}.}
  \ref{KandW}{Katz, G. and Weeks, J.R. \prD{70}{2004}{063527}.}
  \ref{Waerden}{van der Waerden, B.L. {\it Die Gruppen-theoretische
  Methode in der Quantenmechanik} (Springer, Berlin. 1932).}
  \ref{EMOT}{Erdelyi, A., Magnus, W., Oberhettinger, F. and Tricomi, F.G. {
  \it Higher Transcendental Functions} Vol.1 (McGraw-Hill, N.Y. 1953).}
   \ref{EMOT2}{Erdelyi, A., Magnus, W., Oberhettinger, F. and Tricomi, F.G. {
  \it Higher Transcendental Functions} Vol.2 (McGraw-Hill, N.Y. 1953).}
  \ref{Dowzilch}{Dowker, J.S. {\it Proc. Phys. Soc.} {\bf 91} (1967) 28.}
  \ref{DandD}{Dowker, J.S. and Dowker, Y.P. {\it Proc. Phys. Soc.}
  {\bf 87} (1966) 65.}
  \ref{DandD2}{Dowker, J.S. and Dowker, Y.P. \prs{}{}{}.}
  \ref{Dowk3}{Dowker,J.S. \cqg{7}{1990}{1241}.}
  \ref{Dowk5}{Dowker,J.S. \cqg{7}{1990}{2353}.}
  \ref{CoandH}{Courant, R. and Hilbert, D. {\it Methoden der
  Mathematischen Physik} vol.1 \break (Springer, Berlin. 1931).}
  \ref{Applequist}{Applequist, J. \jpa{22}{1989}{4303}.}
  \ref{Torruella}{Torruella, \jmp{16}{1975}{1637}.}
  \ref{Weinberg}{Weinberg, S.W. \pr{133}{1964}{B1318}.}
  \ref{Meyerw}{Meyer, W.F. {\it Apolarit\"at und rationale Curven}
  (Fues, T\"ubingen. 1883.) }
  \ref{Ostrowski}{Ostrowski, A. {\it Jahrsb. Deutsch. Math. Verein.} {\bf
  33} (1923) 245.}
  \ref{Kramers}{Kramers, H.A. {\it Grundlagen der Quantenmechanik}, (Akad.
  Verlag., Leipzig, 1938).}
  \ref{ZandZ}{Zou, W.-N. and Zheng, Q.-S. \prs{459}{2003}{527}.}
  \ref{Weeks1}{Weeks, J.R. {\it Maxwell's multipole vectors
  and the CMB}.  astro-ph/0412231.}
  \ref{Corson}{Corson, E.M. {\it Tensors, Spinors and Relativistic Wave
  Equations} (Blackie, London. 1950).}
  \ref{Rosanes}{Rosanes, J. \jram{76}{1873}{312}.}
  \ref{Salmon}{Salmon, G. {\it Lessons Introductory to the Modern Higher
  Algebra} 3rd. edn. \break (Hodges,  Dublin. 1876.)}
  \ref{Milnew}{Milne, W.P. {\it Homogeneous Coordinates} (Arnold. London. 1910).}
  \ref{Niven}{Niven, W.D. {\it Phil. Trans. Roy. Soc.} {\bf 170} (1879) 393.}
  \ref{Scott}{Scott, C.A. {\it An Introductory Account of
  Certain Modern Ideas and Methods in Plane Analytical Geometry,}
  (MacMillan, N.Y. 1896).}
  \ref{Bargmann}{Bargmann, V. \rmp{34}{1962}{300}.}
  \ref{Maxwell}{Maxwell, J.C. {\it A Treatise on Electricity and
  Magnetism} 2nd. edn. (Clarendon Press, Oxford. 1882).}
  \ref{BandL}{Biedenharn, L.C. and Louck, J.D.
  {\it Angular Momentum in Quantum Physics} (Addison-Wesley, Reading. 1981).}
  \ref{Weylqm}{Weyl, H. {\it The Theory of Groups and Quantum Mechanics}
  (Methuen, London. 1931).}
  \ref{Robson}{Robson, A. {\it An Introduction to Analytical Geometry} Vol I
  (C.U.P., Cambridge. 1940.)}
  \ref{Sommerville}{Sommerville, D.M.Y. {\it Analytical Conics} 3rd. edn.
   (Bell, London. 1933).}
  \ref{Coolidge}{Coolidge, J.L. {\it A Treatise on Algebraic Plane Curves}
  (Clarendon Press, Oxford. 1931).}
  \ref{SandK}{Semple, G. and Kneebone. G.T. {\it Algebraic Projective
  Geometry} (Clarendon Press, Oxford. 1952).}
  \ref{AandC}{Abdesselam A., and Chipalkatti, J. {\it The Higher
  Transvectants are redundant}, arXiv:0801.1533 [math.AG] 2008.}
  \ref{Elliott}{Elliott, E.B. {\it The Algebra of Quantics} 2nd edn.
  (Clarendon Press, Oxford. 1913).}
  \ref{Elliott2}{Elliott, E.B. \qjpam{48}{1917}{372}.}
  \ref{Howe}{Howe, R. \tams{313}{1989}{539}.}
  \ref{Clebsch}{Clebsch, A. \jram{60}{1862}{343}.}
  \ref{Prasad}{Prasad, G. \ma{72}{1912}{136}.}
  \ref{Dougall}{Dougall, J. \pems{32}{1913}{30}.}
  \ref{Penrose}{Penrose, R. \aop{10}{1960}{171}.}
  \ref{Penrose2}{Penrose, R. \prs{273}{1965}{171}.}
  \ref{Burnside}{Burnside, W.S. \qjm{10}{1870}{211}. }
  \ref{Lindemann}{Lindemann, F. \ma{23} {1884}{111}.}
  \ref{Backus}{Backus, G. {\it Rev. Geophys. Space Phys.} {\bf 8} (1970) 633.}
  \ref{Baerheim}{Baerheim, R. {\it Q.J. Mech. appl. Math.} {\bf 51} (1998) 73.}
  \ref{Lense}{Lense, J. {\it Kugelfunktionen} (Akad.Verlag, Leipzig. 1950).}
  \ref{Littlewood}{Littlewood, D.E. \plms{50}{1948}{349}.}
  \ref{Fierz}{Fierz, M. {\it Helv. Phys. Acta} {\bf 12} (1938) 3.}
  \ref{Williams}{Williams, D.N. {\it Lectures in Theoretical Physics} Vol. VII,
  (Univ.Colorado Press, Boulder. 1965).}
  \ref{Dennis}{Dennis, M. \jpa{37}{2004}{9487}.}
  \ref{Pirani}{Pirani, F. {\it Brandeis Lecture Notes on
  General Relativity,} edited by S. Deser and K. Ford. (Brandeis, Mass. 1964).}
  \ref{Sturm}{Sturm, R. \jram{86}{1878}{116}.}
  \ref{Schlesinger}{Schlesinger, O. \ma{22}{1883}{521}.}
  \ref{Askwith}{Askwith, E.H. {\it Analytical Geometry of the Conic
  Sections} (A.\&C. Black, London. 1908).}
  \ref{Todd}{Todd, J.A. {\it Projective and Analytical Geometry}.
  (Pitman, London. 1946).}
  \ref{Glenn}{Glenn. O.E. {\it Theory of Invariants} (Ginn \& Co, N.Y. 1915).}
  \ref{DowkandG}{Dowker, J.S. and Goldstone, M. \prs{303}{1968}{381}.}
  \ref{Turnbull}{Turnbull, H.A. {\it The Theory of Determinants,
  Matrices and Invariants} 3rd. edn. (Dover, N.Y. 1960).}
  \ref{MacMillan}{MacMillan, W.D. {\it The Theory of the Potential}
  (McGraw-Hill, N.Y. 1930).}
   \ref{Hobson}{Hobson, E.W. {\it The Theory of Spherical
   and Ellipsoidal Harmonics} (C.U.P., Cambridge. 1931).}
  \ref{Hobson1}{Hobson, E.W. \plms {24}{1892}{55}.}
  \ref{GandY}{Grace, J.H. and Young, A. {\it The Algebra of Invariants}
  (C.U.P., Cambridge, 1903).}
  \ref{FandR}{Fano, U. and Racah, G. {\it Irreducible Tensorial Sets}
  (Academic Press, N.Y. 1959).}
  \ref{TandT}{Thomson, W. and Tait, P.G. {\it Treatise on Natural Philosophy}
   (Clarendon Press, Oxford. 1867).}
  \ref{Brinkman}{Brinkman, H.C. {\it Applications of spinor invariants in
atomic physics}, North Holland, Amsterdam 1956.}
  \ref{Kramers1}{Kramers, H.A. {\it Proc. Roy. Soc. Amst.} {\bf 33} (1930) 953.}
  \ref{DandP2}{Dowker,J.S. and Pettengill,D.F. \jpa{7}{1974}{1527}}
  \ref{Dowk1}{Dowker,J.S. \jpa{}{}{45}.}
  \ref{Dowk2}{Dowker,J.S. \aop{71}{1972}{577}}
  \ref{DandA}{Dowker,J.S. and Apps, J.S. \cqg{15}{1998}{1121}.}
  \ref{Weil}{Weil,A., {\it Elliptic functions according to Eisenstein
  and Kronecker}, Springer, Berlin, 1976.}
  \ref{Ling}{Ling,C-H. {\it SIAM J.Math.Anal.} {\bf5} (1974) 551.}
  \ref{Ling2}{Ling,C-H. {\it J.Math.Anal.Appl.}(1988).}
 \ref{BMO}{Brevik,I., Milton,K.A. and Odintsov, S.D. \aop{302}{2002}{120}.}
 \ref{KandL}{Kutasov,D. and Larsen,F. {\it JHEP} 0101 (2001) 1.}
 \ref{KPS}{Klemm,D., Petkou,A.C. and Siopsis {\it Entropy
 bounds, monoticity properties and scaling in CFT's}. hep-th/0101076.}
 \ref{DandC}{Dowker,J.S. and Critchley,R. \prD{15}{1976}{1484}.}
 \ref{AandD}{Al'taie, M.B. and Dowker, J.S. \prD{18}{1978}{3557}.}
 \ref{Dow1}{Dowker,J.S. \prD{37}{1988}{558}.}
 \ref{Dowrob}{Dowker,J.S. \cqg{13}{1996}{585}.}
 \ref{Dow30}{Dowker,J.S. \prD{28}{1983}{3013}.}
 \ref{DandK}{Dowker,J.S. and Kennedy,G. \jpa{}{1978}{895}.}
 \ref{Dow2}{Dowker,J.S. \cqg{1}{1984}{359}.}
 \ref{DandKi}{Dowker,J.S. and Kirsten, K. {\it Comm. in Anal. and Geom.
 }{\bf7} (1999) 641.}
 \ref{DandKe}{Dowker,J.S. and Kennedy,G.\jpa{11}{1978}{895}.}
 \ref{Gibbons}{Gibbons,G.W. \pl{60A}{1977}{385}.}
 \ref{Cardy}{Cardy,J.L. \np{366}{1991}{403}.}
 \ref{ChandD}{Chang,P. and Dowker,J.S. \np{395}{1993}{407}.}
 \ref{DandC2}{Dowker,J.S. and Critchley,R. \prD{13}{1976}{224}.}
 \ref{Camporesi}{Camporesi,R. \prp{196}{1990}{1}.}
 \ref{BandM}{Brown,L.S. and Maclay,G.J. \pr{184}{1969}{1272}.}
 \ref{CandD}{Candelas,P. and Dowker,J.S. \prD{19}{1979}{2902}.}
 \ref{Unwin1}{Unwin,S.D. Thesis. University of Manchester. 1979.}
 \ref{Unwin2}{Unwin,S.D. \jpa{13}{1980}{313}.}
 \ref{DandB}{Dowker,J.S. and Banach,R. \jpa{11}{1978}{2255}.}
 \ref{Obhukov}{Obhukov,Yu.N. \pl{109B}{1982}{195}.}
 \ref{Kennedy}{Kennedy,G. \prD{23}{1981}{2884}.}
 \ref{CandT}{Copeland,E. and Toms,D.J. \np {255}{1985}{201}.}
  \ref{CandT2}{Copeland,E. and Toms,D.J. \cqg {3}{1986}{431}.}
 \ref{ELV}{Elizalde,E., Lygren, M. and Vassilevich,
 D.V. \jmp {37}{1996}{3105}.}
 \ref{Malurkar}{Malurkar,S.L. {\it J.Ind.Math.Soc} {\bf16} (1925/26) 130.}
 \ref{Glaisher}{Glaisher,J.W.L. {\it Messenger of Math.} {\bf18}
(1889) 1.} \ref{Anderson}{Anderson,A. \prD{37}{1988}{536}.}
 \ref{CandA}{Cappelli,A. and D'Appollonio,\pl{487B}{2000}{87}.}
 \ref{Wot}{Wotzasek,C. \jpa{23}{1990}{1627}.}
 \ref{RandT}{Ravndal,F. and Tollesen,D. \prD{40}{1989}{4191}.}
 \ref{SandT}{Santos,F.C. and Tort,A.C. \pl{482B}{2000}{323}.}
 \ref{FandO}{Fukushima,K. and Ohta,K. {\it Physica} {\bf A299} (2001) 455.}
 \ref{GandP}{Gibbons,G.W. and Perry,M. \prs{358}{1978}{467}.}
 \ref{Dow4}{Dowker,J.S..}
  \ref{Rad}{Rademacher,H. {\it Topics in analytic number theory,}
Springer-Verlag,  Berlin,1973.}
  \ref{Halphen}{Halphen,G.-H. {\it Trait\'e des Fonctions Elliptiques},
  Vol 1, Gauthier-Villars, Paris, 1886.}
  \ref{CandW}{Cahn,R.S. and Wolf,J.A. {\it Comm.Mat.Helv.} {\bf 51}
  (1976) 1.}
  \ref{Berndt}{Berndt,B.C. \rmjm{7}{1977}{147}.}
  \ref{Hurwitz}{Hurwitz,A. \ma{18}{1881}{528}.}
  \ref{Hurwitz2}{Hurwitz,A. {\it Mathematische Werke} Vol.I. Basel,
  Birkhauser, 1932.}
  \ref{Berndt2}{Berndt,B.C. \jram{303/304}{1978}{332}.}
  \ref{RandA}{Rao,M.B. and Ayyar,M.V. \jims{15}{1923/24}{150}.}
  \ref{Hardy}{Hardy,G.H. \jlms{3}{1928}{238}.}
  \ref{TandM}{Tannery,J. and Molk,J. {\it Fonctions Elliptiques},
   Gauthier-Villars, Paris, 1893--1902.}
  \ref{schwarz}{Schwarz,H.-A. {\it Formeln und
  Lehrs\"atzen zum Gebrauche..},Springer 1893.(The first edition was 1885.)
  The French translation by Henri Pad\'e is {\it Formules et Propositions
  pour L'Emploi...},Gauthier-Villars, Paris, 1894}
  \ref{Hancock}{Hancock,H. {\it Theory of elliptic functions}, Vol I.
   Wiley, New York 1910.}
  \ref{watson}{Watson,G.N. \jlms{3}{1928}{216}.}
  \ref{MandO}{Magnus,W. and Oberhettinger,F. {\it Formeln und S\"atze},
  Springer-Verlag, Berlin 1948.}
  \ref{Klein}{Klein,F. {\it Lectures on the Icosohedron}
  (Methuen, London. 1913).}
  \ref{AandL}{Appell,P. and Lacour,E. {\it Fonctions Elliptiques},
  Gauthier-Villars,
  Paris. 1897.}
  \ref{HandC}{Hurwitz,A. and Courant,C. {\it Allgemeine Funktionentheorie},
  Springer,
  Berlin. 1922.}
  \ref{WandW}{Whittaker,E.T. and Watson,G.N. {\it Modern analysis},
  Cambridge. 1927.}
  \ref{SandC}{Selberg,A. and Chowla,S. \jram{227}{1967}{86}. }
  \ref{zucker}{Zucker,I.J. {\it Math.Proc.Camb.Phil.Soc} {\bf 82 }(1977)
  111.}
  \ref{glasser}{Glasser,M.L. {\it Maths.of Comp.} {\bf 25} (1971) 533.}
  \ref{GandW}{Glasser, M.L. and Wood,V.E. {\it Maths of Comp.} {\bf 25}
  (1971)
  535.}
  \ref{greenhill}{Greenhill,A,G. {\it The Applications of Elliptic
  Functions}, MacMillan. London, 1892.}
  \ref{Weierstrass}{Weierstrass,K. {\it J.f.Mathematik (Crelle)}
{\bf 52} (1856) 346.}
  \ref{Weierstrass2}{Weierstrass,K. {\it Mathematische Werke} Vol.I,p.1,
  Mayer u. M\"uller, Berlin, 1894.}
  \ref{Fricke}{Fricke,R. {\it Die Elliptische Funktionen und Ihre Anwendungen},
    Teubner, Leipzig. 1915, 1922.}
  \ref{Konig}{K\"onigsberger,L. {\it Vorlesungen \"uber die Theorie der
 Elliptischen Funktionen},  \break Teubner, Leipzig, 1874.}
  \ref{Milne}{Milne,S.C. {\it The Ramanujan Journal} {\bf 6} (2002) 7-149.}
  \ref{Schlomilch}{Schl\"omilch,O. {\it Ber. Verh. K. Sachs. Gesell. Wiss.
  Leipzig}  {\bf 29} (1877) 101-105; {\it Compendium der h\"oheren
  Analysis}, Bd.II, 3rd Edn, Vieweg, Brunswick, 1878.}
  \ref{BandB}{Briot,C. and Bouquet,C. {\it Th\`eorie des Fonctions
  Elliptiques}, Gauthier-Villars, Paris, 1875.}
  \ref{Dumont}{Dumont,D. \aim {41}{1981}{1}.}
  \ref{Andre}{Andr\'e,D. {\it Ann.\'Ecole Normale Superior} {\bf 6} (1877)
  265;
  {\it J.Math.Pures et Appl.} {\bf 5} (1878) 31.}
  \ref{Raman}{Ramanujan,S. {\it Trans.Camb.Phil.Soc.} {\bf 22} (1916) 159;
 {\it Collected Papers}, Cambridge, 1927}
  \ref{Weber}{Weber,H.M. {\it Lehrbuch der Algebra} Bd.III, Vieweg,
  Brunswick 190  3.}
  \ref{Weber2}{Weber,H.M. {\it Elliptische Funktionen und algebraische
  Zahlen},
  Vieweg, Brunswick 1891.}
  \ref{ZandR}{Zucker,I.J. and Robertson,M.M.
  {\it Math.Proc.Camb.Phil.Soc} {\bf 95 }(1984) 5.}
  \ref{JandZ1}{Joyce,G.S. and Zucker,I.J.
  {\it Math.Proc.Camb.Phil.Soc} {\bf 109 }(1991) 257.}
  \ref{JandZ2}{Zucker,I.J. and Joyce.G.S.
  {\it Math.Proc.Camb.Phil.Soc} {\bf 131 }(2001) 309.}
  \ref{zucker2}{Zucker,I.J. {\it SIAM J.Math.Anal.} {\bf 10} (1979) 192,}
  \ref{BandZ}{Borwein,J.M. and Zucker,I.J. {\it IMA J.Math.Anal.} {\bf 12}
  (1992) 519.}
  \ref{Cox}{Cox,D.A. {\it Primes of the form $x^2+n\,y^2$}, Wiley,
  New York, 1989.}
  \ref{BandCh}{Berndt,B.C. and Chan,H.H. {\it Mathematika} {\bf42} (1995)
  278.}
  \ref{EandT}{Elizalde,R. and Tort.hep-th/}
  \ref{KandS}{Kiyek,K. and Schmidt,H. {\it Arch.Math.} {\bf 18} (1967) 438.}
  \ref{Oshima}{Oshima,K. \prD{46}{1992}{4765}.}
  \ref{greenhill2}{Greenhill,A.G. \plms{19} {1888} {301}.}
  \ref{Russell}{Russell,R. \plms{19} {1888} {91}.}
  \ref{BandB}{Borwein,J.M. and Borwein,P.B. {\it Pi and the AGM}, Wiley,
  New York, 1998.}
  \ref{Resnikoff}{Resnikoff,H.L. \tams{124}{1966}{334}.}
  \ref{vandp}{Van der Pol, B. {\it Indag.Math.} {\bf18} (1951) 261,272.}
  \ref{Rankin}{Rankin,R.A. {\it Modular forms} C.U.P. Cambridge}
  \ref{Rankin2}{Rankin,R.A. {\it Proc. Roy.Soc. Edin.} {\bf76 A} (1976) 107.}
  \ref{Skoruppa}{Skoruppa,N-P. {\it J.of Number Th.} {\bf43} (1993) 68 .}
  \ref{Down}{Dowker.J.S. {\it Nucl.Phys.}B (Proc.Suppl) ({\bf 104})(2002)153;
  also Dowker,J.S. hep-th/ 0007129.}
  \ref{Eichler}{Eichler,M. \mz {67}{1957}{267}.}
  \ref{Zagier}{Zagier,D. \invm{104}{1991}{449}.}
  \ref{Lang}{Lang,S. {\it Modular Forms}, Springer, Berlin, 1976.}
  \ref{Kosh}{Koshliakov,N.S. {\it Mess.of Math.} {\bf 58} (1928) 1.}
  \ref{BandH}{Bodendiek, R. and Halbritter,U. \amsh{38}{1972}{147}.}
  \ref{Smart}{Smart,L.R., \pgma{14}{1973}{1}.}
  \ref{Grosswald}{Grosswald,E. {\it Acta. Arith.} {\bf 21} (1972) 25.}
  \ref{Kata}{Katayama,K. {\it Acta Arith.} {\bf 22} (1973) 149.}
  \ref{Ogg}{Ogg,A. {\it Modular forms and Dirichlet series} (Benjamin,
  New York,
   1969).}
  \ref{Bol}{Bol,G. \amsh{16}{1949}{1}.}
  \ref{Epstein}{Epstein,P. \ma{56}{1903}{615}.}
  \ref{Petersson}{Petersson.}
  \ref{Serre}{Serre,J-P. {\it A Course in Arithmetic}, Springer,
  New York, 1973.}
  \ref{Schoenberg}{Schoenberg,B., {\it Elliptic Modular Functions},
  Springer, Berlin, 1974.}
  \ref{Apostol}{Apostol,T.M. \dmj {17}{1950}{147}.}
  \ref{Ogg2}{Ogg,A. {\it Lecture Notes in Math.} {\bf 320} (1973) 1.}
  \ref{Knopp}{Knopp,M.I. \dmj {45}{1978}{47}.}
  \ref{Knopp2}{Knopp,M.I. \invm {}{1994}{361}.}
  \ref{LandZ}{Lewis,J. and Zagier,D. \aom{153}{2001}{191}.}
  \ref{DandK1}{Dowker,J.S. and Kirsten,K. {\it Elliptic functions and
  temperature inversion symmetry on spheres} hep-th/.}
  \ref{HandK}{Husseini and Knopp.}
  \ref{Kober}{Kober,H. \mz{39}{1934-5}{609}.}
  \ref{HandL}{Hardy,G.H. and Littlewood, \am{41}{1917}{119}.}
  \ref{Watson}{Watson,G.N. \qjm{2}{1931}{300}.}
  \ref{SandC2}{Chowla,S. and Selberg,A. {\it Proc.Nat.Acad.} {\bf 35}
  (1949) 371.}
  \ref{Landau}{Landau, E. {\it Lehre von der Verteilung der Primzahlen},
  (Teubner, Leipzig, 1909).}
  \ref{Berndt4}{Berndt,B.C. \tams {146}{1969}{323}.}
  \ref{Berndt3}{Berndt,B.C. \tams {}{}{}.}
  \ref{Bochner}{Bochner,S. \aom{53}{1951}{332}.}
  \ref{Weil2}{Weil,A.\ma{168}{1967}{}.}
  \ref{CandN}{Chandrasekharan,K. and Narasimhan,R. \aom{74}{1961}{1}.}
  \ref{Rankin3}{Rankin,R.A. {} {} ().}
  \ref{Berndt6}{Berndt,B.C. {\it Trans.Edin.Math.Soc}.}
  \ref{Elizalde}{Elizalde,E. {\it Ten Physical Applications of Spectral
  Zeta Function Theory}, \break (Springer, Berlin, 1995).}
  \ref{Allen}{Allen,B., Folacci,A. and Gibbons,G.W. \pl{189}{1987}{304}.}
  \ref{Krazer}{Krazer}
  \ref{Elizalde3}{Elizalde,E. {\it J.Comp.and Appl. Math.} {\bf 118}
  (2000) 125.}
  \ref{Elizalde2}{Elizalde,E., Odintsov.S.D, Romeo, A. and Bytsenko,
  A.A and
  Zerbini,S.
  {\it Zeta function regularisation}, (World Scientific, Singapore,
  1994).}
  \ref{Eisenstein}{Eisenstein}
  \ref{Hecke}{Hecke,E. \ma{112}{1936}{664}.}
  \ref{Hecke2}{Hecke,E. \ma{112}{1918}{398}.}
  \ref{Terras}{Terras,A. {\it Harmonic analysis on Symmetric Spaces} (Springer,
  New York, 1985).}
  \ref{BandG}{Bateman,P.T. and Grosswald,E. {\it Acta Arith.} {\bf 9}
  (1964) 365.}
  \ref{Deuring}{Deuring,M. \aom{38}{1937}{585}.}
  \ref{Mordell}{Mordell,J. \prs{}{}{}.}
  \ref{GandZ}{Glasser,M.L. and Zucker, {}.}
  \ref{Landau2}{Landau,E. \jram{}{1903}{64}.}
  \ref{Kirsten1}{Kirsten,K. \jmp{35}{1994}{459}.}
  \ref{Sommer}{Sommer,J. {\it Vorlesungen \"uber Zahlentheorie}
  (1907,Teubner,Leipzig).
  French edition 1913 .}
  \ref{Reid}{Reid,L.W. {\it Theory of Algebraic Numbers},
  (1910,MacMillan,New York).}
  \ref{Milnor}{Milnor, J. {\it Is the Universe simply--connected?},
  IAS, Princeton, 1978.}
  \ref{Milnor2}{Milnor, J. \ajm{79}{1957}{623}.}
  \ref{Opechowski}{Opechowski,W. {\it Physica} {\bf 7} (1940) 552.}
  \ref{Bethe}{Bethe, H.A. \zfp{3}{1929}{133}.}
  \ref{LandL}{Landau, L.D. and Lishitz, E.M. {\it Quantum
  Mechanics} (Pergamon Press, London, 1958).}
  \ref{GPR}{Gibbons, G.W., Pope, C. and R\"omer, H., \np{157}{1979}{377}.}
  \ref{Jadhav}{Jadhav,S.P. PhD Thesis, University of Manchester 1990.}
  \ref{DandJ}{Dowker,J.S. and Jadhav, S. \prD{39}{1989}{1196}.}
  \ref{CandM}{Coxeter, H.S.M. and Moser, W.O.J. {\it Generators and
  relations of finite groups} (Springer. Berlin. 1957).}
  \ref{Coxeter2}{Coxeter, H.S.M. {\it Regular Complex Polytopes},
   (Cambridge University Press, \break Cambridge, 1975).}
  \ref{Coxeter}{Coxeter, H.S.M. {\it Regular Polytopes}.}
  \ref{Stiefel}{Stiefel, E., J.Research NBS {\bf 48} (1952) 424.}
  \ref{BandS}{Brink, D.M. and Satchler, G.R. {\it Angular momentum theory}.
  (Clarendon Press, Oxford. 1962.).}
  \ref{Rose}{Rose}
  \ref{Schwinger}{Schwinger, J. {\it On Angular Momentum}
  in {\it Quantum Theory of Angular Momentum} edited by
  Biedenharn,L.C. and van Dam, H. (Academic Press, N.Y. 1965).}
  \ref{Bromwich}{Bromwich, T.J.I'A. {\it Infinite Series},
  (Macmillan, London, 1926).}
  \ref{Ray}{Ray,D.B. \aim{4}{1970}{109}.}
  \ref{Ikeda}{Ikeda,A. {\it Kodai Math.J.} {\bf 18} (1995) 57.}
  \ref{Kennedy}{Kennedy,G. \prD{23}{1981}{2884}.}
  \ref{Ellis}{Ellis,G.F.R. {\it General Relativity} {\bf2} (1971) 7.}
  \ref{Dow8}{Dowker,J.S. \cqg{20}{2003}{L105}.}
  \ref{IandY}{Ikeda, A and Yamamoto, Y. \ojm {16}{1979}{447}.}
  \ref{BandI}{Bander,M. and Itzykson,C. \rmp{18}{1966}{2}.}
  \ref{Schulman}{Schulman, L.S. \pr{176}{1968}{1558}.}
  \ref{Bar1}{B\"ar,C. {\it Arch.d.Math.}{\bf 59} (1992) 65.}
  \ref{Bar2}{B\"ar,C. {\it Geom. and Func. Anal.} {\bf 6} (1996) 899.}
  \ref{Vilenkin}{Vilenkin, N.J. {\it Special functions},
  (Am.Math.Soc., Providence, 1968).}
  \ref{Talman}{Talman, J.D. {\it Special functions} (Benjamin,N.Y.,1968).}
  \ref{Miller}{Miller, W. {\it Symmetry groups and their applications}
  (Wiley, N.Y., 1972).}
  \ref{Dow3}{Dowker,J.S. \cmp{162}{1994}{633}.}
  \ref{Cheeger}{Cheeger, J. \jdg {18}{1983}{575}.}
  \ref{Cheeger2}{Cheeger, J. \aom {109}{1979}{259}.}
  \ref{Dow6}{Dowker,J.S. \jmp{30}{1989}{770}.}
  \ref{Dow20}{Dowker,J.S. \jmp{35}{1994}{6076}.}
  \ref{Dowjmp}{Dowker,J.S. \jmp{35}{1994}{4989}.}
  \ref{Dow21}{Dowker,J.S. {\it Heat kernels and polytopes} in {\it
   Heat Kernel Techniques and Quantum Gravity}, ed. by S.A.Fulling,
   Discourses in Mathematics and its Applications, No.4, Dept.
   Maths., Texas A\&M University, College Station, Texas, 1995.}
  \ref{Dow9}{Dowker,J.S. \jmp{42}{2001}{1501}.}
  \ref{Dow7}{Dowker,J.S. \jpa{25}{1992}{2641}.}
  \ref{Warner}{Warner.N.P. \prs{383}{1982}{379}.}
  \ref{Wolf}{Wolf, J.A. {\it Spaces of constant curvature},
  (McGraw--Hill,N.Y., 1967).}
  \ref{Meyer}{Meyer,B. \cjm{6}{1954}{135}.}
  \ref{BandB}{B\'erard,P. and Besson,G. {\it Ann. Inst. Four.} {\bf 30}
  (1980) 237.}
  \ref{PandM}{P\'{o}lya,G. and Meyer,B. \cras{228}{1948}{28}.}
  \ref{Springer}{Springer, T.A. Lecture Notes in Math. vol 585 (Springer,
  Berlin,1977).}
  \ref{SeandT}{Threlfall, H. and Seifert, W. \ma{104}{1930}{1}.}
  \ref{Hopf}{Hopf,H. \ma{95}{1925}{313}. }
  \ref{Dow}{Dowker,J.S. \jpa{5}{1972}{936}.}
  \ref{LLL}{Lehoucq,R., Lachi\'eze-Rey,M. and Luminet, J.--P. {\it
  Astron.Astrophys.} {\bf 313} (1996) 339.}
  \ref{LaandL}{Lachi\'eze-Rey,M. and Luminet, J.--P.
  \prp{254}{1995}{135}.}
  \ref{Schwarzschild}{Schwarzschild, K., {\it Vierteljahrschrift der
  Ast.Ges.} {\bf 35} (1900) 337.}
  \ref{Starkman}{Starkman,G.D. \cqg{15}{1998}{2529}.}
  \ref{LWUGL}{Lehoucq,R., Weeks,J.R., Uzan,J.P., Gausman, E. and
  Luminet, J.--P. \cqg{19}{2002}{4683}.}
  \ref{Dow10}{Dowker,J.S. \prD{28}{1983}{3013}.}
  \ref{BandD}{Banach, R. and Dowker, J.S. \jpa{12}{1979}{2527}.}
  \ref{Jadhav2}{Jadhav,S. \prD{43}{1991}{2656}.}
  \ref{Gilkey}{Gilkey,P.B. {\it Invariance theory,the heat equation and
  the Atiyah--Singer Index theorem} (CRC Press, Boca Raton, 1994).}
  \ref{BandY}{Berndt,B.C. and Yeap,B.P. {\it Adv. Appl. Math.}
  {\bf29} (2002) 358.}
  \ref{HandR}{Hanson,A.J. and R\"omer,H. \pl{80B}{1978}{58}.}
  \ref{Hill}{Hill,M.J.M. {\it Trans.Camb.Phil.Soc.} {\bf 13} (1883) 36.}
  \ref{Cayley}{Cayley,A. {\it Quart.Math.J.} {\bf 7} (1866) 304.}
  \ref{Seade}{Seade,J.A. {\it Anal.Inst.Mat.Univ.Nac.Aut\'on
  M\'exico} {\bf 21} (1981) 129.}
  \ref{CM}{Cisneros--Molina,J.L. {\it Geom.Dedicata} {\bf84} (2001)
  \ref{Goette1}{Goette,S. \jram {526} {2000} 181.}
  207.}
  \ref{NandO}{Nash,C. and O'Connor,D--J, \jmp {36}{1995}{1462}.}
  \ref{Dows}{Dowker,J.S. \aop{71}{1972}{577}; Dowker,J.S. and Pettengill,D.F.
  \jpa{7}{1974}{1527}; J.S.Dowker in {\it Quantum Gravity}, edited by
  S. C. Christensen (Hilger,Bristol,1984)}
  \ref{Jadhav2}{Jadhav,S.P. \prD{43}{1991}{2656}.}
  \ref{Dow11}{Dowker,J.S. \cqg{21}{2004}4247.}
  \ref{Dow12}{Dowker,J.S. \cqg{21}{2004}4977.}
  \ref{Dow13}{Dowker,J.S. \jpa{38}{2005}1049.}
  \ref{Zagier}{Zagier,D. \ma{202}{1973}{149}}
  \ref{RandG}{Rademacher, H. and Grosswald,E. {\it Dedekind Sums},
  (Carus, MAA, 1972).}
  \ref{Berndt7}{Berndt,B, \aim{23}{1977}{285}.}
  \ref{HKMM}{Harvey,J.A., Kutasov,D., Martinec,E.J. and Moore,G.
  {\it Localised Tachyons and RG Flows}, hep-th/0111154.}
  \ref{Beck}{Beck,M., {\it Dedekind Cotangent Sums}, {\it Acta Arithmetica}
  {\bf 109} (2003) 109-139 ; math.NT/0112077.}
  \ref{McInnes}{McInnes,B. {\it APS instability and the topology of the brane
  world}, hep-th/0401035.}
  \ref{BHS}{Brevik,I, Herikstad,R. and Skriudalen,S. {\it Entropy Bound for the
  TM Electromagnetic Field in the Half Einstein Universe}; hep-th/0508123.}
  \ref{BandO}{Brevik,I. and Owe,C.  \prD{55}{4689}{1997}.}
  \ref{Kenn}{Kennedy,G. Thesis. University of Manchester 1978.}
  \ref{KandU}{Kennedy,G. and Unwin S. \jpa{12}{L253}{1980}.}
  \ref{BandO1}{Bayin,S.S.and Ozcan,M.
  \prD{48}{2806}{1993}; \prD{49}{5313}{1994}.}
  \ref{Chang}{Chang, P., {\it Quantum Field Theory on Regular Polytopes}.
   Thesis. University of Manchester, 1993.}
  \ref{Barnesa}{Barnes,E.W. {\it Trans. Camb. Phil. Soc.} {\bf 19} (1903) 374.}
  \ref{Barnesb}{Barnes,E.W. {\it Trans. Camb. Phil. Soc.}
  {\bf 19} (1903) 426.}
  \ref{Stanley1}{Stanley,R.P. \joa {49Hilf}{1977}{134}.}
  \ref{Stanley}{Stanley,R.P. \bams {1}{1979}{475}.}
  \ref{Hurley}{Hurley,A.C. \pcps {47}{1951}{51}.}
  \ref{IandK}{Iwasaki,I. and Katase,K. {\it Proc.Japan Acad. Ser} {\bf A55}
  (1979) 141.}
  \ref{IandT}{Ikeda,A. and Taniguchi,Y. {\it Osaka J. Math.} {\bf 15} (1978)
  515.}
  \ref{GandM}{Gallot,S. and Meyer,D. \jmpa{54}{1975}{259}.}
  \ref{Flatto}{Flatto,L. {\it Enseign. Math.} {\bf 24} (1978) 237.}
  \ref{OandT}{Orlik,P and Terao,H. {\it Arrangements of Hyperplanes},
  Grundlehren der Math. Wiss. {\bf 300}, (Springer--Verlag, 1992).}
  \ref{Shepler}{Shepler,A.V. \joa{220}{1999}{314}.}
  \ref{SandT}{Solomon,L. and Terao,H. \cmh {73}{1998}{237}.}
  \ref{Vass}{Vassilevich, D.V. \plb {348}{1995}39.}
  \ref{Vass2}{Vassilevich, D.V. \jmp {36}{1995}3174.}
  \ref{CandH}{Camporesi,R. and Higuchi,A. {\it J.Geom. and Physics}
  {\bf 15} (1994) 57.}
  \ref{Solomon2}{Solomon,L. \tams{113}{1964}{274}.}
  \ref{Solomon}{Solomon,L. {\it Nagoya Math. J.} {\bf 22} (1963) 57.}
  \ref{Obukhov}{Obukhov,Yu.N. \pl{109B}{1982}{195}.}
  \ref{BGH}{Bernasconi,F., Graf,G.M. and Hasler,D. {\it The heat kernel
  expansion for the electromagnetic field in a cavity}; math-ph/0302035.}
  \ref{Baltes}{Baltes,H.P. \prA {6}{1972}{2252}.}
  \ref{BaandH}{Baltes.H.P and Hilf,E.R. {\it Spectra of Finite Systems}
  (Bibliographisches Institut, Mannheim, 1976).}
  \ref{Ray}{Ray,D.B. \aim{4}{1970}{109}.}
  \ref{Hirzebruch}{Hirzebruch,F. {\it Topological methods in algebraic
  geometry} (Springer-- Verlag,\break  Berlin, 1978). }
  \ref{BBG}{Bla\v{z}i\'c,N., Bokan,N. and Gilkey, P.B. {\it Ind.J.Pure and
  Appl.Math.} {\bf 23} (1992) 103.}
  \ref{WandWi}{Weck,N. and Witsch,K.J. {\it Math.Meth.Appl.Sci.} {\bf 17}
  (1994) 1017.}
  \ref{Norlund}{N\"orlund,N.E. \am{43}{1922}{121}.}
   \ref{Norlund1}{N\"orlund,N.E. {\it Differenzenrechnung} (Springer--Verlag, 1924, Berlin.)}
  \ref{Duff}{Duff,G.F.D. \aom{56}{1952}{115}.}
  \ref{DandS}{Duff,G.F.D. and Spencer,D.C. \aom{45}{1951}{128}.}
  \ref{BGM}{Berger, M., Gauduchon, P. and Mazet, E. {\it Lect.Notes.Math.}
  {\bf 194} (1971) 1. }
  \ref{Patodi}{Patodi,V.K. \jdg{5}{1971}{233}.}
  \ref{GandS}{G\"unther,P. and Schimming,R. \jdg{12}{1977}{599}.}
  \ref{MandS}{McKean,H.P. and Singer,I.M. \jdg{1}{1967}{43}.}
  \ref{Conner}{Conner,P.E. {\it Mem.Am.Math.Soc.} {\bf 20} (1956).}
  \ref{Gilkey2}{Gilkey,P.B. \aim {15}{1975}{334}.}
  \ref{MandP}{Moss,I.G. and Poletti,S.J. \plb{333}{1994}{326}.}
  \ref{BKD}{Bordag,M., Kirsten,K. and Dowker,J.S. \cmp{182}{1996}{371}.}
  \ref{RandO}{Rubin,M.A. and Ordonez,C. \jmp{25}{1984}{2888}.}
  \ref{BaandD}{Balian,R. and Duplantier,B. \aop {112}{1978}{165}.}
  \ref{Kennedy2}{Kennedy,G. \aop{138}{1982}{353}.}
  \ref{DandKi2}{Dowker,J.S. and Kirsten, K. {\it Analysis and Appl.}
 {\bf 3} (2005) 45.}
  \ref{Dow40}{Dowker,J.S. \cqg{23}{2006}{1}.}
  \ref{BandHe}{Br\"uning,J. and Heintze,E. {\it Duke Math.J.} {\bf 51} (1984)
   959.}
  \ref{Dowl}{Dowker,J.S. {\it Functional determinants on M\"obius corners};
    Proceedings, `Quantum field theory under
    the influence of external conditions', 111-121,Leipzig 1995.}
  \ref{Dowqg}{Dowker,J.S. in {\it Quantum Gravity}, edited by
  S. C. Christensen (Hilger, Bristol, 1984).}
  \ref{Dowit}{Dowker,J.S. \jpa{11}{1978}{347}.}
  \ref{Kane}{Kane,R. {\it Reflection Groups and Invariant Theory} (Springer,
  New York, 2001).}
  \ref{Sturmfels}{Sturmfels,B. {\it Algorithms in Invariant Theory}
  (Springer, Vienna, 1993).}
  \ref{Bourbaki}{Bourbaki,N. {\it Groupes et Alg\`ebres de Lie}  Chap.III, IV
  (Hermann, Paris, 1968).}
  \ref{SandTy}{Schwarz,A.S. and Tyupkin, Yu.S. \np{242}{1984}{436}.}
  \ref{Reuter}{Reuter,M. \prD{37}{1988}{1456}.}
  \ref{EGH}{Eguchi,T. Gilkey,P.B. and Hanson,A.J. \prp{66}{1980}{213}.}
  \ref{DandCh}{Dowker,J.S. and Chang,Peter, \prD{46}{1992}{3458}.}
  \ref{APS}{Atiyah M., Patodi and Singer,I.\mpcps{77}{1975}{43}.}
  \ref{Donnelly}{Donnelly.H. {\it Indiana U. Math.J.} {\bf 27} (1978) 889.}
  \ref{Katase}{Katase,K. {\it Proc.Jap.Acad.} {\bf 57} (1981) 233.}
  \ref{Gilkey3}{Gilkey,P.B.\invm{76}{1984}{309}.}
  \ref{Degeratu}{Degeratu.A. {\it Eta--Invariants and Molien Series for
  Unimodular Groups}, Thesis MIT, 2001.}
  \ref{Seeley}{Seeley,R. \ijmp {A\bf18}{2003}{2197}.}
  \ref{Seeley2}{Seeley,R. .}
  \ref{melrose}{Melrose}
  \ref{DandW}{Douglas,R.G. and Wojciekowski,K.P. \cmp{142}{1991}{139}.}
  \ref{Dai}{Dai,X. \tams{354}{2001}{107}.}
\end{putreferences}

\bye